\documentclass[11pt]{article}

\def\mathbb{\Bbb}
\usepackage{amsfonts,latexsym,amstext}
\usepackage{amsmath}
\usepackage{amssymb}

\usepackage{comment}
\usepackage{amsthm}
\usepackage[latin1]{inputenc}
\allowdisplaybreaks
\usepackage{color}

\theoremstyle{plain}
\newtheorem{theorem}{Theorem}[section]
\newtheorem{lemma}[theorem]{Lemma}
\newtheorem{proposition}[theorem]{Proposition}
\newtheorem{corollary}[theorem]{Corollary}

\newtheorem{definition}{Definition}[section]

\theoremstyle{hypothesis}

\theoremstyle{remark}
\newtheorem{remark}[theorem]{Remark}

\makeatletter
\@addtoreset{equation}{section}

\makeatother

\def\qed{{\hfill\hbox{\enspace${ \square}$}} \smallskip}
\def\sqr#1#2{{\vcenter{\vbox{\hrule height .#2pt \hbox{\vrule
 width .#2pt height#1pt \kern#1pt \vrule
width .#2pt} \hrule height .#2pt}}}}
\def\square{\mathchoice\sqr54\sqr54\sqr{4.1}3\sqr{3.5}3}

\def\ds{\begin{displaystyle}}
\def\eds{\end{displaystyle}}

\def\<{\langle }
\def\>{\rangle }

\def\R{\mathbb R}
\def\N{\mathbb N}

\def\E{\mathbb E}
\def\P{\mathbb P}
\def\Q{\mathbb Q}
\def\F{\mathbb F}

\def\cala{{\cal A}}
\def\calb{{\cal B}}

\def\cale{{\cal E}}
\def\calf{{\cal F}}
\def\calg{{\cal G}}

\def\calp{{\cal P}}

\newcommand{\sper}[1]{\mathbb{E} \left[ #1 \right]}                               
\newcommand{\spernuxa}[1]{\mathbb{E}^{x,a}_{\nu} \left[ #1 \right]}

\newcommand{\spernuxaepsilon}[1]{\mathbb{E}^{x,a}_{\nu^{\varepsilon}} \left[ #1 \right]}                          
\newcommand{\sperxa}[1]{\mathbb{E}^{x,a} \left[ #1 \right]}                               

\DeclareMathAlphabet{\mathonebb}{U}{bbold}{m}{n}                           %
\newcommand{\one}{\ensuremath{\mathonebb{1}}}                               
\newcommand{\essinf}{\operatornamewithlimits{ess\,inf}}   
\newcommand{\limit}{\operatornamewithlimits{\longrightarrow}}   

\title{Optimal control of Piecewise Deterministic Markov Processes:\\ a   BSDE representation of the value function}
\date{}
\author{Elena BANDINI\thanks{Politecnico di Milano, Dipartimento di Matematica, via Bonardi 9, 20133 Milano, Italy; ENSTA ParisTech, Unit\'e de Math\'ematiques appliqu\'ees, 828, boulevard des Mar\'echaux, F-91120 Palaiseau, France; e-mail: \texttt{elena.bandini@polimi.it}}
}

\begin{document}

\maketitle

\begin{abstract}
We consider an infinite horizon  discounted optimal control problem for  piecewise deterministic Markov processes, where  a piecewise open-loop control acts continuously on the jump dynamics and on the deterministic flow.
For this class of control problems, the value function can in general be characterized as the unique viscosity solution to the corresponding Hamilton-Jacobi-Bellman equation.
We prove that the value function can be represented by means of a  backward stochastic differential equation (BSDE) on infinite horizon, driven by a random measure and with a sign constraint on its martingale part, for which we give existence and uniqueness results. 
This probabilistic representation is known as nonlinear Feynman-Kac formula. 
Finally we show that the constrained BSDE is related to  an auxiliary dominated control problem, whose value function   coincides with  the value function of the original non-dominated control problem.
%
%

\end{abstract}


{\small\textbf{MSC 2010:} 93E20, 60H10, 60J25.}

\section{Introduction}\label{Sec_introduction}
The aim of the present paper is to prove that the value function in an infinite horizon  optimal control problem for  piecewise deterministic Markov processes (PDMPs) can be represented by means of an appropriate Backward Stochastic Differential Equation (BSDE).
Piecewise deterministic Markov processes,  introduced in  
\cite{Da}, evolve
through random jumps at random times,
while the behavior between jumps
is described  by a deterministic flow. 
We  consider optimal control problems of PDMPs where  the control acts continuously on the jump dynamics and
on the deterministic flow
as well.



Let us start by  describing our setting in an informal way.
Let $E$ be a Borel space and
$\mathcal E$ the corresponding $\sigma$-algebra.
A PDMP on 
$(E,\mathcal E)$ can be described by means of three local characteristics, namely a continuous flow $\phi(t,x)$,  a jump rate $\lambda(x)$, and a transition measure $Q(x,dy)$,   according to which the location of the process at the jump time is chosen. 
The  PDMP dynamic can be described as follows: starting from some initial point $x \in E$, the motion of the process follows the flow $\phi(t,x)$ until a random jump $T_1$, verifying
$$
\P(T_1 > s)= \exp\left(- \int_0^s \lambda(\phi(r,x))\,dr\right),\quad s \geq 0.
$$
At time $T_1$ the process jumps to a new point $X_{T_1}$ selected with probability $Q(x,dy)$ (conditionally to $T_1$),  and the motion restarts from this new point as before.

 At this point we introduce a measurable space  $(A, \mathcal A)$ , which will denote the space  of control actions.
 A controlled PDMP is obtained starting from  a jump rate $\lambda(x,a)$ and a transition measure $Q(x,a,dy))$, depending on  
   an additional control parameter  $a\in A$, and  a continuous flow $\phi^{\beta}(t,x)$, depending  on the choice of a measurable function $\beta(t)$ taking values on $(A, \mathcal A)$. A natural way to control a PDMP is to chose a  control strategy among the set of \emph{piecewise open-loop policies}, i.e., measurable functions that depend only on the last jump time and post jump position. We can mention 
\cite{Almudevar}, 
\cite{BauerleRieder}, 
 \cite{CoDu}, 
 \cite{Da}, \cite{DavisArticle}, 
 \cite{Dem}, as a sample of works that use this kind of approach.
Roughly speaking, at each jump time $T_n$, we choose an open loop control $\alpha_n$ depending on the initial
condition $X_{T_n}$ to be used until the next jump time.
A control $\alpha$ in the class of admissible control laws $\mathcal A_{ad}$ has the explicit form
\begin{equation}\label{open_loop_controlsINTRO}
	\alpha_t= \sum_{n=1}^{\infty}\alpha_n(t-T_n,X_{T_n})\,\one_{[T_n,\,T_{n+1})}(t),
\end{equation}
and the controlled process $X$ is
$$
X_t = \phi^{\alpha_n}(t- T_n, E_n),\quad t \in [T_n,\, T_{n+ 1}).
$$
We denote by $\P^{x}_{\alpha}$ the probability measure such that,
for every $n >1$,
the conditional survivor function of
the   jump time $T_{n+1}$
and the distribution of the
post jump position $X_{T_{n+1}}$ are
\begin{align*}
	&\P^x_{\alpha}(T_{n+1} > s\,|\,\mathcal F_{T_n})=
	\exp\left(-\int_{T_n}^s \lambda (\phi^{\alpha_n}(r-T_n,X_{T_n}),\alpha_n(r-T_n,X_{T_n}))\,dr\right),\\
 &\P^x_{\alpha}(X_{T_{n+1}} \in B 
 |\,\mathcal F_{T_n},\,T_{n+1})=
 Q( \phi^{\alpha_n}(T_{n+1}-T_n,X_{T_n}),\alpha_n(T_{n+1}-T_n,X_{T_n}), B),
\end{align*}
on $\{T_n < \infty\}$.

In the classic infinite horizon control problem one wants to minimize over all control laws $\alpha$ a functional cost of the form
\begin{equation}\label{Sec:PDP_functional_costINTRO}
	J(x,\alpha) = \E^{x}_{\alpha}\left[\int_{0}^{\infty} e^{-\delta\,s}\, f(X_{s},\alpha_s)\,ds\right]
\end{equation}
where $\E^{x}_{\alpha}$ denotes the expectation under $\P^{x}_{\alpha}$, $f$ is  a given real function on  $E \times A$  representing the running cost, and  $\delta \in (0,\,\infty)$ is a  discounting factor.
The value function of the control problem is defined in the usual way:
\begin{equation}\label{Sec:PDP_value_functionINTRO}
	V(x) = \inf_{\alpha \in \mathcal{A}_{ad}}J(x,\alpha),\quad x \in E.
\end{equation}

Let now $E$ be an open subset of $\R^d$, and $h(x,a)$ be a bounded Lipschitz continuous function such that $\phi^{\alpha}(t,x)$ is the unique solution of the
ordinary differential equation
$$
\dot x(t) = h(x(t), \alpha(t)), \quad x(0)=x \in E.
$$
We will  assume that $\lambda$ and  $f$ are bounded functions, uniformly continuous,  and $Q$ is a Feller stochastic kernel.
In this case, $V$ is known to be the unique viscosity solution on $[0,\,\infty) \times E$ of the Hamilton-Jacobi-Bellman (HJB) equation
\begin{equation}\label{Sec:PDP_HJBINTRO}
	\delta v(x)= \sup_{a \in A} \left(h(x,a)\cdot \nabla v(x) + \lambda(x,a) \int_E (v(y)-v(x))\,Q(x,a,dy)\right)\quad x \in E.
\end{equation}
The  characterization
of the optimal value function as the viscosity solution of the corresponding integro-differential HJB equation is an important approach to tackle the  optimal control problem of PDMPs, and can  be found for instance in 
 \cite{Da-Fa}, 
 \cite{DemYe1}. 
Alternatively,
the control problem can be reformulated as a discrete-stage Markov decision model, where the stages are the jumps times of the process and the decision at each stage is the control function that solves a deterministic optimal control problem.
The reduction of  the optimal control problem  to a discrete-time Markov decision process is exploited
for instance in  \cite{Almudevar}, \cite{BauerleRieder}, \cite{CoDu}, \cite{Da}, \cite{DavisArticle}.

In the present paper our aim is to represent the value function $V(x)$ by means of an appropriate BSDE. We are interested in the general case when $\{\P^{x}_{\alpha}\}_{\alpha}$ is a non-dominated model, which, roughly speaking, reflects the fully non-linear character of the HJB equation. This basic difficulty has prevented the effective use of BSDE techniques in the context of optimal control of PDMPs until now. In fact, we believe that the present paper is the first one where this difficulty has been coped with and this connection has been established. It is our hope that the great development that BSDE theory has now gained will produce new results in the optimization theory of PDMPs.

In the context of diffusions, probabilistic formulae for the value function for nondominated models have been discovered only in the recent year.
In this sense, a fundamental role is played by   \cite{KhPh}, where  a new class of BSDEs with nonpositive jumps is introduced in order to provide a probabilistic formula, known as nonlinear Feynman-Kac formula, for
fully nonlinear integro-partial differential equations, associated to the classical optimal control for diffusions.
This approach was later applied to many cases within  optimal switching and impulse control problems, see 
\cite{ElKh10}, \cite{ElKh14}, \cite{ElKh14a}, 
\cite{KhMaPhZh}, and developed with extensions and applications, see for instance 
\cite{Co-Chou}, 
\cite{FuPh}. In all the above mentioned cases the controlled processes are   diffusions  constructed as solutions to stochastic differential equations of Ito type driven by a Brownian motion.

We wish  to extend   to the   PDMPs framework the theory developed in the context of  optimal control for diffusions.
The fundamental idea behind the derivation of the Feynman-Kac representation, borrowed from \cite{KhPh},  concerns
the so-called \emph{randomization of the control}, that  we are going to describe below in our framework.
A  first step in  the generalization of this method to the non-diffusive processes context was done in \cite{BandiniFuhrman}, where   a probabilistic representation for the value function associated to an optimal control problem for pure jump Markov processes was provided. As  in the pure jump case, also in the PDMPs framework the correct formulation of the randomization method  requires some efforts, and and can not be modelled on the diffusive case, since the controlled processes are not defined as solutions to stochastic differential equations.
Moreover, the
treatment  in the PDMPs context is more involved and requires different techniques.
For instance, the presence  in the  PDMP's  dynamics  of the controlled flow  
leads to a  differential operator in the HJB equation, which
suggests to  use  the viscosity solution theory;
in addition, since we consider  optimal control problems on infinite horizon, we will need to  deal with BSDEs on infinite horizon as well.

Finally, we notice that we consider PDMPs with state space $E$ with no boundary. This restriction is due to the fact that the presence of the boundary induces technical difficulties on the study of the associated BSDE, which would be driven by a non-quasi left continuous random measure, see Remark \ref{Sec:PDP_R_non_quasi_left}.   
For such general BSDEs, the existence and uniqueness results were at disposal  only in particular frameworks, see e.g. \cite{CohenElliott} for the deterministic case, and counter-examples were provided in the general case, see Section 4.3 in \cite{CFJ}.
Only recently  this problem was faced and solved in a general context in   \cite{BandiniBSDE}:
this fact opens to the possibility of further extensions of the BSDEs approach,
and will be the object of a second work by the author.


Let us now informally describe the randomization method in the PDMPs framework.
The  first step, for any starting point $x \in E$, consists in   replacing the state trajectory and the associated control process $(X_s,\alpha_s)$ by an (uncontrolled) PDMP $(X_s,I_s)$,  in such a way that $I$ is a Poisson process with values in the space of control actions $A$, with an intensity $\lambda_0(db)$ which is arbitrary but finite and with full support, and $X$ is suitably defined.
In particular, the PDMP $(X,I)$ is  constructed in a different probability space by means of a new triplet of local characteristics and takes values on the enlarged space $E \times A$. Let us denote by $\P^{x,a}$ the corresponding law, where $(x,a)$ is the starting point in $E \times A$.
Then we formulate an auxiliary optimal control problem where we control the intensity of the process $I$: for any predictable, bounded and positive random field $\nu_t(b)$, by means of a theorem of Girsanov type, we construct a probability measure $\P^{x,a}_\nu$ under which the compensator of $I$ is the random measure $\nu_t(db)\,\lambda_0(db)\,dt$ (under $\P^{x,a}_\nu$ the law of $X$ is also changed) and we minimize the functional
\begin{equation}\label{Sec:PDP_functional_dualINTRO}
	J(x,a, \nu) = \E^{x,a}_{\nu}\left[\int_{0}^{\infty} e^{-\delta\,s}\, f(X_{s},I_s)\,ds\right].
\end{equation}
over all possible choices of $\nu$. This will be called the \emph{dual} control problem. Notice that, by the Girsanov theorem,  the family $\{\P^{x,a}_{\nu}\}_{\nu}$ is a dominated model. One of our main results states that the value function of the dual control problem, denoted as $V^\ast(x,a)$, can be represented by means of a well-posed constrained BSDE. The latter  is an  equation over an infinite horizon of the form
\begin{align}
	& Y^{x,a}_{s} = Y^{x,a}_T - \delta \,\int_{s}^{T}Y^{x,a}_r\,dr + \int_{s}^{T}f(X_r,I_r)\,dr -( K^{x,a}_T - K^{x,a}_s) \label{Sec:PDP_BSDEINTRO}\\
	&- \int_{s}^{T}\int_{A}Z^{x,a}_r(X_r,\,b)\, \lambda_0(db)\,dr - \int_{s}^{T}\int_{E \times A}Z^{x,a}_r(y,\,b)\, q(dr\,dy\, db), \quad  0\leqslant s \leqslant T<\infty,\nonumber
\end{align}
with unknown triplet $(Y^{x,a},Z^{x,a},K^{x,a})$
where $q$ is the compensated random measure associated to $(X,I)$, $K^{x,a}$ is a predictable increasing c\`adl\`ag process, $Z^{x,a}$ is a predictable random field, where we additionally add the sign constraint
\begin{equation}\label{Sec:PDP_BSDE_constraintINTRO}
	Z^{x,a}_s(X_{s-},b)\geqslant 0.
\end{equation}
The reference filtration is now the  canonical one associated to the pair $(X,I)$.
We prove  that this equation has a unique maximal solution, in an appropriate sense, and that the value of the process $Y^{x,a}$ at the initial time represents the dual value function:
\begin{equation}\label{Sec:PDP_randvalue_id_INTRO}
	Y^{x,a}_0= V^{\ast}(x,a).
\end{equation}
Our main purpose is to show that the maximal solution to \eqref{Sec:PDP_BSDEINTRO}-\eqref{Sec:PDP_BSDE_constraintINTRO} at the initial time  also  provides a Feynman-Kac representation to the value function \eqref{Sec:PDP_value_functionINTRO} of our original  optimal control problem for PDMPs.
To this end, we introduce the deterministic real function on $E \times A$
\begin{equation}\label{Sec:PDP_vINTRO}
	v(x,a):=Y_0^{x,a},
\end{equation}
and we prove that $v$ is a viscosity solution to \eqref{Sec:PDP_HJBINTRO}.  By the uniqueness of the solution to the HJB equation \eqref{Sec:PDP_HJBINTRO}  we conclude that the value of the process $Y$ at the initial time represents both the original and the dual value function:
\begin{equation}\label{Sec:PDP_rapprINTRO}
	Y_0^{x,a}=V^\ast(x,a)=V(x).
\end{equation}
Identity \eqref{Sec:PDP_rapprINTRO} is the desired BSDE representation of the value function for the original control problem and a Feynman-Kac formula for the general HJB equation \eqref{Sec:PDP_HJBINTRO}.

Formula \eqref{Sec:PDP_rapprINTRO}
can be used to design  
algorithms 
based on  the numerical approximation of the solution to the constrained BSDE \eqref{Sec:PDP_BSDEINTRO}-\eqref{Sec:PDP_BSDE_constraintINTRO}, and therefore to get
probabilistic numerical approximations  for the value function of the addressed   optimal control problem.
In the recent years there has been much interest in this  problem, and 
numerical schemes for  constrained BSDEs have been  proposed and analyzed in the diffusive framework, see \cite{KhLaPh},
\cite{KhLaPha}. 
We hope that our results may be used to get similar techniques in the PDMPs context as well.

The paper is organized as follows. Section \ref{Sec:PDP_Sec_preliminaries} is dedicated to  define a setting where the optimal control \eqref{Sec:PDP_value_functionINTRO} is solved by means of the corresponding HJB equation \eqref{Sec:PDP_HJBINTRO}.
We start by recalling the construction of a PDMP given its local characteristics. In order to apply techniques based on BSDEs driven by general random measures, we work in a canonical setting and we use a specific filtration.
The construction is based on the well-posedness of the martingale problem for multivariate marked point processes  studied in Jacod \cite{J}, and is the object of Section \ref{Sec:PDP_Sec_constr_PDMP}. This general procedure is then applied  in Section \ref{Sec:PDP_Sec_control_problem} to formulate in a precise way the optimal control problem we are interested in.
At the end of Section \ref{Sec:PDP_Sec_control_problem} we recall a classical result on existence and uniqueness of the viscosity solution to the HJB equation \eqref{Sec:PDP_HJBINTRO}, and its identification with the value function $V$,  provided by  Davis and Farid \cite{Da-Fa}.

In Section \ref{Sec:PDP_Section_dual_control} we start to develop the control randomization method. Given suitable  local characteristics, we introduce an auxiliary process $(X,I)$  on $E \times A$ by relying on  the construction in Section \ref{Sec:PDP_Sec_constr_PDMP}, and we formulate a dual optimal control problem for it under suitable conditions. The formulation of the randomized process is very different from the diffusive framework, since our data are the local characteristics of the process rather than the coefficients of some stochastic differential equations solved by it.
In particular, we need to choose a specific probability space under which the component $I$ (independent to $X$) is a Poisson process.

In Section \ref{Sec:PDP_Sec_ConstrainedBSDE} we introduce the constrained BSDE \eqref{Sec:PDP_BSDEINTRO}-\eqref{Sec:PDP_BSDE_constraintINTRO} over infinite horizon. By a penalization approach, we prove that under suitable assumptions the above mentioned equation admits a unique maximal solution $(Y,Z,K)$ in a certain class of processes. Moreover, the component $Y$ at the initial time coincides with the value function $V^\ast$ of the dual optimal control problem. This is the first of our main results, and  is the object of Theorem \ref{Sec:PDP_Thm_ex_uniq_maximal_BSDE}.

Finally, in Section \ref{Sec:PDP_Section_nonlinear_IPDE} we prove that the initial value of the maximal solution $Y^{x,a}$ to \eqref{Sec:PDP_BSDEINTRO}-\eqref{Sec:PDP_BSDE_constraintINTRO} provides a viscosity solution to \eqref{Sec:PDP_HJBINTRO}.
 This  is the second main result of the paper, which is stated in Theorem \ref{Sec:PDP_THm_Feynman_Kac_HJB}. 
 As a consequence, by means of a comparison theorem for sub and supersolutions to first-order integro-partial differential equations, 
 we get the desired non-linear Feynman-Kac formula, as well as  the equality between the value functions of the primal and the dual control problems, see Corollary \ref{Sec:PDP_C_final}.   The proof of Theorem \ref{Sec:PDP_THm_Feynman_Kac_HJB} is based on arguments from the viscosity theory, and combines BSDEs techniques with control-theoretic arguments.
A relevant task  is to derive the key property that the function $v$ in \eqref{Sec:PDP_vINTRO} does not  depend on $a$, as consequence of the $A$-nonnegative constrained jumps.
Recalling the identification in Theorem \ref{Sec:PDP_Thm_ex_uniq_maximal_BSDE}, we are able  to give a  direct proof of the non-dependence of $v$ on $a$ by means of control-theoretic techniques, see Proposition \ref{Sec:PDP_P_J_no_dip_a} and the comments below.
This allows us to consider  very general spaces $A$ of control actions.
Moreover, differently to the previous literature, we provide a direct proof of
the viscosity solution property of $v$, which does not need to rely on   a penalized HJB equation. This is achieved by
generalizing  to the setting of the dual control problem the proof that allows to derive the HJB equation from the dynamic programming principle, see Propositions \ref{Sec:PDP_Prop_visc_subsol_property_vdelta} and \ref{Sec:PDP_Prop_visc_supersol_property_vdelta}.

\section{Piecewise Deterministic controlled Markov Processes}\label{Sec:PDP_Sec_preliminaries}

\subsection{The construction of a PDMP given its local characteristics}\label{Sec:PDP_Sec_constr_PDMP}
Given a topological space $F$, in the sequel $\mathcal B(F)$ will denote the Borel $\sigma$-field associated with $F$, and by $\mathbb{C}_b(F)$  the set of all bounded continuous functions on $F$.  The Dirac measure concentrated at some point  $x \in F$ will be denoted  $\delta_x$.

Let $(E,\mathcal E)$ be a Borel measurable space.
We will often need to construct a PDMP in  $E$ with a given triplet of local characteristics $(\phi,\lambda,Q)$. We assume that
$\phi:\R \times E \rightarrow E$ is a continuous function,
$\lambda: E \mapsto \R_+$ is a  nonnegative   continuous function satisfying
\begin{equation}\label{Sec:PDP_nulimitato}
	\sup_{x \in E} \lambda(x) < \infty,
\end{equation}
and that $Q$ maps $E$ 
into the set of probability measures on $(E, \mathcal E)$, and
is a
stochastic Feller kernel,
i.e., for all $v\in \mathbb{C}_b(E)$, the map $x\mapsto \int_{E}v(y)\,Q(x,dy)$ $(x \in E)$ is  continuous.
%

We recall the main steps of the construction of a PDMP given its local characteristics.
The existence of
a Markovian process associated with the triplet  $(\phi,\lambda, Q)$ is a well known fact (see, e.g., \cite{Da}, \cite{CoDu}). Nevertheless, we need special care in the choice
of the corresponding filtration: this will be decisive  
when we will solve
associated BSDEs and implicitly apply
a version of the martingale representation theorem.
For this reason, in the following we
will consider an explicit construction that we are going to present.
Many of the
techniques we are going to use are borrowed from the 
theory of
multivariate (marked) point processes, see 
\cite{J}, 
and also  \cite{JB}
for a more systematic 
treatise.

We start by constructing a suitable sample space to describe
the jumping mechanism of the Markov process.
Let $\Omega'$ denote the set of sequences $\omega'=(t_n,e_n)_{n \geq 1}$
in $ ((0,\infty)\times E)\cup \{(\infty,\Delta)\}$, where $\Delta\notin E$
is   adjoined to  $E$ as an isolated point, satisfying in addition
\begin{equation}\label{Sec:PDP_processo_punto}
	t_n \le t_{n+1}; \qquad
	t_n< \infty\;\Longrightarrow\; t_n < t_{n+1}.
\end{equation}
To describe  the initial condition we will use the measurable space
$(E,\mathcal E)$.
Finally, the sample space for the Markov process will be
$\Omega=E\times \Omega'$. We define canonical
functions $T_n:\Omega\to (0,\infty]$,
$E_n:\Omega\to E\cup\{\Delta\}$
as follows:  writing $\omega=(e,\omega')$ in the form
$\omega=(e,t_1,e_1,t_2,e_2,\ldots)$ we set  for $t\ge 0$ and for $n\ge 1$
$$
T_n (\omega )= t_n ,
\qquad E_n(\omega)=e_n,
\qquad T_\infty  (\omega )= \lim_{n\to\infty} t_n,
\qquad T_0(\omega)=0, \qquad E_0(\omega)=e.
$$
We also introduce the counting process  $N(s,B)=
\sum_{n\in \N}\one_{T_n\le s}\one_{E_n\in B}$, and we define the process $X: \Omega\times [0,\,\infty) \rightarrow E \cup \Delta$ setting
\begin{equation}\label{Sec:PDP_X_def}
	X_t=
	\left\{
	\begin{array}{ll}
		\phi(t-T_n, E_n)& \quad \textup{if}\,\, T_n \leq t < T_{n+1},\,\,\textup{for}\,\,n \in \N,\\
		\Delta&  \quad \textup{if}\,\,t \geq T_{\infty}.
	\end{array}	
	\right.
\end{equation}
In $\Omega$ we introduce for all $t\ge0$ the $\sigma$-algebras
$\calg_t=\sigma(N(s,B)\,:\, s\in (0,t], B\in\cale)$.
To take into account the initial condition we also introduce
the filtration $\F=(\calf_t)_{t\ge 0}$, where
$\calf_0=\cale\otimes \{\emptyset,\Omega'\}$,
and for all $t\ge0$
$\calf_t$ is the $\sigma$-algebra generated by $\calf_0$ and $\calg_t$.
$\F$ is right-continuous and will be called the natural filtration.
In the following all concepts of measurability for stochastic processes
(adaptedness, predictability etc.) refer to $\F$.
We denote by $\calf_\infty$
the $\sigma$-algebra generated by all the $\sigma$-algebras $\calf_t$.
The symbol $\calp$ denotes the $\sigma$-algebra
of $\mathbb{F}$-predictable subsets of $[0,\infty) \times \Omega$.

On the filtered sample space
$(\Omega,\F)$ we have so far introduced
the canonical marked point process $(T_n,E_n)_{n\geq1}$.
The corresponding
random measure $p$ is, for any $\omega\in\Omega$,
a $\sigma$-finite measure on $((0,\infty)\times E,\calb((0,\infty))\otimes \cale)$
defined as
\begin{eqnarray}\label{Sec:PDP_count_measure} 
	p(\omega, ds\, dy)= \sum_{n \in \N}\one_{\{ T_n(\omega) < \infty\}}\,\delta_{(T_n(\omega),E_n(\omega))}(ds\,dy),
\end{eqnarray}
where $\delta_{k}$ denotes the Dirac measure at point $k\in (0,\infty)\times E$.
For notational convenience the dependence on $\omega$ will be suppressed and, instead of $p(\omega, ds\, dy)$, it will be written $p(ds\, dy)$.

\begin{proposition}\label{Sec:PDP_P_prob_compJACOD}
	Assume that \eqref{Sec:PDP_nulimitato} holds, and
	fix $x \in E$.
	Then there exists a unique probability measure on $(\Omega, \mathcal F_{\infty})$, denoted by $\P^x$,
	such that its restriction to  $\mathcal F_0$ is
	$\delta_x$, and  the $\mathbb F$-compensator (or dual predictable projection) of the measure $p$ under $\P^x$ is the random measure
	$$
	\tilde p(ds\,dy) =\sum_{n \in \N}\one_{[T_n,\,T_{n+1})}(s)\,\lambda( \phi(s-T_n, E_n))\,Q( \phi(s-T_n, E_n), dy)\,ds.
	$$
	Moreover, $\P^{x}(T_\infty = \infty)=1$.
\end{proposition}
\proof
The result is a direct application of  Theorem 3.6 in \cite{J}.  The fact that,  $\P^{x}$-a.s., $T_{\infty}= \infty$ follows from the boundedness of  $\lambda$, see Proposition 24.6 in \cite{Da}.
\endproof
For fixed  $x \in E$, the sample path of the process  $(X_t)$ in \eqref{Sec:PDP_X_def} under $\P^{x}$ 
can be defined iteratively, by means of  $(\phi, \lambda, Q)$, in the following way.
Set
$$
F(s,x)=\exp\left(-\int_0^s \lambda (\phi(r,x))\,dr\right),
$$
we have
\begin{align}
	&\P^x(T_1 > s)=F(s,x),\label{Sec:PDP_A}\\
	&\P^x(X_{T_1} \in B 
	|\,T_1)= Q(x, B) 
	,\label{Sec:PDP_B}
\end{align}
on $\{T_1 < \infty\}$,
and, for every $n > 1$,
\begin{align}
	&\P^x(T_{n+1} > s\,|\,\mathcal F_{T_n})=\exp\left(-\int_{T_n}^s \lambda (\phi(r-T_n,X_{T_n}))\,dr\right),\label{Sec:PDP_A_k}\\
	&\P^x(X_{T_{n+1}} \in B 
	\,|\,\mathcal F_{T_n},\,T_{n+1})= Q(\phi(T_{n+1}-T_n,X_{T_n}), B), 
	\label{Sec:PDP_B_k}
\end{align}
on $\{T_n < \infty\}$.

\begin{proposition}\label{Sec:PDP_P_markovianity}
	In the probability space $\{\Omega, \mathcal F_{\infty}, \P^{x}\}$ the process $X$ has distribution $\delta_x$ at time zero, and it is a homogeneous Markov process, i.e.,
	for any $x \in E$,  nonnegative times $t$, $s$, $t \leq s$, and for every  bounded measurable function $f$,
	\begin{equation}\label{Sec:PDP_MarkovProperty}
		\E^x[f(X_{t+s})\,|\,\mathcal{F}_t] = P_s(f(X_t)).
	\end{equation}
	where
	$P_t f(x):= \E^x[f (X_t)]$.
\end{proposition}	
\proof	
From \eqref{Sec:PDP_A_k}, 
taking into account   the semigroup property $\phi(t+s,x)=\phi(t,\phi(s,x))$, we have
\begin{align}\label{Sec:PDP_EE}
	&\P^x(T_{n+1} > t+s\,|\,\mathcal F_t)\,\one_{\{ t \in [T_n,\,T_{n+1})\}} \nonumber\\
	&=
	\frac{\P^x(T_{n+1} > t+s\,|\,\mathcal F_{T_{n})}}
	{\P^x(T_{n+1} > t\,|\,\mathcal F_{T_n})}\,\one_{\{ t \in [T_n,\,T_{n+1})\}}\nonumber\\
	&= \exp\left(-\int_{t}^{t+s} \lambda (\phi(r-T_{n},X_{T_n}))\,dr\right)\,\one_{\{ t \in [T_n,\,T_{n+1})\}}\nonumber\\
	&= \exp\left(-\int_{0}^{s} \lambda (\phi(r+t-T_{n},X_{T_n}))\,dr\right)\,\one_{\{ t \in [T_n,\,T_{n+1})\}}\nonumber\\
	&= \exp\left(-\int_{0}^{s} \lambda (\phi(r, X_t))\,dr\right)\,\one_{\{ t \in [T_n,\,T_{n+1})\}}\nonumber\\
	&=F(s,X_t)\,\one_{\{ t \in [T_n,\,T_{n+1})\}}.
\end{align}
Hence, denoting $N_t = N(t,E)$, it follows from \eqref{Sec:PDP_EE} that
\begin{displaymath}
	\P^x(T_{N_t+1} > t+s\,|\,\mathcal F_t) = F(s,X_t);
\end{displaymath}	
in other words, conditional on $\mathcal F_t$, the jump time after $t$ of a PDMP started at $x$ has the same distribution as the first jump time of a PDMP started at $X_t$.
Since the remaining interarrival times and postjump positions are independent on  the past, we have shown that  \eqref{Sec:PDP_MarkovProperty} holds for every bounded measurable function $f$.
\endproof

\begin{remark}\label{Sec:PDP_R_non_quasi_left}
	In the present paper we restrict the analysis to the case of PDMPs on a domain $E$ with no boundary. 
	Indeed, we have in mind to apply techniques based on BSDEs driven by the compensated random measure  associated to the PDMP (see Section \ref{Sec:PDP_Sec_ConstrainedBSDE}), 
	and  the presence of the jumps at the boundary of the domain would induce
	discontinuities in the compensator, which corresponds to very technical difficulties in the study of the associated BSDE.
	
	More precisely, 
	when the process reaches the boudary of the domain, according to (26.2) in \cite{Da}, the compensator of the counting  measure $p$ in  \eqref{Sec:PDP_count_measure} admits the form
	$$
	\tilde p(ds\,dy)=(\lambda(X_{s-})\,ds + d p^{\ast}_s)\,Q(X_{s-}, dy),
	$$
	where
	$$
	p^\ast_s= \sum_{n=1}^{\infty} \one_{\{s \geq T_n\}}\,\one_{\{X_{{T_n}-} \in \Gamma\}}
	$$
	is the process counting the number of jumps of $X$ from the active boundary $\Gamma\in \partial E$ (for the precise definition of $\Gamma$ see page 61 in \cite{Da}).
	In particular, the compensator $\tilde p$ can be rewritten as $\tilde p(ds\,dy)= d A_s\,Q(X_{s-}, dy)$, where $A_s= \lambda(X_{s-})\,ds + d p^{\ast}_s$ is a predictable and discontinuous process, with jumps $\Delta A_s= \one_{X_{s-}\in \Gamma}$. 
	
	For  BSDEs driven by random measures with discontinuous compensator, existence and uniqueness results were at disposal  only in particular frameworks, see e.g. \cite{CohenElliott} for the deterministic case, and counter-examples were provided in the general case, see Section 4.3 in \cite{CFJ}.
	Only recently  this problem was faced and solved in a general context in   \cite{BandiniBSDE}, were  a technical condition is provided in order to achieve  existence and uniqueness of the BSDE.
	The mentioned condition turns out to be verified in the case of control problems related to PDMPs with discontinuities at the boundary of the domain, see Remark 4.5 in \cite{BandiniBSDE}.
	This fact opens to the possibility to apply   the BSDEs techniques also in this context, which  will be analyzed in a second work by the author.

%
\end{remark}

\subsection{Optimal control of PDMPs}\label{Sec:PDP_Sec_control_problem}
In the present section we aim at formulating
an  optimal control problem for piecewise deterministic Markov processes, and to discuss its solvability.
The PDMP state space $E$ will be an open subset of  $\R^d$, and $\mathcal E$ the corresponding $\sigma$-algebra. In addition,  we introduce  a 
Borel space $A$, endowed with its $\sigma$-algebra $\mathcal{A}$, called the space of control actions.
The additional hypothesis that $A$ is compact  is not necessary for the majority of  the results, and will be explicitly asked  whenever it will be needed.
The other data of the problem
consist in three functions $f$, $h$ and $\lambda$ on $E \times A$, and in a  probability transition $Q$ from $(E \times A, \mathcal E \otimes \mathcal A)$ to $(E, \mathcal E)$, satisfying the following conditions.

\medskip

\noindent \textbf{(H$\textup{h$\lambda$Q}$)}
\begin{itemize}
	\item[(i)]
	$h: E \times A \mapsto E$ is a bounded uniformly continuous function satisfying
	\begin{equation*}
		\left\{
		\begin{array}{ll}
			\forall x,\,x' \in E,\,\,\textup{and}\,\,\forall a, a' \in A,\quad \quad & |h(x,a)-h(x',a')| \leqslant L_h \,(|x-x'|+|a-a'|),\\
			\forall x\in E\,\,\textup{and}\,\,\forall a \in A, \quad & |h(x,a)| \leqslant M_h,
		\end{array}
		\right. 
	\end{equation*}
	where $L_h$ and $M_h$ are constants independent of $a,a' \in A$, $x,x' \in E$.
	\item[(ii)]
	$\lambda: E \times A \mapsto \R^+$ is a  nonnegative bounded uniformly continuous function, satisfying
	\begin{equation}\label{Sec:PDP_lambdaContrfinite}
		\sup_{x \in E, a \in A} \lambda(x,a) < \infty.
	\end{equation}	
	\item[(iii)]
	$Q$ maps 
	$E \times A$ 
	into the set of probability measures on $(E, \mathcal E)$, and
	is a
	stochastic Feller kernel.
	i.e., for all $v\in \mathbb{C}_b(E)$, the map $(x,a)\mapsto \int_{\R^d}v(y)\,Q(x,a,dy)$ is continuous (hence it belongs to $\mathbb{C}_b(E \times A)$). 
\end{itemize}

\noindent \textbf{(H$\textup{f}$)}\quad
$f: E \times A \mapsto \R^+$ is a  nonnegative bounded uniformly continuous function.
In particular, there exists a  positive constant $M_f$ such that
$$
0 \leqslant f(x,a)\leqslant M_f, \quad \forall \,x \in E,\,a \in  A.
$$

\medskip

The requirement that $Q(x,a, \{x\})=0$ for all $x \in E$, $a \in A$ is natural in many applications, but here is not needed.
$h$, $\lambda$ and  $Q$  depend on the control parameter $a \in A$ and play respectively the role  of and controlled drift, controlled jump rate and controlled probability transition. Roughly speaking, we may control the dynamics of the process by changing dynamically its deterministic drift,  its jump intensity and its post jump distribution.

Let us  give a more precise definition of the  optimal control problem under study.
To this end, we first construct $\Omega$, $\mathbb F=(\mathcal F_t)_{t \geq 0}$, $\mathcal F_{\infty}$ as in the previous paragraph.

We will consider the class of
\emph{piecewise open-loop
	controls}, first introduced in 
\cite{Ver} and  often adopted in this context,  see for instance \cite{Da}, 
\cite{CoDu}, \cite{Almudevar}.
Let $X$ be the (uncontrolled) process constructed in a canonical way from a marked point process $(T_n,E_n)$ as in Section \ref{Sec:PDP_Sec_constr_PDMP}.
The class of admissible control law $\mathcal A_{ad}$ is the set of all Borel-measurable maps $\alpha: [0,\,\infty) \times E \rightarrow A$, and the control applied to  $X$ is of the form:
\begin{equation}\label{Sec:PDP_open_loop_controls}
	\alpha_t= \sum_{n=1}^{\infty}\alpha_n(t-T_n,E_n)\,\one_{[T_n,\,T_{n+1})}(t).
\end{equation}
In other words, at each jump time $T_n$, we choose an open loop control $\alpha_n$ depending on the initial
condition 
$E_n$ to be used until the next jump time.

We define the controlled  process $X: \Omega \times [0,\,\infty) \rightarrow E \cup \{\Delta\}$ setting
\begin{equation}\label{Sec:PDP_controlledX}
	X_t = \phi^{\alpha_n}(t- T_n, E_n),\quad t \in [T_n,\, T_{n+ 1})
\end{equation}
where $\phi^\beta(t, x)$
is the unique solution to the ordinary differential equation
$$
\dot x(t) = h(x(t), \beta(t)),\quad x(0)=x \in E,
$$
with $\beta$ a $\mathcal A$-measurable function.
Then, for every starting point $x \in E$ and for each $\alpha \in \mathcal A_{ad}$, by Proposition \ref{Sec:PDP_P_prob_compJACOD} there exists a unique probability measure on $(\Omega, \mathcal F_{\infty})$,  denoted by $\P^{x}_{\alpha}$, such that its restriction to $\mathcal F_0$ is  $\delta_x$, and
the  $\mathbb F$-compensator 
under $\P^{x}_{\alpha}$ of the measure $p(ds\,dy)$ is
$$
\tilde{p}^{\alpha}(ds\,dy)=
\sum_{n=1}^{\infty}\one_{[T_n,\,T_{n+1})}(s)\,\lambda(X_{s},\alpha_n(s-T_n,  E_n)
)\,Q(X_s,\alpha_n(s-T_n,  E_n)
, dy)\,ds.
$$
According to Proposition \ref{Sec:PDP_P_markovianity}, under $\P^{x}_{\alpha}$ the process $X$ in \eqref{Sec:PDP_controlledX} is 
markovian  with respect to  $\mathbb F$.

Denoting by $\E^{x}_{\alpha}$ the expectation under $\P^{x}_{\alpha}$, we finally define, for $x \in E$ and $\alpha \in \mathcal A_{ad}$, the  functional cost
\begin{equation}\label{Sec:PDP_functional_cost}
	J(x,\alpha) =
	\E^{x}_{\alpha}\left[\int_{0}^{\infty} e^{-\delta\,s}\, f(X_{s},\alpha_s)\,ds\right]
\end{equation}
and the value function of the control problem
\begin{equation}\label{Sec:PDP_value_function}
	V(x) = \inf_{\alpha \in \mathcal{A}_{ad}}J(x,\alpha),
\end{equation}
where   $\delta \in (0,\,\infty)$ is a  discounting factor that will be fixed from here on.
By the boundedness assumption 
on $f$, both $J$ and $V$ are well defined and bounded.

Let us consider the  Hamilton-Jacobi-Bellman equation (for short, HJB equation) associated to the optimal control problem: this is the following elliptic nonlinear equation on $[0,\,\infty)\times E$:
\begin{eqnarray}\label{Sec:PDP_HJB}
	H^{v}(x,v,Dv) = 0,
\end{eqnarray}
where
\begin{align*}
	H^\psi(z, v, p) &= \sup_{a \in A}\left\{ \delta\,v  -h(z,a) \cdot  p -\int_{E} (\psi(y)-\psi(z))\, \lambda(z, a)\,  Q(z,a, dy) - f(z,a) \right\}.
\end{align*}
\begin{remark}\label{Sec:PDP_Rem_form_HJB}
	The HJB
	equation \eqref{Sec:PDP_HJB} can be rewritten as
	\begin{equation}\label{Sec:PDP_HJB2}
		\delta\,v(x) =\sup_{a \in A}\left\{\mathcal L^a v(x)  + f(x,a) \right\}=0,
	\end{equation}
	where 
	$\mathcal L^a$ is the operator depending on $a \in A$
	defined  as
	\begin{equation}\label{Sec:PDP_ext_gen_HJB}
		\mathcal L^a v(x) := h(x,a) \cdot \nabla v(x) +\lambda(x, a)\int_{E} (v(y)-v(x))\,   Q(x,a, dy). 
	\end{equation}
\end{remark}
Let us recall the following facts.
Given a  locally bounded function $z: E \rightarrow \R$, we define its lower semicontinuous (l.s.c. for short) envelope $z_\ast$, and its upper semicontinuous (u.s.c. for short) envelope $z^\ast$, by
\begin{displaymath}
	z_{\ast}(x) = \liminf_{\substack{y \rightarrow x\\ y \in E}}z(y),\qquad z^{\ast}(x) = \limsup_{\substack{y \rightarrow x\\ y \in E}}z(y), \quad \text{for all}\,\, x \in E.
\end{displaymath}

\begin{definition}\label{Sec:PDP_Def_viscosity_sol_HJB} Viscosity solution to \eqref{Sec:PDP_HJB}.
	\begin{itemize}
		\item[(i)] A locally bounded u.s.c. function $w$ on $E$ is called a \textbf{viscosity supersolution} (resp. \textbf{viscosity subsolution}) of \eqref{Sec:PDP_HJB} if
		$$
		H^w(x_0,w(x_0),D\varphi(x_0)) \geqslant \,\,(resp.\,\, \leqslant ) \,\,0.
		$$
		for any  $x_0\in E$ and for any  $\varphi \in C^1(E)$ such that
		$$
		(u-\varphi)(x_0)= \min_{E}(u-\varphi)\,\,(resp.\,\, \max_{E}(u-\varphi)).
		$$
		\item[(ii)] A function $z$ on $E$ is called  a \textbf{viscosity solution} of \eqref{Sec:PDP_HJB} if it is locally bounded  and its u.s.c. and l.s.c. envelopes are respectively subsolution and supersolution of \eqref{Sec:PDP_HJB}. 
	\end{itemize}
\end{definition}

The HJB equation \eqref{Sec:PDP_HJB} admits a unique continuous solution, which coincides with the value function $V$ in \eqref{Sec:PDP_value_function}.
The following result is stated in Theorem 7.5  in \cite{Da-Fa}.

\begin{theorem}\label{Sec:PDP_Thm_unique_viscosity_sol_HJB}
	Let \textup{\textbf{(H$\textup{h$\lambda$Q}$)}} and  \textup{\textbf{(H$\textup{f}$)}} hold, and assume that $A$ is compact.
	Then,  the value function  $V$ 
	of the PDMPs optimal control problem is 
	the unique continuous viscosity solution to \eqref{Sec:PDP_HJB}.
\end{theorem}

\section{Control randomization and dual optimal control problem}\label{Sec:PDP_Section_dual_control}
In this section we start to  implement the control randomization method.
In the first step, for an initial time $t \geq 0$ and starting point $x \in E$, we construct an (uncontrolled) Markovian pair of PDMPs $(X,I)$ by specifying its local characteristics, see \eqref{Sec:PDP_h_XI}-\eqref{Sec:PDP_lambda_XI}-\eqref{Sec:PDP_Q_XI} below. Next we formulate an auxiliary optimal control problem where,
roughly speaking, we optimize a  functional cost by modifying the intensity of the process $I$ over a suitable family.

This dual problem  is studied in Section \ref{Sec:PDP_Sec_ConstrainedBSDE} by means of a suitable class of BSDEs.
In Section \ref{Sec:PDP_Section_nonlinear_IPDE} we will show that
the same class of BSDEs provides a probabilistic representation of the value function introduced in the previous section.
As a byproduct, we also get that the dual value function coincides with the one associated to the original optimal control problem.

\subsection{A dual control system}\label{Sec:PDP_Section_control_rand}


Let $E$ still denote an open subset of $\R^d$ with $\sigma$-algebra $\mathcal E$, and  $A$ be a Borel space with corresponding $\sigma$-algebra $\mathcal A$. Let moreover $h$, $\lambda$ and  $Q$ be respectively two real  functions on $E \times A$ and  a   probability transition  from $(E \times A, \mathcal E \otimes \mathcal A)$, satisfying  \textbf{(H\textup{h$\lambda$Q})} as before.
We denote by $\phi(t, x,a)$
the unique solution to the ordinary differential equation
$$
\dot x(t) = h(x(t), a), \quad x(0)=x \in E, \,\,a \in A.
$$
In particular, $\phi(t,x,a)$ corresponds to the function $\phi^\beta(t,x)$, introduced in Section \ref{Sec:PDP_Sec_control_problem}, when $\beta(t)\equiv a$.
Let us now introduce another finite measure $\lambda_0$ on $(A, \mathcal A)$ satisfying the following assumption:

\medskip

\noindent \noindent \textbf{(H\textup{$\lambda_0$})}\quad $\lambda_0$ is a finite measure on $(A, \mathcal A)$ with full topological support.

\medskip

\noindent The existence of such a measure is guaranteed by the fact that the space $A$ is metric separable. We define
\begin{eqnarray}
	\tilde \phi(t,x,a)&:=& (\phi(t,x,a) \quad a),
	\label{Sec:PDP_h_XI}\\ 
	\tilde{\lambda}(x,a) &:=& \lambda(x,a)+ \lambda_0(A),\label{Sec:PDP_lambda_XI}\\ 
	\tilde{Q}(x,a,dy\,db) &:=& \frac{\lambda(x,a)\,Q(x,a,dy)\,\delta_{a}(db) + \lambda_0(db)\,\delta_{x}(dy)}{\tilde{\lambda}(x,a)}.\label{Sec:PDP_Q_XI} 
\end{eqnarray}
We wish to construct a PDMP $(X,I)$ as in Section \ref{Sec:PDP_Sec_constr_PDMP} but with enlarged state space $E \times A$ and  local characteristics $(\tilde \phi,\tilde \lambda, \tilde Q)$.


Firstly, we need to introduce  a suitable sample space to describe the jump mechanism of the  process $(X,I)$ on $E \times A$.
Accordingly, we  set  $\Omega'$ as the set of sequences $\omega'=(t_n, e_n, a_n)_{n \geq 1}$ contained in $((0,\,\infty) \times E \times A) \cup \{ (\infty, \Delta, \Delta')\}$, where $\Delta \notin E$ (resp. $\Delta' \notin A$) is adjoined to $E$ (resp. to $A$) as an isolated point, satisfying \eqref{Sec:PDP_processo_punto}. In the sample space $\Omega = \Omega' \times E \times A$ we defined the random variables $T_n : \Omega \rightarrow (0,\,\infty]$, $E_n : \Omega \rightarrow E \cup \{\Delta\}$, $A_n : \Omega \rightarrow A \cup \{\Delta'\}$, as follows: writing $\omega =(e,a,\omega')$ in the form $\omega = (e,a,t_1, e_1, a_1, t_2, e_2, a_2,...)$ we set for $t \geq 0$ and for $n \geq 1$
\begin{align*}
	& T_n (\omega )= t_n, \qquad T_\infty  (\omega )= \lim_{n\to\infty} t_n,
	\qquad T_0(\omega)=0,\\
	& E_n(\omega)=e_n,\qquad A_n(\omega)=a_n,
	\qquad E_0(\omega)=e, \qquad A_0(\omega)=a.
\end{align*}
We  define the process $(X,I)$ on $(E \times A) \cup \{\Delta, \Delta'\}$ setting
\begin{align}\label{Sec:PDP_XI_def}
	(X,I)_t&=
	\left\{
	\begin{array}{ll}
		(\phi(t-T_n, E_n,A_n),A_n)& \quad \textup{if}\,\, T_n \leq t < T_{n+1},\,\,\textup{for}\,\,n \in \N,\\
		(\Delta,\Delta')&  \quad \textup{if}\,\,t \geq T_{\infty}.
	\end{array}	
	\right.
\end{align}

In $\Omega$ we introduce for all $t\geq 0$ the $\sigma$-algebras
$\mathcal G_t=\sigma(N(s,B)\,:\, s\in (0,t], B\in\mathcal E \otimes \mathcal A)$ generated
by the counting processes  $N(s,A)=
\sum_{n \in \N}\one_{T_n\le s}\one_{E_n\in A}$ and the
$\sigma$-algebra $\mathcal F_t$ generated by $\calf_0$ and $\calg_t$, where
$\calf_0=\mathcal E\otimes \cala \otimes \{\emptyset,\Omega'\}$.
We still denote by $\F = (\mathcal F_t)_{t \geq 0}$ and $\mathcal P$ the corresponding filtration and predictable $\sigma$-algebra.
The random measure $p$ is now defined on $(0,\,\infty) \times E \times A$ as
\begin{equation}\label{Sec:PDP_p_dual}
	p(ds\,dy\,db)= \sum_{n\in \N} \one_{\{T_n,E_n, A_n\}}(ds\,dy\,db).
\end{equation}
Given any starting point $(x,a) \in E \times A$,
by  Proposition \ref{Sec:PDP_P_prob_compJACOD}, there exists a unique
probability measure    on $(\Omega, \mathcal F_{\infty})$, denoted by $\P^{x,a}$, such that its restriction to  $\mathcal F_0$ is $\delta_{(x,a)}$  and the $\mathbb F$-compensator  of the measure $p(ds\,dy\,db)$ under $\P^{x,a}$ is the random measure
$$
\tilde p(ds\,dy\,db)= \sum_{n \in \N}\one_{[T_n,\,T_{n+1})}(s)\,\Lambda( \phi(s-T_n, E_n,A_n),A_n, dy\,db)\,ds,
$$
where
$$
\Lambda(x,a,dy\,db)= \lambda(x,a)\,Q(x,a,dy)\,\delta_a(db)+ \lambda_0(db)\,\delta_x(dy), \quad \forall (x,a)\in E \times A.
$$
We set $q=p-\tilde{p}$, the compensated martingale measure associated to $p$.

As in Section \ref{Sec:PDP_Sec_constr_PDMP},  the sample path of a process $(X,I)$ with values in $E \times A$, starting from a fixed initial point $(x,a) \in E \times A$ at time zero, can be defined iteratively by means of its local characteristics $(\tilde h, \tilde \lambda, \tilde Q)$ in the following way. Set
$$
F(s,x,a)=\exp\left(-\int_0^s (\lambda (\phi(r,x,a),a) + \lambda_0(A))\,dr\right),
$$
we have
\begin{align}
	&\P^{x,a}(T_1 > s)=F(s,x,a),\label{Sec:PDP_Abis}\\
	&\P^{x,a}(X_{T_1} \in B, I_{T_1}\in C, T_1 < \infty\, |\,T_1)= \tilde Q(x, B \times C)\,\one_{\{T_1 < \infty\}},\label{Sec:PDP_Bbis}
\end{align}
and, for every $n > 1$,
\begin{align}
	&\P^{x,a}(T_{n+1} > s\,|\,\mathcal F_{T_n})=\exp\left(-\int_{T_n}^s (\lambda (\phi(r-T_n,X_{T_n},I_{T_n}), I_{T_n})+ \lambda_0(A))\,dr\right),\label{Sec:PDP_A_kbis}\\
	&\P^{x,a}(X_{T_{n+1}} \in B, I_{T_{n+1}} \in C 
	|\,\mathcal F_{T_n},\,T_{n+1})= \tilde Q(\phi(T_{n+1}-T_{n},X_{T_n},I_{T_n}), I_{T_n}, B \times C) 
	,\label{Sec:PDP_B_kbis}
\end{align}
on $\{T_n < \infty\}$.

Finally, an application of Proposition \ref{Sec:PDP_P_markovianity} provides that $(X,I)$ is a Markov process on $[0,\,\infty)$ with respect to $\F$.
For
every 
real function taking values in $E \times A$, the infinitesimal generator is given by
\begin{align*}
	\mathcal{L}\varphi(x,a) :=& h(x,a)\cdot \nabla_x \varphi(x,a) + \int_{E}(\varphi(y,a) - \varphi(x,a))\, \lambda(x,a)\,Q(x,a,dy)\\
	&+\int_{A}(\varphi(x,b) - \varphi(x,a))\, \lambda_0(db).
\end{align*}
For our purposes, it will be not necessary to specify the domain of the previous operator (for its formal definition   we  refer to Theorem 26.14 in \cite{Da}); in the sequel the operator $L$  will be applied  to test functions with suitable regularity. 

\subsection{The dual optimal control problem}\label{Sec:PDP_Section_dual_optimal_control}

We now introduce a dual  optimal control problem associated to the process $(X,I)$, and formulated in a weak form. For fixed $(x,a)$, we consider a family of probability measures $\{\P_{\nu}^{x,a},\,\nu \in \mathcal V\}$ in the space $(\Omega, \mathcal F_{\infty})$ 
whose effect   is to change the stochastic intensity of the process $(X,I)$. 

Let us proceed with precise definitions.
We still  assume that {\bf (Hh$\lambda$Q)}, {\bf (H$\lambda_0$)} and {\bf (H$f$)} hold.
We recall that $\mathbb{F}= (\mathcal{F}_t)_{t \geqslant 0}$ is the augmentation of the natural filtration generated by $p$ in \eqref{Sec:PDP_p_dual}.
We define
$$
\mathcal{V} = \{ \nu: \Omega \times [0,\,\infty) \times A \rightarrow (0,\,\infty)\,\, \mathcal{P}\otimes \mathcal A\text{-measurable and bounded}\}.
$$
For every $\nu \in \mathcal{V}$,
we consider the predictable random measure
\begin{align}\label{Sec:PDP_dual_comp}
	\tilde{p}^\nu(ds\,dy\,db)&:= \nu_s(b) \, \lambda_0 (db) \, \delta_{\{ X_{s-} \}}(dy) \, ds\nonumber\\
	& +\,
	\lambda(X_{s-},\,I_{s-})\,Q(X_{s-},\,I_{s-},\,dy)\, \delta_{\{ I_{s-} \}}(db)\, ds.
\end{align}

In particular, by the Radon Nikodym theorem one can find two nonnegative functions $d_1$, $d_2$ defined on $\Omega \times [0,\,\infty) \times E \times A$, $\mathcal P \otimes \mathcal E \otimes \mathcal A$, such that
\begin{eqnarray*}
	\lambda_0(db) \, \delta_{\{ X_{t-} \}}(dy) \, dt &=& d_1(t,y,b)\, \tilde{p}(dt\,dy\,db)
	\\
	\lambda(X_{t-},\,I_{t-},\,dy) \, \delta_{\{ I_{t-} \}}(db) \, dt &=&
	d_2(t,y,b)\, \tilde{p}(dt\,dy\,db),
	\\
	d_1(t,y,b) + d_2(t,y,b) &=& 1, \qquad \tilde{p}(dt\,dy\,db)-a.e.
\end{eqnarray*}
and we have
$
d \tilde{p}^{\nu} = (\nu \, d_1  + d_2)\, d\tilde{p}$.
For any $\nu \in \mathcal{V}$, consider then the Dol\'eans-Dade exponential local martingale
$L^\nu$ defined setting
\begin{align}
	L_s^\nu &= \exp\bigg(\int_0^s\!\int_{E\times A}\log(\nu_r(b) \, d_1(r,y,b) \, + d_2(r,y,b))\,p(dr \, dy \,db)\nonumber\\
	& \qquad \qquad - \int_{0}^{s}\!\int_{A}(\nu_r(b) - 1)\lambda_0(db)\,dr\bigg)
	\nonumber\\
	&= e^{\int_{0}^{s}\int_{A}(1 - \nu_r(b) )\lambda_0(db)\,dr} \prod_{n \geqslant 1: T_{n} \leqslant s}(\nu_{T_{n}}(A_n)\,d_1(T_{n},E_n,A_n) + d_2(T_{n},E_n,A_n)),\label{Sec:PDP_Lnu}
\end{align}
for $s\geq 0$. 
When $(L^\nu_t)_{t \geq 0}$ is a true martingale, 
for every time $T>0$ we can define  a probability measure $\P^{x,a}_{\nu, T}$ equivalent to $\P^{x,a}$ on $(\Omega,\,\mathcal{F}_T)$ setting
\begin{equation}\label{Sec:PDP_PnuT}
	\P^{x,a}_{\nu, T}(d\omega)=L_T^\nu(\omega)
	\,\P^{x,a}(d\omega).
\end{equation}
By the Girsanov theorem for point processes (see Theorem 4.5 in \cite{J})
the restriction of the random measure $p$ to $(0,T]\times E\times A$
admits $\tilde{p}^\nu=(\nu\,d_1 + d_2)\,\tilde{p}$
as compensator   under $\P^{x,a}_{\nu, T}$.
We set $q^\nu: = p-\tilde{p}^\nu$. 
and we denote by $\mathbb{E}_{\nu, T}^{x,a}$ the expectation operator under $\P_{\nu, T}^{x,a}$.
We recall the following result, which   is a direct consequence  of  
Lemma 3.2 in \cite{BandiniFuhrman}.
\begin{lemma}\label{Sec:PDP_lemma_P_nu_martingale}
	Let assumptions {\bf (Hh$\lambda$Q)} and  {\bf (H$\lambda_0$)} hold. Then, for every  $(x,a)\in E\times A$  and
	$\nu \in \mathcal{V}$, under the probability
	$\P^{x,a}$, 
	the process $(L^\nu_t)_{t \geq 0}$ is a   martingale.
	Moreover, for every time $T>0$,  $L_T^\nu$ is square integrable,
	and, for every  $\mathcal{P}_T \otimes \mathcal E\otimes \cala$-measurable  function $H: \Omega \times [0,T] \times E \times A \rightarrow \R$ such that $\sperxa{\int_0^T \int_{E\times A} |H_s(y,b)|^2\, \tilde{p}(ds\,dy\,db)}$ $< \infty$,  the process $\int_{0}^{\cdot}\int_{E \times A}H_s(y,b)\, q^\nu(ds\,dy\, db)$ is a $\P^{x,a}_{\nu, T}$-martingale on $[0,T]$.
\end{lemma}

Our  aim is to extend the previous construction to get a suitable probability measure  on $(\Omega, \mathcal F_{\infty})$.
We have the following result.
\begin{proposition}\label{P_Prob_infty}
	Let assumptions {\bf (Hh$\lambda$Q)} and  {\bf (H$\lambda_0$)} hold. Then, for every  $(x,a)\in E\times A$  and
	$\nu \in \mathcal{V}$, there exists a unique probability
	$\P^{x,a}_{\nu}$ on $(\Omega, \mathcal F_{\infty})$, under which the random measure  $\tilde{p}^\nu$ in \eqref{Sec:PDP_dual_comp} is the compensator of the measure $p$ in \eqref{Sec:PDP_p_dual} on $(0,\,\infty)\times E \times A$. Moreover, for any time $T >0$,  the restriction of $\P^{x,a}_{\nu}$ on $(\Omega, \mathcal F_T)$ is given by the  probability measure 
	$\P^{x,a}_{\nu,T}$	in \eqref{Sec:PDP_PnuT}. 
\end{proposition}


\proof
For simplicity, in the sequel we will drop the dependence of  $\P^{x,a}$ and $\P^{x,a}_\nu$  on $(x,a)$, which will be denoted respectively by  $\P$ and $\P^\nu$.

We notice that $\mathcal F_{T_n}$ $=$ $ \sigma(T_1,E_1,A_1,...,T_n,E_n,A_n)$ defines an increasing family of sub
$\sigma$-fields of $\mathcal F_{\infty}$ such that $\mathcal F_{\infty}$ is generated by $\bigcup_n \mathcal F_{T_n}$.
The idea is then  to provide
a  family
$\{\P^\nu_{n}\}_n$ of probability measures on $(\Omega, \mathcal F_{T_n})$ under which $\tilde{p}^\nu$ is the compensator of the measure $p$ on $(0,\,T_n]\times E \times A$, and  which is consistent (i.e., $\P^\nu_{n+1} \big|_{\mathcal F_{T_n}} = \P^\nu_{n}$).
Indeed, if we have at disposal such a family of probabilities,
we can  naturally define on $\bigcup_n \mathcal F_{T_n}$ a set function  $\P^{\nu}$ verifying the desired property, by setting $\P^{\nu}(B):= \P^\nu_{n}(B)$ for every $B \in \mathcal F_{T_n}$,  $n \geq 1$.
Finally,  to conclude  we would need to  show  that $\P^{\nu}$ is countably additive on  $\bigcup_n \mathcal F_{T_n}$, and therefore  can be extended uniquely to $\mathcal F_{\infty}$.

Let us proceed by steps. For every $n \in \N$, we set
\begin{equation}\label{Sec:PDP_Pn}
	d\P^\nu_{n}:= L^\nu_{T_{n}}\,d \P \quad \textup{on } (\Omega, \mathcal F_{T_n}),
\end{equation}
where $L^\nu$ is given by \eqref{Sec:PDP_Lnu}.
Notice that, for every $n \in \N$, the probability $\P^\nu_{n}$ is well defined. Indeed, recalling  the boundedness properties of $\nu$ and $\lambda_0$, we have
\begin{align}
	L_{T_n}^\nu&= e^{\int_{0}^{T_n}\int_{A}(1 - \nu_r(b) )\lambda_0(db)\,dr} \prod_{k=1}^n (\nu_{T_{k}}(A_k)\,d_1(T_{k},E_k,A_k) + d_2(T_{k},E_k,A_k))\nonumber\\
	&\leq  (||\nu||_{\infty})^n \, e^{\lambda_0(A)\,T_n},\label{Lwedge_unif_int}
\end{align}
and since $T_n$ is exponentially distributed (see \eqref{Sec:PDP_A_k}),  we get
\begin{align*}
	\sper{L_{T_n}^\nu}\leq (||\nu||_{\infty})^n \, \sper{e^{\lambda_0(A)\,T_n}}< \infty.
\end{align*}
Then, arguing as in the proof of  the Girsanov theorem for point process (see, e.g., the comments after Theorem 4.5 in \cite{J}), it can be proved that the restriction of the random measure $p$ to $(0,T_n]\times E\times A$
admits $\tilde{p}^\nu=(\nu\,d_1 + d_2)\,\tilde{p}$
as compensator   under $\P^{\nu}_{n}$.
Moreover,
$\{\P^{\nu}_{n}\}_n$ is a consistent family of probability measures on $(\Omega, \mathcal F_{T_n})$, namely
\begin{equation}\label{Sec:PDP_consistencyPn}
	\P^{\nu}_{n+1} \big|_{\mathcal F_{T_n}} = \P^{\nu}_{n},\quad n \in \N.
\end{equation}
Indeed, taking into account definition \eqref{Sec:PDP_Pn}, it is easy to see that identity \eqref{Sec:PDP_consistencyPn} is equivalent to
\begin{equation}\label{Sec:PDP_consistencyPn2}
	\sper{L^\nu_{T_{n}}|{\mathcal F_{T_{n-1}}}} = L^\nu_{T_{n-1}},\quad n \in \N.
\end{equation}
By  Corollary 3.6, Chapter II, in \cite{Re-Yor}, and taking into account the estimate  \eqref{Lwedge_unif_int},  it follows that the process  $(L^\nu_{t \wedge T_n})_{t \geq 0}$ is a uniformly integrable martingale. Then, identity \eqref{Sec:PDP_consistencyPn2} follows  from the optional stopping theorem for uniformly integrable martingales (see, e.g., Theorem 3.2, Chapter II, in \cite{Re-Yor}).

At this point, we define the following probability measure on $\bigcup_n \mathcal F_{T_n}$:
\begin{equation}
	\label{Sec:PDP_Pinfty}
	\P^{\nu}(B):= \P^{\nu}_{n}(B), \quad  B \in \mathcal F_{T_n}, \,\,n \in \N.
\end{equation}
In order to  get the
desired probability measure on $(\Omega, \mathcal F_{\infty})$,
we need to show that $\P^{\nu}$  in \eqref{Sec:PDP_Pinfty}  is $\sigma$-additive on $\bigcup_n \mathcal F_{T_n}$: in this case, $\P^{\nu}$ can indeed be  extended uniquely to $\mathcal F_{\infty}$, see  Theorem 6.1 in \cite{jacodprotter}.

Let us then prove that $\P^{\nu}$ in \eqref{Sec:PDP_Pinfty} is countably additive on $\bigcup_n \mathcal F_{T_n}$.
To this end, let us introduce the product space $\tilde E_{\Delta}^{\N}:= (E \times A \times [0,\,\infty)\cup\{(\Delta, \Delta', \infty)\})^{\N}$, with associated Borel $\sigma$-algebra $\mathcal{\tilde E}_{\Delta}^{\N \otimes }$.  For every $n \in \N$, we define the following probability measure on $(E_{\Delta}^{n},\mathcal{\tilde E}_{\Delta}^{n \otimes })$:
\begin{equation}\label{Sec:PDP_Qn}
	\Q^\nu_n(A):=\P^\nu_n(\omega: \pi_n(\omega)\in A),\quad A \in \tilde E_{\Delta}^{n},
\end{equation}
where $\pi_n = (T_1, E_1, A_1,..., T_n, E_n, A_n)$. The consistency property \eqref{Sec:PDP_consistencyPn} of the family $(\P^\nu_n)_n$ implies that
\begin{equation}\label{Sec:PDP_coerenceQn}
	\Q^\nu_{n+1}(A \times \tilde E_{\Delta}) = \Q^\nu_{n+1}(A), \quad A \in \tilde E^n_{\Delta}.
\end{equation}
Let us now define
\begin{align}
	&\mathcal A:= \{A \times  \tilde E_{\Delta}  \times  \tilde E_{\Delta}  \times ...: A \in \tilde E^n_{\Delta}, \,\,n \geq 0\}, \nonumber\\
	&\Q^\nu(A \times  \tilde E_{\Delta}  \times  \tilde E_{\Delta} \times ...) := \Q^\nu_{n}(A), \quad A \in \tilde E^n_{\Delta}, n \geq 0.\label{Sec:PDP_Qinfty}
\end{align}
By the Kolmogorov extension theorem for product spaces (see Theorem 1.1.10  in \cite{Stroock_Varadhan}), it follows that $\Q^\nu$ is $\sigma$-additive on $\mathcal A$. Then, collecting  \eqref{Sec:PDP_Pinfty}, \eqref{Sec:PDP_Qn} and  \eqref{Sec:PDP_Qinfty}, it is easy to see that the $\sigma$-additivity of $\Q^\nu$ on $\mathcal A$ implies the $\sigma$-additivity of $\P^\nu$ on $\bigcup_n \mathcal F_{T_n}$.

Finally, we need to show that
$$
\P^{\nu}\big|_{\mathcal F_{T}}= L^\nu_{T}\,\P\quad \forall\, T >0,
$$
or, equivalently, that
$$
\E[L^\nu_{T}\,\psi]= \E^\nu[\psi]\,\quad \forall \psi\,\,\mathcal F_T\textup{-measurable function}.
$$
To this end, fix $T>0$, and let $\psi$ be a $\mathcal F_{T \wedge T_n}$-measurable bounded  function. In particular,  $\psi$ is $\mathcal F_{T \wedge T_m}$-measurable, for every $m \geq n$.
Since by definition
$
\P^{\nu}\big|_{\mathcal F_{T_n}}= L^\nu_{T_n}\,\P,  n \in \N,
$
we have
\begin{align*}
	\E^{\nu}[\psi]&= \E[L^\nu_{T_m}\,\psi]\\
	&= \E[\E[L^\nu_{T_m}\,\psi|\mathcal F_{T \wedge T_m}]]\\
	&= \E[\psi\,\E[L^\nu_{T_m}|\mathcal F_{T \wedge T_m}]]\\
	&= \E[\psi\,L^\nu_{T \wedge T_m}]\quad \forall m \geq n.
\end{align*}
Since $L^\nu_{T \wedge T_m} \underset{m \rightarrow  \infty}{\longrightarrow} L_T^\nu$ a.s., and $(L^\nu_s)_{s \in [0,\,T]}$ is a uniformly integrable martingale, by Theorem 3.1, Chapter II, in \cite{Re-Yor}, we get
\begin{align*}
	\E^{\nu}[\psi]= \lim_{m \rightarrow \infty}\E[L^\nu_{T \wedge T_m}\,\psi]= \E[L^\nu_{T}\,\psi],\quad \forall \psi\in \bigcup_n\mathcal F_{T\wedge T_n}.
\end{align*}
Then, by the monotone class theorem, recalling that  $\bigvee_n\mathcal F_{T\wedge T_n}=\mathcal F_{\bigvee_n\mathcal F_{T\wedge T_n}}$ (see, e.g.,  Corollary 3.5, point 6, in \cite{chineseBook}), we get
\begin{align*}
	\E^{\nu}[\psi]=  \E[L^\nu_{T}\,\psi],\quad \forall \psi\in \bigvee_n\mathcal F_{T\wedge T_n}=\mathcal F_{\bigvee_n\mathcal F_{T\wedge T_n}}= \mathcal F_T.
\end{align*}
This concludes the proof.
\endproof


Finally, for every $x\in E$, $a \in A$ and $\nu \in \mathcal V$, we introduce the dual  functional cost
\begin{equation}\label{Sec:PDP_dual_functional_cost}
	J(x,a,\nu) := \spernuxa{\int_{0}^{\infty} e^{-\delta\,t}\, f(X_t,I_t)\,dt},
\end{equation}
and the dual value function
\begin{equation}\label{Sec:PDP_dual_value_function}
	V^{\ast}(x,a) := \inf_{ \nu \in \mathcal{V}} J(x,a,\nu),
\end{equation}
where  $\delta >0$ in \eqref{Sec:PDP_dual_functional_cost} is  the discount factor introduced in Section \ref{Sec:PDP_Sec_control_problem}.


\section{Constrained BSDEs and  representation of the dual value function}\label{Sec:PDP_Sec_ConstrainedBSDE}
In this section we introduce a BSDE with a sign constrain on its martingale part, for which we prove the existence and uniqueness of a minimal solution, in an appropriate sense.
This constrained BSDE is then used to give a probabilistic representation formula for the dual value function introduced in \eqref{Sec:PDP_dual_value_function}.

Throughout this section we still  assume that {\bf (Hh$\lambda$Q)}, {\bf (H$\lambda_0$)} and {\bf (H$f$)} hold.
The random measures $p$, $\tilde p$ and $q$, as well as  the dual control setting $\Omega, \mathbb F, X, \P^{x,a}$, are the same as in Section \ref{Sec:PDP_Section_control_rand}.
We recall that $\mathbb{F}= (\mathcal{F}_t)_{t \geqslant 0}$ is the augmentation of the natural filtration generated by $p$, and that  $\mathcal{P}_T$, $T >0$, denotes   the $\sigma$-field of $\mathbb{F}$-predictable subsets of $[0,\,T] \times \Omega$.

For any $(x,a) \in E \times A$ we introduce the following notation.
\begin{itemize}
	\item $\textbf{L}^\textbf{2}_{\textbf{x,a}}(\mathcal{F}_\tau)$,  the set of $\mathcal{F}_\tau$-measurable random variables $\xi$ such that $\sperxa{|\xi|^2} < \infty$; here  $\tau \geqslant 0$ is an $\mathbb F$-stopping time.
	\item $\textbf{S}^\infty$  the set of real-valued c\`adl\`ag 	adapted processes $Y = (Y_t)_{t \geqslant 0}$ which are uniformly bounded.
	\item $\textbf{S}^\textbf{2}_{\textbf{x,a}}\textbf{(0,\,T)}$, $T >0$, the set of real-valued c\`adl\`ag adapted processes $Y = (Y_t)_{0 \leqslant t \leqslant T}$ satisfying
	\begin{displaymath}
		||Y||_{\textbf{S}^\textbf{2}_{\textbf{x,a}}\textbf{(0,\,T)}}:=\sperxa{\sup_{0 \leqslant t \leqslant T}|Y_t|^2} < \infty.
	\end{displaymath}
	\item $\textbf{L}_\textbf{x,a}^\textbf{2}(\textbf{0,\,T})$, $T >0$,  the set of real-valued progressive processes $\Phi = (\Phi_t)_{0 \leqslant t \leqslant T}$ such that
	\begin{displaymath}
		||Y||^2_{\textbf{L}^\textbf{2}_\textbf{x,a}\textbf{(0,\,T)}}:=\sperxa{\int_0^T |Y_t|^2\,dt} < \infty.
	\end{displaymath}
	We also define $\textbf{L}^\textbf{2}_{\textbf{x,a,\text{loc}}} :=\cap_{T >0} \textbf{L}_\textbf{x,a}^\textbf{2}(\textbf{0,\,T})$.
	\item $\textbf{L}_\textbf{x,a}^{\textbf{2}}(\textup{q};\textbf{0,\,T})$, $T >0$,  the set of $\mathcal{P}_T\otimes \mathcal{B}(E) \otimes \mathcal A$-measurable maps $Z: \Omega \times [0,\,T] \times E \times A \rightarrow \R$ such that
	\begin{align*}
		||Z||^2_{\textbf{L}_{\textbf{x,a}}^{\textbf{2}}(\textup{q};\textbf{0,\,T})}&:=\mathbb{E}^{x,a}\Big[ \int_{0}^{T}\int_{E \times A} |Z_t(y,b)|^2 \, \tilde{p}(dt\,dy\,db) \Big]\\
		&  = \mathbb{E}^{x,a}\Big[  \int_{0}^{T}\int_{E} |Z_t(y,I_t)|^2 \, \lambda(X_t,I_t)\,Q(X_t,I_t,dy)\,ds\\ 
		&+ \int_{0}^{T}\int_{A} |Z_t(X_t,b)|^2 \, \lambda_0(db)\,ds \Big] < \infty.
	\end{align*}
	We also define $\textbf{L}^{\textbf{2}}_{\textbf{x,a,\textup{loc}}}(\textup{q}) :=\cap_{T >0} \textbf{L}^{\textbf{2}}_{\textbf{x,a}}(\textup{q};\textbf{0,\,T})$.
	\item $\textbf{L}^{\textbf{2}}(\lambda_0)$,
	the set of $\mathcal{A}$-measurable maps $\psi: A \rightarrow \R$ such that
	\begin{eqnarray*}
		|\psi|^2_{\textbf{L}^{\textbf{2}}(\lambda_0)}:=\int_{A} |\psi(b)|^2 \, \lambda_0(db)< \infty.
	\end{eqnarray*}
	\item $\textbf{L}^{\textbf{2}}_\textbf{x,a}(\lambda_0;\textbf{0,\,T})$, $T >0$,  
	the set of $\mathcal{P}_T \otimes \mathcal A$-measurable maps $W: \Omega \times [0,\,T] \times A \rightarrow \R$ such that
	\begin{eqnarray*}
		|W|^2_{\textbf{L}^{\textbf{2}}_{\textbf{x,a}}(\lambda_0;\textbf{0,\,T})}:=
		\sperxa{\int_{0}^{T}\int_{A} |W_t(b)|^2 \, \lambda_0(db)\,ds}< \infty.
	\end{eqnarray*}
	We also define $\textbf{L}^{\textbf{2}}_{\textbf{x,a,\textup{loc}}}(\lambda_0) :=\cap_{T >0} \textbf{L}^{\textbf{2}}_\textbf{x,a}(\lambda_0;\textbf{0,\,T})$.
	\item $\textbf{K}^\textbf{2}_\textbf{x,a}(\textbf{0,\,T})$, $T>0$,  the set of nondecreasing c\`adl\`ag predictable  processes $K = (K_t)_{0 \leqslant t \leqslant T}$ such that $K_0 = 0$ and $\sperxa{|K_T|^2}< \infty$. We also define $\textbf{K}^\textbf{2}_\textbf{x,a,\text{loc}} :=\cap_{T >0} \textbf{K}^\textbf{2}_\textbf{x,a}(\textbf{0,\,T})$.
\end{itemize}
We are interested in studying the following family of BSDEs with partially nonnegative jumps over an infinite horizon, parametrized by $(x,a)$: 
${\mathbb{P}}^{x,a}$-a.s.,
\begin{align}
	Y^{x,a}_{s} &= Y^{x,a}_T - \delta \,\int_{s}^{T}Y^{x,a}_r\,dr + \int_{s}^{T}f(X_r,I_r)\,dr -( K^{x,a}_T - K^{x,a}_s) \label{Sec:PDP_BSDE}\\
	&\hspace{-10mm}- \int_{s}^{T}\int_{A}Z^{x,a}_r(X_r,\,b)\, \lambda_0(db)\,dr - \int_{s}^{T}\int_{E \times A}Z^{x,a}_r(y,\,b)\, q(dr\,dy\, db), \quad  0\leqslant s \leqslant T<\infty,\nonumber
\end{align}
with
\begin{equation}\label{Sec:PDP_BSDE_constraint}
	Z_s^{x,a}(X_{s-},b)\geqslant 0, \qquad  ds \otimes d\P^{x,a} \otimes \lambda_0(db), \text{-a.e. on } [0,\,\infty)\times \Omega \times A,
\end{equation}
where  $\delta$ is the positive  parameter introduced in Section \ref{Sec:PDP_Sec_control_problem}.

We look for a \emph{maximal solution} $(Y^{x,a},Z^{x,a},K^{x,a})\in \textbf{S}^{\infty}\times \textbf{L}^{\textbf{2}}_\textbf{x,a,\text{loc}}(\textup{q})\times \textbf{K}^\textbf{2}_\textbf{x,a,\text{loc}}$ to \eqref{Sec:PDP_BSDE}-\eqref{Sec:PDP_BSDE_constraint},
in the sense that for any other solution $(\tilde{Y}, \tilde{Z},\tilde{K})\in \textbf{S}^{\infty}\times \textbf{L}^{\textbf{2}}_\textbf{x,a,\text{loc}}(\textup{q})\times \textbf{K}^\textbf{2}_\textbf{x,a,\text{loc}}$ to \eqref{Sec:PDP_BSDE}-\eqref{Sec:PDP_BSDE_constraint}, we  have $Y_t^{x,a} \geqslant \tilde{Y}_t$, $\P^{x,a}$-a.s., for all $t\geqslant 0$.

\begin{proposition}\label{Sec:PDP_Prop_uniq_max_sol}
	Let Hypotheses  \textup{\textbf{(H$\textup{h$\lambda$Q}$)}}, \textup{\textbf{(H$\lambda_0$)}} and  \textup{\textbf{(H$\textup{f}$)}}  hold. Then, for any $(x,a) \in E \times A$, there exists at most one maximal solution $(Y^{x,a},Z^{x,a},K^{x,a})\in \textup{\textbf{S}}^{\infty}\times \textup{\textbf{L}}^{\textbf{2}}_\textbf{x,a,\textup{loc}}(\textup{q})\times \textup{\textbf{K}}^\textbf{2}_\textbf{x,a,\textup{loc}}$ to the BSDE with partially nonnegative jumps \eqref{Sec:PDP_BSDE}-\eqref{Sec:PDP_BSDE_constraint}.
\end{proposition}
\proof
Let $(Y,Z,K)$ and $(Y',Z',K')$ be two maximal solutions of \eqref{Sec:PDP_BSDE}-\eqref{Sec:PDP_BSDE_constraint}. 
By definition, we clearly have the uniqueness of the component $Y$.
Regarding the other components, taking the difference between the two backward equations we obtain: $\P^{x,a}$-a.s.
\begin{eqnarray*}
	0 &=& -(K_t -K'_t) -\int_{0}^{t} \int_{A} (Z_{s}(X_s,b)- Z'_{s}(X_s,b))\, \lambda_0(db)\,ds\nonumber\\
	&&\quad -\int_{0}^{t} \int_{E \times A} (Z_{s}(y,b)- Z'_{s}(y,b))\, q(ds\,dy\,db), \quad 0\leqslant t\leqslant T<\infty,
\end{eqnarray*}
that can be rewritten as
\begin{eqnarray}\label{Sec:PDP_E}
	&&\int_{0}^{t} \int_{E \times A} (Z_{s}(y,b)- Z'_{s}(y,b))\, p(ds\,dy\,db) = -(K_t -K'_t)\nonumber \\
	&&  +\int_{0}^{t} \int_{E} (Z_{s}(y,I_s)- Z'_{s}(y,I_s))\, \lambda(X_s,I_s)\,Q(X_s,I_s,dy)\,ds,\quad 0\leqslant t\leqslant T<\infty.
\end{eqnarray}
The right-hand side of \eqref{Sec:PDP_E} is a predictable process, therefore it has no totally inaccessible jumps (see, e.g., Proposition 2.24, Chapter I, in \cite{JS}); on the other hand, the left side is a pure jump process with totally inaccessible jumps.
This implies that $Z=Z'$, and as a consequence the component $K$ is unique as well.
\endproof
In the sequel, we prove by a penalization approach the existence of the maximal solution to \eqref{Sec:PDP_BSDE}-\eqref{Sec:PDP_BSDE_constraint}. In particular, this will  provide a probabilistic representation of the dual value function $V^\ast$ introduced in Section \ref{Sec:PDP_Section_dual_optimal_control}.

\subsection{Penalized BSDE and associated dual control problem}
Let us introduce the family of penalized BSDEs on $[0,\,\infty)$ associated to \eqref{Sec:PDP_BSDE}-\eqref{Sec:PDP_BSDE_constraint}, parametrized by the integer  $n \geqslant 1$: $\P^{x,a}$-a.s.,
\begin{eqnarray}\label{Sec:PDP_BSDE_penalized}
	Y_{s}^{n,x,a}&=& Y_{T}^{n,x,a}-\delta\int_s^T Y_r^{n,x,a}\,dr + \int_{s}^{T} f(X_{r},I_{r})\, dr\nonumber\\
	&& -n\int_s^T\int_A [Z_r^{n,x,a}(X_r,b)]^-\,\lambda_0(db)\,dr - \int_{s}^{T} \int_{A} Z_{r}^{n,x,a}(X_{r},b) \,\lambda_{0}(db)\, dr \nonumber\\
	&& - \int_{s}^{T} \int_{E \times A} Z_{r}^{n,x,a}(y,b) \, q(dr\,dy\,db),\quad 0\leqslant s\leqslant T<\infty,
\end{eqnarray}
where $[z]^- =\max(-z,0)$ denotes the negative part of $z$.

We shall prove that there exists a unique solution to equation \eqref{Sec:PDP_BSDE_penalized}, 
and provide an explicit representation to \eqref{Sec:PDP_BSDE_penalized} in terms of a family of dual control problems. To this end, we start by considering, for fixed  $T>0$,  the family of  BSDEs on $[0,\,T]$: $\P^{x,a}$-a.s.,
\begin{eqnarray}\label{Sec:PDP_BSDE_penalized_T}
	Y_{s}^{T,n,x,a}&=&  -\delta\int_s^T Y_r^{T,n,x,a}\,dr + \int_{s}^{T} f(X_{r},I_{r})\, dr\nonumber\\
	&& -n\int_s^T\int_A [Z_r^{T,n,x,a}(X_r,b)]^-\,\lambda_0(db)\,dr - \int_{s}^{T} \int_{A} Z_{r}^{T,n,x,a}(X_{r},b) \,\lambda_{0}(db)\, dr \nonumber\\
	&& - \int_{s}^{T} \int_{E \times A} Z_{r}^{T,n,x,a}(y,b) \, q(dr\,dy\,db),\quad 0\leqslant s\leqslant T,
\end{eqnarray}
with zero  final cost at time $T>0$.

\begin{remark}\label{f^n_lipsch_in_z}
	The penalized BSDE \eqref{Sec:PDP_BSDE_penalized_T} can be rewritten in the equivalent form: $\P^{x,a}$-a.s.,
	$$
	Y_{s}^{T,n,x,a}=\int_{s}^{T} f^n(X_{r},\,I_{r},\,Y_r^{T,n,x,a},\,Z_r^{T,n,x,a})\, ds  - \int_{s}^{T} \int_{E \times A} Z_{r}^{T,n,x,a}(y,b) \, q(dr\,dy\,db),
	$$
	$s \in [0,\,T]$, where the generator $f^n$ is defined by
	\begin{equation}\label{Sec:PDP_f^n}
		f^n(x,a,u,\psi) := f(x,a)  -\delta u -\int_{A}\left\{ n\,[\psi(a)]^{-} + \psi(b)\right\} \,\lambda_{0}(db),
	\end{equation}
	for all $(x,a,u,\psi)\in  E \times A \times \R\times \textbf{L}^{\textbf{2}}(\lambda_0)$.
	
	We notice that, under  Hypotheses \textup{\textbf{(H$\textup{h$\lambda$Q}$)}}, \textup{\textbf{(H$\lambda_0$)}} and  \textup{\textbf{(H$\textup{f}$)}}, $f^n$ is Lipschitz continuous in $\psi$ with respect to the norm of $\textbf{L}^{\textbf{2}}(\lambda_0)$,  uniformly in $(x,a,u)$, i.e., for every $n \in \N$, there exists a constant $L_n$, depending only on $n$, such that for every $(x,a,u)\in  E \times A \times \R$ and $\psi,\,\psi' \in \textbf{L}^{\textbf{2}}(\lambda_0)$,
	\begin{displaymath}
		|f^n(x,a,u,\psi')-f^n(x,a,u,\psi)| \leqslant L^n |\psi-\psi'|_{\textbf{L}^{\textbf{2}}(\lambda_0)}.
	\end{displaymath}
\end{remark}
For every integer $n \geqslant 1$, let $\mathcal{V}^n$ denote the subset of elements $\nu \in \mathcal{V}$ valued in $(0,\,n]$.
\begin{proposition}\label{Sec:PDP_Prop_ex_uniq_BSDE_pen_T}
	Let Hypotheses \textup{\textbf{(H$\textup{h$\lambda$Q}$)}}, \textup{\textbf{(H$\lambda_0$)}} and  \textup{\textbf{(H$\textup{f}$)}} hold.
	For every  $(x,a,n,T) \in E \times A\times \N\times (0,\,\infty)$,
	there exists a unique solution $(Y^{T,n,x,a},Z^{T,n,x,a}) \in \textup{\textbf{S}}^\infty
	\times \textup{\textbf{L}}^{\textbf{2}}_\textbf{x,a}(\textup{q};\textbf{\textup{0},\,\textup{T}})$ to \eqref{Sec:PDP_BSDE_penalized_T}. 
	Moreover, 
	the following uniform estimate  holds: $\P^{x,a}$-a.s.,
	\begin{equation}\label{Sec:PDP_Sinfty_estimate}
		Y^{T,n,x,a}_s \leqslant   \frac{M_f}{\delta}, \quad \forall \,s \in [0,\,T].
	\end{equation}
\end{proposition}
\proof
The existence and uniqueness of a solution $(Y^{T,n,x,a},Z^{T,n,x,a}) \in \textup{\textbf{S}}^{\textbf{2}}_\textbf{x,a}(\textbf{\textup{0},\,\textup{T}})
\times \textup{\textbf{L}}^{\textbf{2}}_\textbf{x,a}(\textup{q};\textbf{\textup{0},\,\textup{T}})$ to \eqref{Sec:PDP_BSDE_penalized_T} is based on a fixed point argument, and uses  integral   representation results for $\mathbb F$-martingales, with $\mathbb F$ the natural filtration (see, e.g., Theorem 5.4 in \cite{J}).
This procedure  is standard and   we  omit it  (similar proofs can be found in the proofs of Theorem 3.2 in \cite{Xia},  Proposition 3.2 in \cite{Bech}, Theorem 3.4 in \cite{CoFu-m} ).
It remains to  prove the uniform estimate \eqref{Sec:PDP_Sinfty_estimate}.
To this end, let us apply It\^o's formula to $e^{-\delta\,r}\,Y^{x,a,n,T}_r$ between $s$ and $T$. We get:
$\P^{x,a}$-a.s.
\begin{align}
	& Y^{T,n,x,a}_{s} =   \int_{s}^{T} e^{-\delta\,(r-s)}\,f(X_{r},I_{r})\, dr - \int_{s}^{T}\,\int_{E \times A}  e^{-\delta\,(r-s)}\,Z^{T,n,x,a}_{r}(y,b)\, q(dr\,dy\,db)\nonumber\\
	&- \int_{s}^{T} \int_{A} e^{-\delta\,(r-s)}\,\{n [Z^{T,n,x,a}_{r}(X_r,b)]^- + Z^{T,n,x,a}_{r}(X_r,b)\}\, \lambda_0(db)\, dr, \quad s \in [0,\,T].\label{Sec:PDP_BSDE_proof1}
\end{align}
Now for any $\nu \in \mathcal{V}^n$, let us introduce the compensated martingale measure $q^\nu(ds \, dy \,db)= q(ds \, dy \,db)- (\nu_s(b)-1)\,d_1(s,y,b)\,\tilde{p}(ds\,dy\,db)$ under $\P^{x,a}_\nu$.
Taking  the expectation in \eqref{Sec:PDP_BSDE_proof1} under $\P^{x,a}_\nu$,
conditional to $\mathcal{F}_s$, and
since   $Z^{T,n,x,a}$ is in  $\textbf{L}^{\textbf{2}}_\textbf{x,a}(\textup{q};\textup{\textbf{0}},\textup{\textbf{T}})$, from Lemma \ref{Sec:PDP_lemma_P_nu_martingale} we get  that, $\P^{x,a}$-a.s.,
\begin{align}
	Y^{T,n,x,a}_{s}
	&= - \spernuxa{\int_{s}^{T} \int_{A} e^{-\delta\,(r-s)}\,\{n [Z^{T,n,x,a}_{r}(X_r,b)]^- + \nu_r(b)\,Z^{T,n,x,a}_{r}(X_r,b)\}\, \lambda_0(db)\, dr\Big | \mathcal{F}_s}\nonumber \\
	&+ \spernuxa{\int_{s}^{T}
		e^{-\delta\,(r-s)}\,f(X_{r},I_{r})\, dr\Big| \mathcal{F}_s}, \quad s \in [0,\,T]. \label{Sec:PDP_BSDE_nu}
\end{align}
From the elementary numerical inequality: $n[z]^- + \nu z \geqslant 0$ for all $z \in \R$, $\nu \in (0,\,n]$, we deduce by \eqref{Sec:PDP_BSDE_nu} that, for all $\nu\in\mathcal V^n$,
\begin{eqnarray*}
	Y^{T,n,x,a}_s \leqslant \spernuxa{\int_{s}^{T} e^{-\delta\,(r-s)}\,f(X_{r},I_{r})\, dr\Big| \mathcal{F}_s},\quad s \in [0,\,T].
\end{eqnarray*}
Therefore, $\P^{x,a}$-a.s.,
\begin{eqnarray*}
	Y^{T,n,x,a}_s 
	\leqslant   \spernuxa{\int_{s}^{\infty} e^{-\delta\,(r-s)}\,|f(X_{r},I_{r})|\, dr\Big| \mathcal{F}_s}
	\leqslant    \frac{M_f}{\delta}, \quad s \in [0,\,T].
\end{eqnarray*}
\endproof
\begin{proposition}\label{Sec:PDP_Prop_exist_uniq_pen_BSDE}
	Let Hypotheses \textup{\textbf{(H$\textup{h$\lambda$Q}$)}}, \textup{\textbf{(H$\lambda_0$)}} and  \textup{\textbf{(H$\textup{f}$)}} hold.
	Then, for every  $(x,a,n) \in E \times A\times \N$, there exists a unique solution $(Y^{n,x,a},Z^{n,x,a}) \in \textup{\textbf{S}}^\infty \times \textup{\textbf{L}}^{\textbf{2}}_\textbf{x,a,\textup{loc}}(\textup{q})$ to the penalized  BSDE \eqref{Sec:PDP_BSDE_penalized}. 
\end{proposition}
\proof
\emph{Uniqueness.} Fix $n\in  \N$, $(x,a)\in E \times A$, and consider two solutions $(Y^{1},Z^{1})=(Y^{1,n,x,a},Z^{1,n,x,a})$, $(Y^{2},Z^{2})= (Y^{2,n,x,a},Z^{2,n,x,a})\in \textup{\textbf{S}}^\infty \times \textup{\textbf{L}}^{\textbf{2}}_\textbf{x,a,\textup{loc}}(\textup{q})$ of \eqref{Sec:PDP_BSDE_penalized}. 
Set $\bar{Y} = Y^{2}-Y^{1}$, $\bar{Z} = Z^{2}-Z^{1}$. Let $0\leqslant s \leqslant T< \infty$. Then, an application of It\^o's formula to $e^{-2\,\delta\,r}|\bar{Y}_r|^2$ between $s$ and $T$ yields: $\P^{x,a}$-a.s.,
\begin{align}\label{Sec:PDP_eq_uniq}
	e^{-2\,\delta\,s}|\bar{Y}_s|^2 &= e^{-2\,\delta\,T}|\bar{Y}_T|^2 \nonumber\\
	&-2n \int_s^T\int_{A}e^{-2\,\delta\,r}\,\bar{Y}_r\,\{[Z_r^2(X_s,b)]^--[Z_r^1(X_s,b)]^-\}\,\lambda_0(db)\,dr\nonumber\\
	&-2 \int_s^T\int_{A}e^{-2\,\delta\,r}\,\bar{Y}_r\,\bar{Z}_r(X_s,b)\,\lambda_0(db)\,dr\nonumber\\
	&-2 \int_s^T\int_{E \times A}e^{-2\,\delta\,r}\,\bar{Y}_r\,\bar{Z}_r(y,b)\,q(dr\,dy\,db)\nonumber\\
	&- \int_s^T\int_{E \times A}e^{-2\,\delta\,r}\,|\bar{Z}_r(y,b)|^2\,p(dr\,dy\,db).
\end{align}
Notice that
\begin{eqnarray*}
	&&-n\,\int_{s}^{T} \int_{A} e^{-\delta\,(r-s)}\,\bar{Y}_r\,\{[Z_r^2(X_r,b)]^--[Z_r^1(X_r,b)]^-\}\,\lambda_0(db)\,dr \\
	&&= \int_{s}^{T} \int_{A} e^{-\delta\,(r-s)}\,\bar{Y}_r\,\{Z_r^2(X_r,b)-Z_r^1(X_r,b)\}\,\nu_r^{\varepsilon}\,\lambda_0(db)\,dr\\
	&& \,\,  -\varepsilon\,\int_{s}^{T} \int_{A} e^{-\delta\,(r-s)}\,\bar{Y}_r\,\{Z_r^2(X_r,b)-Z_r^1(X_r,b)\}\,\one_{\{|\bar{Y}_r|\leqslant 1\}}\cdot\\
	&& \hspace{40mm} \cdot\, \one_{\{[Z_r^2(X_r,b)]^-=[Z_r^1(X_r,b)]^- ,\,|\bar{Z}_r(X_r,b)|\leqslant 1\}}\,\lambda_0(db)\, dr\\
	&& \,\,  -\varepsilon\,\int_{s}^{T} \int_{A} e^{-\delta\,(r-s)}\,\one_{\{|\bar{Y}_r|> 1\}}\, \one_{\{[Z_r^2(X_r,b)]^-=[Z_r^1(X_r,b)]^- ,\,|\bar{Z}_r(X_r,b)|>1 \}}\,\lambda_0(db)\, dr,
\end{eqnarray*}
where $\nu^{\varepsilon}: \R_{+}\times \Omega \times A$ is given by
\begin{eqnarray}
	\nu^\varepsilon_r(b) &=& -n\,\frac{[Z_r^2(X_r,b)]^--[Z_r^1(X_r,b)]^-}{\bar{Z}_r(X_r,b)}\one_{\{Z_r^2(X_r,b)]^--[Z_r^1(X_r,b)]^- \neq 0 \}}\label{Sec:PDP_nuepsilon}\\
	&& + \varepsilon\,\one_{\{|\bar{Y}_r|\leqslant 1\}}\,\one_{\{[Z_r^2(X_r,b)]^-=[Z_r^1(X_r,b)]^- ,\,|\bar{Z}_r(X_r,b)|\leqslant 1\}}\nonumber \\
	&&+ \varepsilon\,(\bar{Y}_r)^{-1}\,(\bar{Z}_r(X_r^{x,a},b))^{-1}\,\one_{\{|\bar{Y}_r|> 1\}}\, \one_{\{[Z_r^2(X_r^{x,a},b)]^-=[Z_r^1(X_s^{x,a},b)]^- ,\,|\bar{Z}_r(X_r^{x,a},b)|>1 \}},\nonumber
\end{eqnarray}
for arbitrary $\varepsilon \in (0,\,1)$.
In particular, $\nu^{\varepsilon}$ is a $\mathcal{P}\otimes \mathcal{A}$-measurable map satisfying $\nu_r^{\varepsilon}(b) \in [\varepsilon,\,n]$, $dr \otimes d\P^{x,a} \otimes \lambda_0(db)$-almost everywhere.
Consider the probability measure $\P^{x,a}_{\nu^\varepsilon}$ on $(\Omega, \mathcal F_{\infty})$, whose restriction to $(\Omega, \mathcal F_T)$ has 
Radon-Nikodym 
density:
\begin{equation}\label{Sec:PDP_RN_density}
	L_s^{\nu^\varepsilon} := \mathcal{E}\left(\int_{0}^{\cdot}\int_{E \times A}(\nu^\varepsilon_t(b) \, d_1(t,y,b) \, + d_2(t,y,b) - 1 )\,q(dt\,dy\,db)\right)_s
\end{equation}
for all $0 \leqslant s \leqslant T$, where $\mathcal{E}(\cdot)_s$ is the Dol\'eans-Dade exponential. 
The existence of such a probability is guaranteed by Proposition \ref{P_Prob_infty}.
From Lemma 
\ref{Sec:PDP_lemma_P_nu_martingale}
it follows that $(L_s^{\nu^\varepsilon})_{s \in [0,\,T]}$ is a uniformly integrable martingale. 
Moreover, $L_T^{\nu^\varepsilon}\in \textup{\textbf{L}}^\textup{\textbf{p}}(\mathcal{F}_T)$, for any $p \geqslant 1$. Under the probability measure $\P^{x,a}_{\nu^{\varepsilon}}$, by Girsanov's theorem, the compensator of $p$ on $[0,\,T]\times E\times A$ is $(\nu_s^{\varepsilon}(b)\,d_1(s,y,b)$ $+ d_2(s,y,b))\,\tilde{p}(ds\,dy\,db)$. We denote by $q^{\nu^{\varepsilon}}(ds\,dy\,db)$ $:= p(ds\,dy\,db) - (\nu_s^{\varepsilon}(b)\,d_1(s,y,b)+ d_2(s,y,b))$ $\tilde{p}(ds\,dy\,db)$ the compensated martingale measure of $p$ under $\P^{x,a}_{\nu^{\varepsilon}}$.
Therefore equation \eqref{Sec:PDP_eq_uniq} becomes: $\P^{x,a}$-a.s.,
$$
e^{-2\,\delta\,s}|\bar{Y}_s|^2 \leqslant e^{-2\,\delta\,T}|\bar{Y}_T|^2
-2 \int_s^T\int_{A}e^{-2\,\delta\,r}\,\bar{Y}_r\,\bar{Z}_r(X_s,b)\,q^{\nu^{\varepsilon}}(ds\,dy\,db) +2\,\frac{\varepsilon}{\delta}\,\lambda_0(A),
$$
for all $\varepsilon \in (0,\,1)$. Moreover, from the arbitrariness of $\varepsilon$, we obtain
\begin{equation}\label{Sec:PDP_ineq_subm}
	e^{-2\,\delta\,s}|\bar{Y}_s|^2 \leqslant e^{-2\,\delta\,T}|\bar{Y}_T|^2
	-2 \int_s^T\int_{A}e^{-2\,\delta\,r}\,\bar{Y}_r\,\bar{Z}_r(X_s,b)\,q^{\nu^{\varepsilon}}(ds\,dy\,db).
\end{equation}
From Lemma \ref{Sec:PDP_lemma_P_nu_martingale}, we see that the stochastic integral in \eqref{Sec:PDP_ineq_subm} is a martingale, so that, taking the expectation $\E^{x,a}_{\nu^{\varepsilon}}$, conditional on $\mathcal{F}_s$, with respect to $\P^{x,a}_{\nu^{\varepsilon}}$, we achieve
\begin{eqnarray}\label{Sec:PDP_estimate_difference}
	e^{-2\,\delta\,s}|\bar{Y}_s|^2 &\leqslant& e^{-2\,\delta\,T}\,\E^{x,a}_{\nu^{\varepsilon}}[|\bar{Y}_T|^2|\mathcal{F}_s].
\end{eqnarray}
In particular, $(e^{-2\,\delta\,s}|\bar{Y}_s|^2)_{t\geqslant 0}$ is a submartingale. Since $\bar{Y}$ is uniformly bounded, we see that $(e^{-2\,\delta\,s}|\bar{Y}_s|^2)_{t\geqslant 0}$ is a uniformly integrable submartingale, therefore $e^{-2\,\delta\,s}|\bar{Y}_s|^2\rightarrow \xi_{\infty}\in \textup{\textbf{L}}^{\textbf{1}}(\Omega,\mathcal{F},\P^{x,a}_{\nu^{\varepsilon}})$, as $s \rightarrow \infty$. Using again the boundedness of $\bar{Y}$, we obtain that $\xi_{\infty}=0$, which implies $\bar{Y}=0$. Finally, plugging $\bar{Y}=0$ into \eqref{Sec:PDP_eq_uniq} we conclude that $\bar{Z} = 0$.

\medskip

\noindent\emph{Existence.}
Fix $(x,a,n)\in E \times A  \times \N$.
For $T>0$, let $(Y^{T,n,x,a},Z^{T,n,x,a})=(Y^{T},Z^{T})$ denote the unique solution to the penalized BSDE \eqref{Sec:PDP_BSDE_penalized_T} on $[0,\,T]$.\\
\emph{Step 1.} \emph{Convergence of $(Y^T)_T$.}
Let $T,T'>0$, with $T< T'$, and $s \in [0,\,T]$. We have
\begin{eqnarray}\label{Sec:PDP_Ybar_conv}
	|Y_s^{T'}-Y_s^{T}|^2
	\leqslant e^{-2\,\delta\,(T-s)}\,\spernuxaepsilon{|Y_T^{T'}-Y_T^{T}|^2|\mathcal{F}_s}\overset{T\rightarrow \infty}{\longrightarrow} 0,
\end{eqnarray}
where the convergence result follows from \eqref{Sec:PDP_Sinfty_estimate}.
Let us now consider the sequence of real-valued c\`adl\`ag adapted processes $(Y^T)_T$. It follows from \eqref{Sec:PDP_Ybar_conv} that, for any $t\geqslant 0$, the sequence $(Y_t^T(\omega))_T$ is  Cauchy  for almost every $\omega$, so that it converges $\P^{x,a}$-a.s. to some $\mathcal{F}_t$-measurable random variable $Y_t$, which is bounded from the right-hand side of \eqref{Sec:PDP_Sinfty_estimate}. Moreover, using again \eqref{Sec:PDP_Ybar_conv} and \eqref{Sec:PDP_Sinfty_estimate}, we see that, for any $0\leqslant S< T\wedge T'$, with $T,T' >0$, we have
\begin{equation}\label{Sec:PDP_conv_Y}
	\sup_{0\leqslant t\leqslant S}|Y_t^{T'}-Y_t^{T}| \leqslant e^{-\delta\,(T\wedge T' -S)}\,\frac{M_f}{\delta}\overset{T,T'\rightarrow \infty}{\longrightarrow} 0.
\end{equation}
In other words, the sequence $(Y^{T})_{T >0}$ converges $\P^{x,a}$-a.s. to $Y$ uniformly on compact subsets of $\R_{+}$. Since each $Y^T$ is a c\`adl\`ag process, it follows that $Y$ is c\`adl\`ag, as well. Finally, from estimate \eqref{Sec:PDP_Sinfty_estimate} we see that  $Y$ is uniformly bounded and therefore belongs to $\textup{\textbf{S}}^\infty$.

\medskip

\noindent \emph{Step 2. Convergence of $(Z^T)_T$.}
Let 
$S,\,T,\,T'>0$, with $S<T<T'$. Then, applying It\'o's formula to $e^{-2\,\delta\,r}\,|Y^{T'}_r-Y^{T}_r|^2$ between $0$ and $S$, and taking the expectation, we find
\begin{eqnarray*}
	&&\sperxa{\int_0^S\int_{E \times A}e^{-2\,\delta\,r}\,|Z^{T'}_r(y,b)-Z^{T}_r(y,b)|^2\,\tilde{p}(dr\,dy\,db)} \\
	&&  = e^{-2\,\delta\,S}\sperxa{|Y^{T'}_S-Y^{T}_S|^2} - |Y^{T'}_0-Y^{T}_0|^2 \\
	&& -2n\, \sperxa{\int_0^S\,\int_{A}e^{-2\,\delta\,r}\,(Y^{T'}_r-Y^{T}_r)\,\{[Z_r^2(X_r,b)]^--[Z_r^1(X_r,b)]^-\}\,\lambda_0(db)\,dr}\nonumber\\
	&&   -2 \sperxa{\int_0^S\int_{A}e^{-2\,\delta\,r}\,(Y^{T'}_r-Y^{T}_r)\,(Z^{T'}_r(X_r,b)-Z^{T}_r(X_r,b))\,\lambda_0(db)\,dr}.\nonumber\\
\end{eqnarray*}
Recalling the elementary inequality $bc \leqslant b^2+c^2/4$, for any $b,\,c\,\in \R$, we get
\begin{eqnarray*}
	&&\sperxa{\int_0^S\int_{E \times A}e^{-2\,\delta\,r}\,|Z^{T'}_r(y,b)-Z^{T}_r(y,b)|^2\,\tilde{p}(dr\,dy\,db)} \\
	&&  \leqslant e^{-2\,\delta\,S}\sperxa{|Y^{T'}_S-Y^{T}_S|^2}  +4(n^2+1)\,\lambda_0(A)\, \sperxa{\int_0^S e^{-2\,\delta\,r}\,|Y^{T'}_r-Y^{T}_r|^2\,dr}\\
	&& +\frac{1}{4} \sperxa{\int_0^S\int_{A}e^{-2\,\delta\,r}\,|[Z_r^2(X_r,b)]^--[Z_r^1(X_r,b)]^-|^2\,\lambda_0(db)\,dr}\\
	&& +\frac{1}{4} \sperxa{\int_0^S\int_{A}e^{-2\,\delta\,r}\,|Z^{T'}_r(X_r,b)-Z^{T}_r(X_r,b)|^2\,\lambda_0(db)\,dr}.
\end{eqnarray*}
Multiplying the previous inequality by $e^{2\,\delta\,s}$, and recalling the form 
of the compensator $\tilde{p}$, we get
\begin{align*}
	&\frac{1}{2}\sperxa{\int_0^S\int_{E \times A}e^{-2\,\delta\,r}\,|Z^{T'}_r(y,b)-Z^{T}_r(y,b)|^2\,\tilde{p}(dr\,dy\,db)} \\
	&  \leqslant e^{-2\,\delta\,S}\sperxa{|Y^{T'}_S-Y^{T}_S|^2}  +4(n^2+1)\,\lambda_0(A) \,\sperxa{\int_0^S e^{-2\,\delta\,r}\,|Y^{T'}_r-Y^{T}_r|^2\,dr}\\
	&\overset{T,T'\rightarrow \infty}{\longrightarrow}0,
\end{align*}
where the convergence to zero follows from estimate \eqref{Sec:PDP_conv_Y}.
Then, for any $S>0$, we see  that $(Z^T_{|[0,\,S]})_{T>S}$ is a Cauchy sequence in the Hilbert space $\textbf{L}^{\textbf{2}}_\textbf{x,a}(\textup{q};\textbf{0},\,\textbf{S})$. Therefore, we deduce that there exists $\tilde{Z}^S \in \textbf{L}^{\textbf{2}}_\textbf{x,a}(\textup{q};\textbf{0},\,\textbf{S})$ such that $(Z^T_{|[0,\,S]})_{T>S}$ converges to  $\tilde{Z}^S$ in $\textbf{L}^{\textbf{2}}_\textbf{x,a}(\textup{q};\textbf{0},\,\textbf{S})$, i.e.,
\begin{displaymath}
	\sperxa{\int_0^S\int_{E \times A}e^{-2\,\delta\,r}\,|Z^{T}_r(y,b)-\tilde{Z}^{S}_r(y,b)|^2\,\tilde{p}(dr\,dy\,db)}\overset{T\rightarrow \infty}{\longrightarrow} 0.
\end{displaymath}
Notice that $\tilde{Z}^{S'}_{|[0,\,S]} = \tilde{Z}^S$, for any $0 \leqslant S \leqslant S'< \infty$. Indeed, $\tilde{Z}^{S'}_{|[0,\,S]}$, as $\tilde{Z}^S$, is the limit in $\textbf{L}^{\textbf{2}}_\textbf{x,a}(\textup{q};\textbf{0},\,\textbf{S})$ of $(Z^{T}_{|[0,\,S]})_{T>S}$. Hence, we define $Z_s = \tilde{Z}_s^S$ for all $s \in [0,\,S]$ and for any $S >0$. Observe that $Z\in \textbf{L}^{\textbf{2}}_\textbf{x,a,\textup{loc}}(\textup{q})$. Moreover, for any $S>0$, $(Z^{T}_{|[0,\,S]})_{T>S}$ converges to $Z_{|[0,\,S]}$ in $\textbf{L}^{\textbf{2}}_\textbf{x,a}(\textup{q};\textbf{0},\,\textbf{S})$, i.e.,
\begin{equation}\label{Sec:PDP_conv_Z}
	\sperxa{\int_0^S\int_{E \times A}e^{-2\,\delta\,r}\,|Z^{T}_r(y,b)-Z_r(y,b)|^2\,\tilde{p}(dr\,dy\,db)}\overset{T\rightarrow \infty}{\longrightarrow} 0.
\end{equation}
Now, fix $S \in [0,\,T]$ and consider the BSDE satisfied by $(Y^T,Z^T)$ on $[0,\,S]$: $\P^{x,a}$-a.s.,
\begin{eqnarray*}
	Y_{t}^{T}&=& Y_{S}^{T}-\delta \int_t^S Y_r^{T}\,dr + \int_{t}^{S} f(X_{r},I_{r})\, dr \nonumber\\
	&& -n\int_t^S\int_A [Z_r^{T}(X_r,b)]^-\,\lambda_0(db)\,dr - \int_{t}^{S} \int_{A} Z_{r}^{T}(X_{r},b) \,\lambda_{0}(db)\, dr, \nonumber\\
	&& - \int_{t}^{S} \int_{E \times A} Z_{r}^{T}(y,b) \, q(dr\,dy\,db),\quad \quad 0 \leqslant t\leqslant S.
\end{eqnarray*}
From \eqref{Sec:PDP_conv_Z} and \eqref{Sec:PDP_conv_Y}, we can pass to the limit in the above BSDE by letting $T \rightarrow \infty$ keeping $S$ fixed. Then we deduce that $(Y,Z)$ solves the penalized BSDE \eqref{Sec:PDP_BSDE_penalized} on $[0,\,S]$. Since $S$ is arbitrary, it follows that $(Y,Z)$ solves equation \eqref{Sec:PDP_BSDE_penalized} on $[0,\,\infty)$.
\endproof
The penalized  BSDE \eqref{Sec:PDP_BSDE_penalized} can be represented by means of  an suitable  family of  dual control problems.
\begin{lemma}\label{Sec:PDP_Lemma_rep_Y_n}
	Let Hypotheses \textup{\textbf{(H$\textup{h$\lambda$Q}$)}}, \textup{\textbf{(H$\lambda_0$)}} and  \textup{\textbf{(H$\textup{f}$)}} hold.
	Then, for every $(x,a,n) \in E\times A \times \N$, $\P^{x,a}$-a.s., the solution $(Y^{n,x,a},Z^{n,x,a})$
	to \eqref{Sec:PDP_BSDE_penalized} 
	admits the following explicit representation:
	\begin{equation}\label{Sec:PDP_rep_Y_n}
		Y_s^{n,x,a} = \essinf_{\nu \in \mathcal{V}^n}\spernuxa{\int_s^\infty e^{-\delta\,(r-s)}\,f(X_{r},I_{r})\, dr\Big| \mathcal{F}_s}, \quad s \geqslant 0.
	\end{equation}
\end{lemma}
\proof
Fix $n \in \N$, and for any $\nu \in \mathcal{V}^n$, let us introduce the compensated martingale measure $q^\nu(ds \, dy \,db)= q(ds \, dy \,db)- (\nu_s(b)-1)d_1(s,y,b)\tilde{p}(ds\,dy\,db)$ under $\P^{x,a}_\nu$.  Fix $T\geqslant s$ and apply It\^o's formula to $e^{-\delta\,r}\,Y_r^{n,x,a}$ between $s$ and $T$. Then we obtain:
\begin{eqnarray}\label{Sec:PDP_BSDE_proof}
	Y_{s}^{n,x,a} &=&  e^{-\delta\,(T-s)}\,Y_{T}^{n,x,a} + \int_{s}^{T} e^{-\delta\,(r-s)}\,f(X_{r},I_{r})\,dr \nonumber\\
	&&- \int_{s}^{T} \int_{A} e^{-\delta\,(r-s)}\,\{n [Z_{r}^{n,x,a}(X_r,b)]^- + \nu_r(a) \,Z_{r}^{n,x,a}(X_r,b)\}\, \lambda_0(db)\, dr \nonumber\\
	&& - \int_{s}^{T}\,\int_{E \times A}  e^{-\delta\,(r-s)}\,Z_{r}^{n,x,a}(y,b)\, q^\nu(dr\,dy\,db), \quad s \in [t,\,T].
\end{eqnarray}
Taking the expectation in \eqref{Sec:PDP_BSDE_proof} under $\P^{x,a}_\nu$,
conditional to $\mathcal{F}_s$, and since by Proposition \ref{Sec:PDP_Prop_exist_uniq_pen_BSDE} $Z^{n,x,a}$ is in $\textbf{L}^{\textbf{2}}_\textbf{\textup{loc},x,a}(\textup{q})$, we get from Lemma \ref{Sec:PDP_lemma_P_nu_martingale} that, $\P^{x,a}$-a.s.,
\begin{align}
	Y_{s}^{n,x,a} &= \spernuxa{e^{-\delta\,(T-s)}\,Y_{T}^{n,x,a} + \int_{s}^{T} e^{-\delta\,(r-s)}\,f(X_{r},I_{r})\, dr \Big | \mathcal{F}_s}\label{Sec:PDP_BSDE_spernu}\\
	&- \spernuxa{\int_{s}^{T} \int_{A} e^{-\delta\,(r-s)}\,\{n [Z_{r}^{n,x,a}(X_r,b)]^- + \nu_r(a)\,Z_{r}^{n,x,a}(X_r,b)\}\, \lambda_0(db)\, dr\Big | \mathcal{F}_s}.\nonumber
\end{align}
From the elementary numerical inequality: $n[z]^- + \nu z \geqslant 0$ for all $z \in \R$, $\nu \in (0,\,n]$, we deduce by \eqref{Sec:PDP_BSDE_spernu} that, for all $\nu\in\mathcal V^n$,
\begin{eqnarray*}
	Y_s^{n,x,a} &\leqslant& \spernuxa{e^{-\delta\,(T-s)}\,Y_{T}^{n,x,a} + \int_{s}^{T} e^{-\delta\,(r-s)}\,f(X_{r},I_{r})\, dr \Big | \mathcal{F}_s}\nonumber\\
	&\leqslant& \spernuxa{e^{-\delta\,(T-s)}\,Y_{T}^{n,x,a} + \int_{s}^{\infty} e^{-\delta\,(r-s)}\,f(X_{r},I_{r})\, dr \Big | \mathcal{F}_s}.
\end{eqnarray*}
Since $Y^{n,x,a}$ is in $\textbf{S}^{\infty}$ by Proposition \ref{Sec:PDP_Prop_exist_uniq_pen_BSDE}, sending $T \rightarrow \infty$, we obtain from the conditional version of  Lebesgue dominated convergence theorem that
\begin{eqnarray*}
	Y_s^{n,x,a} &\leqslant&  \spernuxa{\int_{s}^{\infty} e^{-\delta\,(r-s)}\,f(X_{r},I_{r})\, dr \Big | \mathcal{F}_s},
\end{eqnarray*}
for all $\nu\in\mathcal V^n$. Therefore,
\begin{eqnarray}
	Y_s^{n,x,a} \leqslant  \essinf_{\nu \in \mathcal{V}^n} \spernuxa{\int_{s}^{\infty} e^{-\delta\,(r-s)}\,f(X_{r},I_{r})\, dr \Big | \mathcal{F}_s}.
	\label{Sec:PDP_ineq_esssup}
\end{eqnarray}
On the other hand, for $\varepsilon \in (0,\,1)$, let us consider the process $\nu^\varepsilon \in \mathcal{V}^n$ defined by:
\begin{align*}
	\nu^\varepsilon_s(b) =& n\,\one_{\{Z_s^{n,x,an}(X_{s-},b) \leqslant 0 \}} + \varepsilon\,\one_{\{0 < Z_s^{n,x,a}(X_{s-},b) < 1\}} \\
	&+ \varepsilon\,Z_s^{n,x,a}(X_{s-},b)^{-1}\,\one_{\{Z_s^{n,x,a}(X_{s-},b) \geqslant 1\}}
\end{align*}
(notice that we can not take $\nu_s(b) = n \one_{\{Z_s^{n}(X_{s-},b)\leqslant 0\}}$, since this process does not belong to $\mathcal{V}^n$ because of the requirement of strict positivity).
By construction, we have
\begin{displaymath}
	n [Z_{s}^{n}(X_{s-},b)]^-+ \nu^{\varepsilon}_s(b)\,Z_{s}^{n}(X_{s-},b) \leqslant \varepsilon, \qquad s \geqslant 0,\,b \in A,
\end{displaymath}
and thus for this choice of $\nu = \nu^\varepsilon$ in \eqref{Sec:PDP_BSDE_spernu}:
\begin{align*}
	Y_{s}^{n,x,a} &\geqslant \spernuxaepsilon{e^{-\delta\,(T-s)}\,Y_{T}^{n,x,a} + \int_{s}^{T} e^{-\delta\,(r-s)}\,f(X_{r},I_{r})\, dr \Big | \mathcal{F}_s}\\
	&-\varepsilon\, \frac{1-e^{-\delta (T-s)}}{\delta}\,\lambda_0(A).
\end{align*}
Letting $T \rightarrow \infty$, since $f$ is bounded by $M_f$ and  $Y^{n,x,a}$ is in $\textbf{S}^{\infty}$, it follows from the conditional version of Lebesgue dominated convergence theorem that
\begin{eqnarray*}
	Y_{s}^{n,x,a} &\geqslant&  \spernuxaepsilon{\int_{s}^{\infty} e^{-\delta\,(r-s)}\,f(X_{r},I_{r})\, dr\Big | \mathcal{F}_s}-\frac{\varepsilon}{\delta}\,\lambda_0(A),\\
	&\geqslant&  \essinf_{\nu \in \mathcal{V}^n} \spernuxa{\int_{s}^{\infty} e^{-\delta\,(r-s)}\,f(X_{r},I_{r})\, dr \Big | \mathcal{F}_s}-\frac{\varepsilon}{\delta}\,\lambda_0(A).
\end{eqnarray*}
From the arbitrariness of $\varepsilon$, together with \eqref{Sec:PDP_ineq_esssup}, this is enough to prove the required representation of $Y^{n,x,a}$.
\endproof
%
Let us define
\begin{displaymath}
	K_t^{n,x,a} := n\int_{0}^{t}\int_{A}[Z_{s}^{n,x,a}(X_s,b)]^{-}\,\lambda_0(db)\,ds,\quad t \geqslant 0.
\end{displaymath}
The following a priori uniform estimates on the sequence $(Z^{n,x,a},K^{n,x,a})_{n \geqslant 0}$ holds:
\begin{lemma}\label{Sec:PDP_lemma_BSDE_estimations}
	Assume that hypotheses \textup{\textbf{(H$\textup{h$\lambda$Q}$)}}, \textup{\textbf{(H$\lambda_0$)}} and  \textup{\textbf{(H$\textup{f}$)}} hold. 
	For every  $(x,a,n) \in E \times A\times \N$ and for every $T \in (0,\,\infty)$,
	there exists a constant $C$ depending only on $M_f$, $\delta$ and  $T$
	such that
	\begin{eqnarray}\label{Sec:PDP_estimate_bsde}
		||Z^{n,x,a}||^2_{\textbf{\textup{L}}^{\textbf{2}}_\textbf{x,a}(\textup{q};\textup{\textbf{0}},\,\textup{\textbf{T}})} + ||K^{n,x,a}||^2_{\textup{\textbf{K}}^{\textbf{2}}_\textbf{x,a}(\textup{\textbf{0}},\,\textup{\textbf{T}})} \leqslant C.
	\end{eqnarray}
\end{lemma}
\proof
In what follows we shall denote by $C >0$ a generic positive constant depending on $M_f$, $\delta$ and  $T$,
which may vary from line to line.
Fix $T >0$ and apply  It\^o's formula to $|Y_r^{n,x,a}|^2$ between $0$ and $T$. Noticing that $K^{n,x,a}$ is continuous and $\Delta Y^{n,x,a}_r = \int_{E \times A}Z_r^{n,x,a}(y,b)\,p(\{r\}\,dy \,db)$, we get:
$\P^{x,a}$-a.s.,
\begin{align*}
	\sperxa{|Y_0^{n,x,a}|^2} &= \sperxa{|Y_T^{n,x,a}|^2} -2\sperxa{\int_0^T\,|Y_r^{n,x,a}|^2\,\,dr} \\
	&-2 \sperxa{\int_s^T Y_r^{n,x,a}\,dK^{n,x,a}_r} +2\sperxa{\int_0^T\,Y_r^{n,x,a}\,f(X_r,I_r)\,dr}\\
	&-2\sperxa{ \int_0^T\int_{A}\,Y_r^{n,x,a}\,Z_r^{n,x,a}(X_r,b)\,\lambda_0(db)\,dr}\nonumber\\
	&-\sperxa{ \int_0^T\int_{E \times A}\,|Z_r^{n,x,a}(y,b)|^2\,\tilde{p}(dr\,dy\,db)}.
\end{align*}
Let us now denote
$
C_Y := \frac{M_f}{\delta}.
$
Recalling the uniform estimate \eqref{Sec:PDP_Sinfty_estimate}   on $Y^{n,x,a}$, and using elementary inequalities, we get
\begin{eqnarray}\label{Sec:PDP_partial_estimate_Z}
	&&\sperxa{\int_{0}^{T}\int_{E \times A}|Z^{n,x,a}_s(y,\,b)|^2 \, \tilde{p}(ds\,dy\, db)}  \nonumber\\
	&&\leqslant C_Y^2 +2 \, T\,C_Y^2\, +2\,T\,C_Y\,M_f
	+ 2\,C_Y\,T\,\sperxa{ |K_T^{n,x,a}|}\nonumber\\
	&&\quad+ \frac{C_Y}{\alpha}\,T\,\lambda_0(A) + \alpha \,C_Y\,\sperxa{\int_0^T\int_{A}\,|Z_r^{n,x,a}(X_s,b)|^2\,\lambda_0(db)\,dr},
\end{eqnarray}
for any $\alpha >0$.
Now, from relation \eqref{Sec:PDP_BSDE_penalized}, we obtain:
\begin{eqnarray}\label{eq_Kn}
	K_T^{n,x,a}  &=& Y_0^{n,x,a} - Y_T^{n,x,a} -\delta\, \int_{0}^{T}\int_A \,Y_s^{n,x,a}\,ds \nonumber\\
	&& + \int_0^T f(X_s,I_s)ds + \int_{0}^{T}\int_A \,Z_s^{n,x,a}(X_s,b)\,\lambda_0(db)\,ds\nonumber\\
	&& + \int_{0}^{T}\int_{E \times A} \,Z_s^{n,x,a}(y,b)\,q(ds\,dy\,db).
\end{eqnarray}
Then, using  the inequality  $2bc \leqslant \frac{1}{\beta}b^2 + \beta c^2$, for any $\beta >0$, and taking the expected value we have
\begin{eqnarray}\label{Sec:PDP_partial_estim_K^n}
	2\,\sperxa{|K_T^{n,x,a}|}  &\leqslant&  2\,\delta\,C_Y\,T
	+ 2\,M_f\,T+ \frac{T}{\beta}\,\lambda_0(A)\nonumber\\
	&&+ \beta\,\sperxa{\int_{0}^{T}\int_A \,|Z_s^{n,x,a}(X_s,b)|^2\,\lambda_0(db)\,ds}.
\end{eqnarray}
Plugging \eqref{Sec:PDP_partial_estim_K^n} into \eqref{Sec:PDP_partial_estimate_Z},
we get
\begin{align*}
	&\sperxa{\int_{0}^{T}\int_{E \times A}|Z^{n,x,a}_s(y,\,b)|^2 \, \tilde{p}(ds\,dy\, db)}\\
	&\leqslant
	C +C_Y\,\left(2\,T\,\beta+\alpha\right)\,\int_{0}^{T}\int_A \,|Z_s^{n,x,a}(X_s,b)|^2\,\lambda_0(db)\,ds.
\end{align*}
Hence,
choosing $\alpha + 2\,T\,\beta = \frac{1}{2\,C_Y}$, we get
\begin{eqnarray*}
	\frac{1}{2}\,\sperxa{\int_{0}^{T}\int_{E \times A}|Z^{n,x,a}_s(y,\,b)|^2 \, \tilde{p}(ds\,dy\, db)}
	&\leqslant&
	C,
\end{eqnarray*}
which gives the required uniform estimate  for $(Z^{n,x,a})$,
and also $(K^{n,x,a})$ by \eqref{Sec:PDP_partial_estim_K^n}.
\endproof

\subsection{BSDE representation of the dual value function}

To prove  the main result of this section we will need  the following Lemma.
\begin{lemma}\label{Sec:PDP_Lem_intermediate}
	Assume that Hypotheses \textup{\textbf{(H$\textup{h$\lambda$Q}$)}}, \textup{\textbf{(H$\lambda_0$)}} and  \textup{\textbf{(H$\textup{f}$)}} hold.
	For every  $(x,a) \in E \times A$, let $(Y^{x,a},Z^{x,a},K^{x,a})\in \textup{\textbf{S}}^{\infty}\times \textup{\textbf{L}}^{\textbf{2}}_\textbf{x,a,\textup{loc}}({\textup{q}})\times \textup{\textbf{K}}^\textbf{2}_\textbf{x,a,\textup{loc}}$  be a solution to the BSDE with partially nonnegative jumps \eqref{Sec:PDP_BSDE}-\eqref{Sec:PDP_BSDE_constraint}. 
	Then,
	\begin{equation}\label{Sec:PDP_rep_inf_Y}
		Y_s^{x,a} \leqslant \essinf_{\nu \in \mathcal{V}}\spernuxa{\int_s^\infty e^{-\delta (r-s)}\,f(X_{r},I_{r})\, dr \Big | \mathcal{F}_s}, \,\, \forall \,\,s \geqslant 0.
	\end{equation}
\end{lemma}
\proof
Let $(x,a)\in E \times A$, and
consider a triplet $(Y^{x,a},Z^{x,a},K^{x,a}) \in \textbf{S}^{\infty} \times \textbf{L}^{\textbf{2}}_\textbf{x,a,\textup{loc}}(\textup{q}) \times \textbf{K}^\textbf{2}_\textbf{x,a,\textup{loc}}$ satisfying \eqref{Sec:PDP_BSDE}-\eqref{Sec:PDP_BSDE_constraint}.
Applying It\^o's formula to $e^{-\delta\,r}\,Y_r^{x,a}$ between $s$ and  $T>s$, and recalling that $K^{x,a}$ is nondecreasing, we have
\begin{align}\label{Sec:PDP_Y_max}
	Y_{s}^{x,a}&\leqslant e^{-\delta\,(T-s)}\,Y_{T}^{x,a} + \int_{s}^{T} e^{-\delta\,(r-s)}\,f(X_{r},I_{r})\, dr \nonumber\\
	&- \int_{s}^{T} \int_{A}e^{-\delta\,(r-s)}\, Z_{r}^{x,a}(X_{r},b) \,\lambda_{0}(db)\, dr\nonumber\\
	& - \int_{s}^{T} \int_{E \times A}e^{-\delta\,(r-s)}\, Z_{r}^{x,a}(y,b) \, \tilde{q}(dr\,dy\,db),\quad 0\leqslant s\leqslant T< \infty.
\end{align}
Then
for any $\nu \in \mathcal{V}$, let us introduce the compensated martingale measure $q^\nu(ds \, dy \,da)= q(ds \, dy \,db)- (\nu_s(b)-1)\,d_1(s,y,b)\,\tilde{p}(ds\,dy\,db)$ under $\P^{x,a}_\nu$.
Taking expectation in \eqref{Sec:PDP_Y_max} under $\P^{x,a}_\nu$, conditional to $\mathcal{F}_s$, and  recalling that $Z^{x,a}$ is in $\textbf{L}^{\textbf{2}}_\textbf{x,a,\textup{loc}}(\textup{q})$,
we get from Lemma \ref{Sec:PDP_lemma_P_nu_martingale} that, $\P^{x,a}$-a.s.,
\begin{eqnarray}\label{Sec:PDP_BSDE_spernu1}
	Y_{s}^{x,a} &\leqslant& \spernuxa{e^{-\delta\,(T-s)}\,Y_{T}^{x,a} + \int_{s}^{T} e^{-\delta\,(r-s)}\,f(X_{r},I_{r})\, dr \Big | \mathcal{F}_s}\nonumber\\
	&&- \spernuxa{\int_{s}^{T} \int_{A} e^{-\delta\,(r-s)}\, \nu_r(a)\,\bar{Z}_{r}^{x,a}(X_r,b)\, \lambda_0(db)\, dr\Big | \mathcal{F}_s}.
\end{eqnarray}
Furthermore,  since $\nu$ is strictly positive and $Z^{x,a}$ satisfies the nonnegative constraint \eqref{Sec:PDP_BSDE_constraint}, from inequality \eqref{Sec:PDP_BSDE_spernu1} we get
\begin{eqnarray*}
	Y_{s}^{x,a}
	&\leqslant& \spernuxa{e^{-\delta\,(T-s)}\,Y_{T}^{x,a} + \int_{s}^{T} e^{-\delta\,(r-s)}\,f(X_{r},I_{r})\, dr \Big | \mathcal{F}_s}\nonumber\\
	&\leqslant& \spernuxa{e^{-\delta\,(T-s)}\,Y_{T}^{x,a} + \int_{s}^{\infty} e^{-\delta\,(r-s)}\,f(X_{r},I_{r})\, dr \Big | \mathcal{F}_s}.
\end{eqnarray*}
Finally, sending $T \rightarrow \infty$ and recalling that  $Y^{x,a}$ is in $\textbf{S}^{\infty}$,
the conditional version of  Lebesgue dominated convergence theorem leads to
\begin{eqnarray*}
	Y_s^{x,a} &\leqslant&  \spernuxa{\int_{s}^{\infty} e^{-\delta\,(r-s)}\,f(X_{r},I_{r})\, dr\Big | \mathcal{F}_s}
\end{eqnarray*}
for all $\nu\in\mathcal V$, and the conclusion follows from the arbitrariness of  $\nu \in \mathcal{V}$, .
\endproof
Now we are  ready to state the main result of the section.
\begin{theorem}\label{Sec:PDP_Thm_ex_uniq_maximal_BSDE}
	Under Hypotheses  \textup{\textbf{(H$\textup{h$\lambda$Q}$)}}, \textup{\textbf{(H$\lambda_0$)}} and  \textup{\textbf{(H$\textup{f}$)}},
	for every  $(x,a) \in E \times A$, there exists a unique maximal solution $(Y^{x,a},Z^{x,a},K^{x,a})\in \textup{\textbf{S}}^{\infty}\times \textup{\textbf{L}}^{\textbf{2}}_\textbf{x,a,\textup{loc}}(\textup{q})\times \textup{\textbf{K}}^\textbf{2}_\textbf{x,a,\textup{loc}}$  to the BSDE with partially nonnegative jumps \eqref{Sec:PDP_BSDE}-\eqref{Sec:PDP_BSDE_constraint}. 
	In particular,
	\begin{itemize}
		\item[(i)] $Y^{x,a}$ is the nondecreasing limit of $(Y^{n,x,a})_n$;
		\item[(ii)] $Z^{x,a}$ is the weak limit of $(Z^{n,x,a})_n$ in $\textup{\textbf{L}}^{\textbf{2}}_\textbf{x,a,\textup{loc}}(\textup{q})$;
		\item[(iii)] $K_{s}^{x,a}$ is the weak limit of $(K_{s}^{n,x,a})_n$ in $\textup{\textbf{L}}^{\textbf{2}}(\mathcal{F}_{s})$, for any $s \geqslant 0$;
	\end{itemize}
	Moreover, $Y^{x,a}$  has the explicit representation:
	\begin{equation}\label{Sec:PDP_rep_Y}
		Y_s^{x,a} = \essinf_{\nu \in \mathcal{V}}\spernuxa{\int_s^\infty e^{-\delta (r-s)}\,f(X_{r},I_{r})\, dr\Big| \mathcal{F}_s}, \,\, \forall \,\,s \geqslant 0.
	\end{equation}
	In particular, setting $s=0$, 
	we have the following representation formula for the value function of the dual control problem:
	\begin{equation}\label{Sec:PDP_Vstar_Y0}
		V^{\ast}(x,a)
		= Y_0^{x,a}, \quad (x,a)\in E \times A.
	\end{equation}
\end{theorem}
\proof
Let $(x,a) \in E \times A$ be fixed.
From the representation formula \eqref{Sec:PDP_rep_Y_n} it  follows that
$Y_s^{n} \geqslant Y_s^{n+1}$ for all $s \geqslant 0$ and all $n \in \N$,  since by definition $\mathcal{V}^n \subset \mathcal{V}^{n+1}$ and $(Y^n)_n$ are càdlàg processes.
Moreover, recalling  the boundedness of $f$, from \eqref{Sec:PDP_rep_Y_n} we see that  $(Y^n)_n$ is lower-bounded by a constant which does not depend $n$.
Then $(Y^{n,x,a})_n \in \textbf{S}^\infty$  converges decreasingly to some adapted process $Y^{x,a}$, which is moreover uniformly bounded by Fatou's lemma.
Furthermore,  for every $T >0$,  the Lebesgue's dominated convergence theorem insures that the convergence of $(Y^{n,x,a})_n$ to $Y$ also holds in $\textbf{L}^\textbf{2}(\textbf{0},\,\textbf{T})$.

Let us  fix $T \geqslant 0$.
By the uniform estimates in Lemma \ref{Sec:PDP_lemma_BSDE_estimations}, the sequence $(Z^{n,x,a}_{|[0,\,T]})_n$ is bounded in the Hilbert space $\textbf{L}^{\textbf{2}}_\textbf{x,a}(\textup{q};\textbf{0},\,\textbf{T})$. Then, we can extract a subsequence which weakly converges to some $Z^{T}$ in $\textbf{L}^{\textbf{2}}_\textbf{x,a}(\textup{q};\textbf{0},\,\textbf{T})$.
Let us then define the following mappings
\begin{eqnarray*}
	I_{\tau}^1 := \qquad  \qquad  Z &\longmapsto& \int_0^\tau \int_{E \times A}\, Z_s(y,b)\, q(ds\,dy \,db)\\
	\textup{\textbf{L}}^{\textbf{2}}_\textbf{x,a}(\textup{q};\textup{\textbf{0}},\,\textup{\textbf{T}})
	&\longrightarrow& \textup{\textbf{L}}^\textbf{2}(\mathcal{F}_{\tau}),
\end{eqnarray*}
\begin{eqnarray*}
	I_{\tau}^2 := \qquad  \qquad  Z(X_s,\cdot) &\longmapsto& \int_0^\tau \int_{A }\, Z_s(X_s,b)\, \lambda_0(db)\,ds\\
	\textup{\textbf{L}}^\textbf{2}_\textbf{x,a}(\lambda_0;\textup{\textbf{0}},\,\textup{\textbf{T}})
	&\longrightarrow& \textup{\textbf{L}}^\textbf{2}(\mathcal{F}_{\tau}),
\end{eqnarray*}
for every stopping time $0 \leqslant \tau \leqslant T$.
We notice that $I^1_{\tau}$ (resp., $I^2_{\tau}$) defines a linear continuous operator (hence weakly continuous) from $\textup{\textbf{L}}^{\textbf{2}}_\textbf{x,a}(\textup{q};\textup{\textbf{0}},\,\textup{\textbf{T}})$ (resp., $\textup{\textbf{L}}^\textbf{2}_\textbf{x,a}(\lambda_0;\textup{\textbf{0}},\,\textup{\textbf{T}})$) to $\textup{\textbf{L}}^\textbf{2}(\mathcal{F}_{\tau})$. Therefore $I_{\tau}^1 Z^{n,x,a}_{|[0,\,T]}$ (resp., $I_{\tau}^2 Z^{n,x,a}_{|[0,\,T]}(X,\cdot)$) weakly converges to  $I_{\tau}^1 \tilde{Z}^{T}$ (resp., $I_{\tau}^2 \tilde{Z}^{T}(X,\cdot)$) in $\textup{\textbf{L}}^\textbf{2}(\mathcal{F}_{\tau})$.
Since
\begin{eqnarray*}
	K_\tau^{n,x,a} &=& Y_{\tau}^{n,x,a}- Y_{0}^{n,x,a} -\delta\int_0^\tau Y_r^{n,x,a}\,dr + \int_{0}^{\tau} f(X_{r},I_{r})\, dr \nonumber\\
	&&- \int_{0}^{\tau} \int_{A} Z_{r}^{n,x,a}(X_{r},b) \,\lambda_{0}(db)\, dr \\
	&&- \int_{0}^{\tau} \int_{E \times A} Z_{r}^{n,x,a}(y,b) \, q(dr\,dy\,db),\,\,\forall\,\,\tau \in [0,\,T],
\end{eqnarray*}
we also have the following weak convergence in $\textbf{L}^\textbf{2}(\mathcal{F}_\tau)$:
\begin{eqnarray*} 
	K_\tau^{n,x,a} \rightharpoonup  \tilde{K}_\tau^{T} &:=& Y_{\tau}^{x,a}- Y_{0}^{x,a} -\delta\int_0^\tau Y_r^{x,a}\,dr + \int_{0}^{\tau} f(X_{r},I_{r})\, dr\nonumber\\
	&&- \int_{0}^{\tau} \int_{A} Z_{r}^{x,a}(X_{r},b) \,\lambda_{0}(db)\, dr \\
	&&- \int_{0}^{\tau} \int_{E \times A} Z_{r}^{x,a}(y,b) \, q(dr\,dy\,db),\,\,\forall \,\,\tau \in [0,\,T].
\end{eqnarray*}
Since the process $(K^{n,x,a}_{s})_{s \in [0,\,T]}$ is nondecreasing and predictable and $K^{n,x,a}_{0}=0$, the limit process $\tilde{K}_\tau^{T}$ on $[0,\,T]$ remains nondecreasing and predictable with $\sperxa{|\tilde{K}_T^{T}|^2}< \infty$ and $\tilde{K}_0^{T}=0$. Moreover, by Lemma 2.2. in \cite{Pe}, $\tilde{K}_\tau^{T}$ and $\tilde{Y}_\tau^{T}$  are càdlàg, therefore  $\tilde{K}_\tau^{T}\in \textbf{K}^\textbf{2}_\textbf{x,a}(\textbf{0},\,\textbf{T})$ and $\tilde{Y}_\tau^{T}\in \textbf{S}^\infty$.

Then we notice that $\tilde{Z}^{T'}_{|[0,\,T]} = \tilde{Z}^T$, $\tilde{K}^{T'}_{|[0,\,T]} = \tilde{K}^T$, for any $0 \leqslant T \leqslant T'< \infty$.
Indeed, for $i=1,2$, $I^{i}\,\tilde{Z}^{T'}_{|[0,\,T]}$, as $I^{i}\,\tilde{Z}^T$, is the weak limit in $\textbf{L}^{\textbf{2}}(\mathcal{F}_s)$
of $(I^{i}\,Z^{n,x,a}_{|[0,\,T]})_{n\geqslant 0}$, while $\tilde{K}^{T'}_{|[0,\,T]}$, as $\tilde{K}^T$, is the weak limit in $\textbf{L}^{\textbf{2}}(\mathcal{F}_s)$ of $(K^{n,x,a}_{|[0,\,T]})_{n \geqslant 0}$, for every $s \in [0,\,T]$. Hence, we define $Z_s^{x,a} = \tilde{Z}_s^T$, $K_s^{x,a} = \tilde{K}_s^T$ for all $s \in [0,\,T]$ and for any $T >0$. Observe that $Z^{x,a}\in \textbf{L}^{\textbf{2}}_\textbf{x,a,\textup{loc}}(\textup{q})$ and $K^{x,a}\in \textbf{K}^{\textbf{2}}_\textbf{x,a,\textup{loc}}$. Moreover, for any $T>0$, for $i=1,2$, $(I^{i}\,Z^{n,x,a}_{|[0,\,T]})_{n \geqslant 0}$ weakly converges to $I^{i}\,Z^{x,a}_{|[0,\,T]}$ in 
$\textbf{L}^{\textbf{2}}(\mathcal{F}_s)$, and $(K^{n,x,a}_{|[0,\,T]})_{n \geqslant 0}$ weakly converges   to $K^{x,a}_{|[0,\,T]}$ in $\textbf{L}^{\textbf{2}}(\mathcal{F}_s)$, for $s \in [0,\,T]$.
In conclusion, we have: $\P^{x,a}$-a.s.,
\begin{eqnarray*}
	Y_{s}^{x,a}&=& Y_{T}^{x,a}-\delta\int_s^T Y_r^{x,a}\,dr + \int_{s}^{T} f(X_{r},I_{r})\, dr  -( K_T^{x,a} - K_s^{x,a,\delta}) \nonumber\\
	&&- \int_{s}^{T} \int_{A} Z_{r}^{x,a}(X_{r},b) \,\lambda_{0}(db)\, dr\\
	&&  - \int_{s}^{T} \int_{E \times A} Z_{r}^{x,a}(y,b) \, q(dr\,dy\,db),\quad 0\leqslant s\leqslant T.
\end{eqnarray*}
Since $T$ is arbitrary, it follows that $(Y^{x,a},Z^{x,a},K^{x,a})$ solves equation \eqref{Sec:PDP_BSDE} on $[0,\,\infty)$.

To show that the jump constraint \eqref{Sec:PDP_BSDE_constraint} is satisfied, we consider the functional $G:\textbf{L}^\textbf{2}_\textbf{x,a}(\lambda_0;\textbf{0},\,\textbf{T})$ $\rightarrow \R$ given by
\begin{eqnarray*}
	G(V(\cdot)) \,\, := \,\, \sper{\int_{0}^{T}\int_A \, [V_s(b)]^-\, \lambda_0(db)\, ds},\quad \forall \,\,V \in \textbf{L}^\textbf{2}_\textbf{x,a}(\lambda_0;\textbf{0},\,\textbf{T}).
\end{eqnarray*}
Notice that $G(Z^{n,x,a}(X,\cdot)) = \sperxa{K_{T}^{n,x,a}/n}$, for any $n \in \N$. From uniform estimate \eqref{Sec:PDP_estimate_bsde}, we see that
$G(Z^{n,x,a}(X,\cdot))\rightarrow 0$ as $n \rightarrow \infty$. Since $G$ is convex and strongly continuous in the strong topology of $\textbf{L}^\textbf{2}_\textbf{x,a}(\lambda_0;\textbf{0},\,\textbf{T})$, then $G$ is lower semicontinuous in the weak topology of $\textbf{L}^\textbf{2}_\textbf{x,a}(\lambda_0;\textbf{0},\,\textbf{T})$, see, e.g., Corollary 3.9 in \cite{Bre}. Therefore, we find
\begin{displaymath}
	G(Z^{x,a}(X,\cdot))\leqslant \liminf_{n \rightarrow \infty}G(Z^{n,x,a}(X,\cdot))=0,
\end{displaymath}
which implies the validity of jump constraint \eqref{Sec:PDP_BSDE_constraint} on $[0,\,T]$, and the conclusion follows from the arbitrary of $T$.

Hence, $(Y^{x,a},Z^{x,a},K^{x,a})$ is a solution to the constrained BSDE \eqref{Sec:PDP_BSDE}-\eqref{Sec:PDP_BSDE_constraint} on $[0,\,\infty)$.

It remains to prove  the representation formula \eqref{Sec:PDP_rep_Y} and the maximality property for $Y^{x,a}$. 
Firstly, since  by definition $\mathcal{V}^n \subset \mathcal{V}$ for all $n \in \N$, it is clear from representation formula \eqref{Sec:PDP_rep_Y_n} that
\begin{eqnarray*}
	Y_s^{n,x,a} &=&  \essinf_{\nu \in \mathcal{V}^{n}}\spernuxa{\int_{s}^{\infty} e^{-\delta\,(r-s)}\,f(X_{r},I_{r})\, dr \Big | \mathcal{F}_s}\\
	&\geqslant&  \essinf_{\nu \in \mathcal{V}}\spernuxa{\int_{s}^{\infty} e^{-\delta\,(r-s)}\,f(X_{r},I_{r})\, dr \Big | \mathcal{F}_s},
\end{eqnarray*}
for all $n \in \N$, for all $s \geqslant 0$.
Moreover, being $Y^{x,a}$  the pointwise  limit of $Y^{n,x,a}$, we deduce that
\begin{equation}\label{Sec:PDP_Y_ineq}
	Y_s^{x,a} \geqslant \essinf_{\nu \in \mathcal{V}}\spernuxa{\int_{s}^{\infty} e^{-\delta\,(r-s)}\,f(X_{r},I_{r})\, dr| \mathcal{F}_s}, \quad s \geqslant 0.
\end{equation}
On the other hand, $Y^{x,a}$ satisfies the opposite inequality \eqref{Sec:PDP_rep_inf_Y} from Lemma \ref{Sec:PDP_Lem_intermediate}, and thus we achieve the representation formula \eqref{Sec:PDP_rep_Y}.

Finally, to show that $Y^{x,a}$  is the maximal solution,  let
consider  a triplet $(\bar{Y}^{x,a},\bar{Z}^{x,a},\bar{K}^{x,a}) \in \textbf{S}^{\infty} \times \textbf{L}^{\textbf{2}}_\textbf{x,a,\textup{loc}}(\textup{q}) \times \textbf{K}^\textbf{2}_\textbf{x,a,\textup{loc}}$ solution to \eqref{Sec:PDP_BSDE}-\eqref{Sec:PDP_BSDE_constraint}.
By Lemma \ref{Sec:PDP_Lem_intermediate}, $(\bar{Y}^{x,a},\bar{Z}^{x,a},\bar{K}^{x,a})$ satisfies inequality \eqref{Sec:PDP_rep_inf_Y}.
Then, from the representation formula \eqref{Sec:PDP_rep_Y} it follows   that $\bar{Y}^{x,a}_s \leqslant Y^{x,a}_s$, $\forall \, s\geqslant 0$, $\P^{x,a}$-a.s., i.e., the maximality property holds.
The uniqueness of the minimal solution directly follows from Proposition \ref{Sec:PDP_Prop_uniq_max_sol}.
\endproof

\section{A BSDE representation for the value function}\label{Sec:PDP_Section_nonlinear_IPDE}
Our main purpose is to show how maximal solutions to BSDEs with nonnegative jumps of the form \eqref{Sec:PDP_BSDE}-\eqref{Sec:PDP_BSDE_constraint} provide actually a Feynman-Kac representation to the value function $V$ associated to our optimal control problem for  PDMPs.
We know from Theorem \ref{Sec:PDP_Thm_ex_uniq_maximal_BSDE} that, under Hypotheses \textup{\textbf{(H$\textup{h$\lambda$Q}$)}}, \textup{\textbf{(H$\lambda_0$)}} and  \textup{\textbf{(H$\textup{f}$)}} , there exists a unique maximal solution $(Y^{x,a},Z^{x,a},K^{x,a})$ on $(\Omega, \mathcal{F}, \mathbb F, \mathbb P^{x,a})$ to \eqref{Sec:PDP_BSDE}-\eqref{Sec:PDP_BSDE_constraint}.  
Let us introduce a deterministic function $v: E \times A \rightarrow \R$ as
\begin{equation}
v(x,a):= Y_0^{x,a}, \quad (x,a) \in E\times A. \label{Sec:PDP_def_v}
\end{equation}
%
%
%
Our main result is as follows:
\begin{theorem}\label{Sec:PDP_THm_Feynman_Kac_HJB}
	Assume that Hypotheses \textup{\textbf{(H$\textup{h$\lambda$Q}$)}}, \textup{\textbf{(H$\lambda_0$)}}, 
	and  \textup{\textbf{(H$\textup{f}$)}}  hold. Then the function $v$  
	in \eqref{Sec:PDP_def_v} does not depend on the variable $a$:
	$$
	v(x,a)= v(x,a'),\quad \forall a,a' \in A,
	$$
	for all $x \in E$. Let us define by misuse of notation the function $v$ on $E$ by
	$$
	v(x)= v(x,a),\quad \forall x \in E,
	$$
	for any $a \in A$. Then $v$ is a (discontinuous) viscosity  solution to \eqref{Sec:PDP_HJB}.
\end{theorem}

%


To conclude that $v(x)$ actually provides the unique solution to \eqref{Sec:PDP_HJB} (and therefore by Theorem \ref{Sec:PDP_Thm_unique_viscosity_sol_HJB} coincides with the value function $V$), we need to use a 
comparison theorem for viscosity sub and supersolutions to the fully nonlinear IPDE of HJB type. 
To this end, we  introduce  the  following additional condition on  $Q$. 

\vspace{3mm}

\textbf{(H$\textup{Q'}$)} 
\begin{itemize}
	\item[(i)] $\sup_{(x,a)\in E\times A} \int_{E}|y-x|\,\lambda(x,a)\, Q(x,a,dy) < \infty$;
	\item[(ii)]$\exists$ $c, C>0$:  for every $\psi \in W^{1,\,\infty}(E)$, $\psi(0)=0$, and for every  $K\subset E$ compact set,
	\begin{align*}
	&\Big|\int_{K+x_1} \psi(y-x_1)\,\lambda(x_1,a)\,(Q(x_1,a,dy)-\int_{K+x_2} \psi(y-x_2)\,\lambda(x_2,a)\,Q(x_2,a,dy)\Big|\\
	&\leqslant c ||\nabla \psi||_{\infty}||x_1-x_2|| 
	\end{align*}
	and
	\begin{align*}
	&\Big|\int_{K^c+x_1} \psi(y-x_1)\,\lambda(x_1,a)\,Q(x_1,a,dy)-\int_{K^c+x_2} \psi(y-x_2)\,\lambda(x_2,a)\,Q(x_2,a,dy)\Big|\\
	&\leqslant C ||\nabla \psi||_{\infty}||x_1-x_2||, 
	\end{align*}
	for every $x_1,x_2\in E$, $a \in A$.
\end{itemize}
\begin{corollary}\label{Sec:PDP_C_final}
	Let Hypotheses \textup{\textbf{(H$\textup{h$\lambda$Q}$)}}, \textup{\textbf{(H$\lambda_0$)}}, \textup{\textbf{(H$\textup{Q'}$)}}
	and  \textup{\textbf{(H$\textup{f}$)}}  hold, and assume that $A$ is compact. Then 
	the value function $V$ of the optimal control problem defined  in \eqref{Sec:PDP_value_function} admits the Feynman-Kac representation formula:
	\begin{equation*} \label{Sec:PDP_Feynman-Kac}
	V(x)=Y^{x,a}_0,\quad (x,a)\in E \times A.
	\end{equation*}
	Moreover, the  value function $V$ coincides with the dual value function $V^\ast$ defined in \eqref{Sec:PDP_dual_value_function}, namely
	\begin{equation}\label{Sec:PDP_equality_value_functions}
	V(x)= V^\ast(x,a)=Y^{x,a}_0,\quad (x,a)\in E \times A.
	\end{equation}
\end{corollary}
\proof
Under the additional assumption \textbf{(H$\textup{Q'}$)}, 
a  comparison theorem   for viscosity  super and subsolutions for  elliptic  
IPDEs of the form \eqref{Sec:PDP_HJB} holds, see 
Theorem IV.1 in Sayah \cite{Sa}. 
Then, it follows from 
Theorem \ref{Sec:PDP_THm_Feynman_Kac_HJB} that the function $v$ in \eqref{Sec:PDP_def_v} is the unique viscosity soluton to \eqref{Sec:PDP_HJB}, and it is continuous. 
In particular, by Theorem \ref{Sec:PDP_Thm_unique_viscosity_sol_HJB},  $v$ coincides with   the value function $V$ of the PDMPs optimal control problem, which admits therefore the probabilistic representation 
\eqref{Sec:PDP_Feynman-Kac}.
Finally,   Theorem \ref{Sec:PDP_Thm_ex_uniq_maximal_BSDE} implies that the dual value function $V^\ast$ coincides with the value function $V$ of the original control problem, so that \eqref{Sec:PDP_equality_value_functions} holds.
\endproof

The rest of the chapter is devoted to prove Theorem \ref{Sec:PDP_THm_Feynman_Kac_HJB}.

\subsection{The identification property of the penalized BSDE}
For every $n\in \N$, 
let us introduce the deterministic function $v^{n}$ defined on $E \times A$ by
\begin{equation}
	v^{n}(x,a)= Y_0^{n,x,a}, \quad (x,a) \in E\times A. \label{Sec:PDP_def_vn}
\end{equation}
We  investigate the properties of the function $v^{n}$.
Firstly, it straightly follows from \eqref{Sec:PDP_def_vn} and 
\eqref{Sec:PDP_Sinfty_estimate} that
$$
|v^{n}(x,a)|\leqslant \frac{M_f}{\delta}, \quad \forall \,(x,a)\in E \times A.
$$
Moreover, we have the following  result.
\begin{lemma}\label{Sec:PDP_Lem_identification_vn}
	Under Hypotheses \textup{\textbf{(H$\textup{h$\lambda$Q}$)}}, \textup{\textbf{(H$\lambda_0$)}} 
	and  \textup{\textbf{(H$\textup{f}$)}}, for any 
	$n \in \N$, the function $v^{n}$ is such that, for any $(x,a)\in E \times A$, we have
	\begin{equation}\label{Sec:PDP_Y_identification}
		Y_s^{n,x,a}= v^{n}(X_{s}, I_s),\quad s \geqslant 0\quad  d\P^{x,a} \otimes ds\textup{ -a.e.}
	\end{equation}
\end{lemma}
\begin{remark}\label{Sec:PDP_Rem_iterative_meth_for_ident}
When the pair of Markov processes $(X,I)$  is the unique strong solution to some system of  stochastic differential equations, $(X,I)$ often satisfies a  stochastic flow property, and 
the fact that $Y_s^{n,x,a}$ is a deterministic function of  $(X_s,I_s)$ 
straight follows from the uniqueness of the BSDE
(see, e.g., Remark 2.4 in  Barles, Buckdahn and Pardoux \cite{BaBuPa}). In our framework,   we deal with the local characteristics of the state process $(X,I)$ rather than with the stochastic differential equation solved by it. As a consequence,  a stochastic flow property for  $(X,I)$ is no more directly available. The idea is then  to prove the identification \eqref{Sec:PDP_Y_identification} using an  iterative construction of the solution of standard BSDEs. This alternative approach  is based on the fact that,  when $f$ does not depend on $y,z$, the desired identification follows from the Markov property of the state process $(X,I)$, 
	and  it is inspired by the proof of  the Theorem 4.1. in El Karoui, Peng and Quenez \cite{EPQ}.
\end{remark}
\proof
Fix $(x,a,n)\in E \times A \times \N$.
Let  $(Y^n,Z^n)= (Y^{n,x,a},Z^{n,x,a})$ be the solution to the penalized BSDE \eqref{Sec:PDP_BSDE_penalized}.
From Proposition \ref{Sec:PDP_Prop_exist_uniq_pen_BSDE} we know that there exists a sequence $(Y^{n,T},Z^{n,T})_T=(Y^{n,T,x,a}, Z^{n,T,x,a})_T$ in $\textbf{S}^{\infty}\times \textbf{L}^{\textbf{2}}_\textbf{x,a,\text{loc}}(\textup{q})$ such that, when $T$ goes to infinity,  $(Y^{n,T})_T$ converges $\P^{x,a}$-a.s. to $(Y^{n})$ and $(Z^{n,T})_T$ converges to $(Z^n)$ in $\textbf{L}^{\textbf{2}}_\textbf{x,a,\text{loc}}(\textup{q})$.
Let us now fix $T,S>0$, $S<T$, and consider the BSDE solved by $(Y^{n,T},Z^{n,T})$ on $[0,\,S]$:
\begin{eqnarray*}
	Y_{t}^{n,T}&=& Y_{S}^{n,T}-\delta\int_t^S Y_r^{n,T}\,dr + \int_{t}^{S} f(X_{r},I_{r})\, dr \nonumber\\
	&& -n\int_t^S\int_A [Z_r^{n,T}(X_r,b)]^-\,\lambda_0(db)\,dr - \int_{t}^{S} \int_{A} Z_{r}^{n,T}(X_{r},b) \,\lambda_{0}(db)\, dr, \nonumber\\
	&& - \int_{t}^{S} \int_{E \times A} Z_{r}^{n,T}(y,b) \, q(dr\,dy\,db),\quad \quad 0 \leqslant t\leqslant S.
\end{eqnarray*}
Then, it follows from Proposition \ref{Sec:PDP_Prop_ex_uniq_BSDE_pen_T} that there exists a sequence $(Y^{n,T,k},Z^{n,T,k})_k=(Y^{n,T,k,x,a}, Z^{n,T,k,x,a})_k$ in $\textbf{L}^{\textbf{2}}_\textbf{x,a}( \textbf{0,\,S})\times \textbf{L}^{\textbf{2}}_\textbf{x,a}(\textup{q}, \textbf{0,\,S})$ converging to $(Y^{n,T},Z^{n,T})$ in $\textbf{L}^{\textbf{2}}_\textbf{x,a}( \textbf{0,\,S})\times \textbf{L}^{\textbf{2}}_\textbf{x,a}(\textup{q}, \textbf{0,\,S})$, such that $(Y^{n,T,0}, Z^{n,T,0})= (0,0)$ and
\begin{eqnarray*}
	Y_{t}^{n,T,k+1}&=& Y_{S}^{n,T,k}-\delta\int_t^S Y_r^{n,T,k}\,dr + \int_{t}^{S} f(X_{r},I_{r})\, dr \nonumber\\
	&& -n\int_t^S\int_A [Z_r^{n,T,k}(X_r,b)]^-\,\lambda_0(db)\,dr - \int_{t}^{S} \int_{A} Z_{r}^{n,T,k}(X_{r},b) \,\lambda_{0}(db)\, dr, \nonumber\\
	&& - \int_{t}^{S} \int_{E \times A} Z_{r}^{n,T,k+1}(y,b) \, q(dr\,dy\,db),\quad \quad 0 \leqslant t\leqslant S.
\end{eqnarray*}
Let us define
$$
v^{n,T}(x,a):= Y_{0}^{n,T}, \quad v^{n,T,k}(x,a):= Y_{0}^{n,T,k}, \quad k \geq0.
$$
We start by noticing that, for $k=0$, we have, $\P^{x,a}$-a.s.,
$$
Y_{t}^{n,T,1}= \sperxa{\int_{t}^{S} f(X_{r},I_{r})\, dr \Big |\mathcal{F}_t}, \quad t \in [0,\,S].
$$
Then, from the Markov property of $(X,I)$ we get
\begin{equation}\label{Sec:PDP_id_1}
	Y_{t}^{n,T,1}= v^{n,T,1}(X_{t},I_{t}),\quad  d\P^{x,a} \otimes dt\textup{ -a.e.}
\end{equation}
Furthermore, identification \eqref{Sec:PDP_id_1} implies
\begin{equation}\label{id_2}
	Z_{t}^{n,T,1}(y,b)= v^{n,T,1}(X_{t-},I_{t-})-v^{n,T,1}(y,b),
\end{equation}
where \eqref{id_2} has to be understood as an equality (almost everywhere) between elements of the space $\textbf{L}^{\textbf{2}}_\textbf{x,a}(\textup{q}; \textbf{0,\,S})$.
At this point we consider the inductive step: $1\leqslant k \in \N$. Assume that, $\P^{x,a}$-a.s.,
\begin{eqnarray*}
	Y_{t}^{n,T,k} &=& v^{n,T,k}(X_{t},I_{t})\\
	Z_{t}^{n,T,k}(y,b) &=& v^{n,T,k}(y,b)-v^{n,T,k}(X_{t-},I_{t-}).
\end{eqnarray*}
Then
\begin{eqnarray*}
	Y_{t}^{n,T,k+1}&=& \mathbb E^{x,a}\bigg[v_{\delta}^{n,T,k}(X_{S},I_{S})
	-\delta\int_t^S v^{n,T,k}(X_{r},I_{r})\,dr + \int_{t}^{S} f(X_{r},I_{r})\, dr \nonumber\\
	&& -n\int_t^S\int_A [v^{n,T,k}(X_{t},b)-v^{n,T,k}(X_{t},I_{t})]^-\,\lambda_0(db)\,dr \nonumber\\
	&&- \int_{t}^{S} \int_{A} v^{n,T,k}(X_{t},b)-v^{n,T,k}(X_{t},I_{t}) \,\lambda_{0}(db)\, dr \Big |\mathcal{F}_t\bigg], \quad \quad 0 \leqslant t\leqslant S.
\end{eqnarray*}
Using again the Markov property of $(X,I)$, we achieve that
\begin{equation}\label{Sec:PDP_id_k}
	Y_{t}^{n,T,k+1}= v^{n,T,k+1}(X_{t},I_{t}),\quad  d\P^{x,a} \otimes dt\textup{ -a.e.}
\end{equation}
Then, applying the It\^o formula to $|Y_{t}^{n,T,k}-Y_{t}^{n,T}|^2$ and taking the supremum of $t$ between $0$ and $S$, one can show that
$$
\sperxa{\sup_{0 \leqslant t \leqslant S}\Big|Y_{t}^{n,T,k}-Y_{t}^{\delta,n,T}\Big|^2} \rightarrow 0 \quad \textup{as } k \textup{ goes to infinity.}
$$
Therefore, $v^{n,T,k}(x,a) \rightarrow v^{n,T}(x,a)$ as $k$ goes to infinity, for all $(x,a) \in E \times A$, from which it follows that
\begin{equation}\label{Sec:PDP_id_3}
	Y_{t}^{n,T,x,a}= v^{n,T}(X_{t},I_{t}),\quad  d\P^{x,a} \otimes dt\textup{ -a.e.}
\end{equation}
Finally, from \eqref{Sec:PDP_conv_Y} we have that $(Y^{n,T,x,a})_T$ converges $\P^{x,a}$-a.s. to $(Y^{n,x,a})$ uniformly on compact sets of $\R$.
Thus, $v^{n,T}(x,a) \rightarrow v^{n}(x,a)$ as $T$ goes to infinity, for all $(x,a) \in E \times A$, and  this gives the requested identification $Y_t^{n,x,a} = v^{n}(X_t,I_t)$, $d\P^{x,a}\otimes dt${ -a.e.}
\endproof

\begin{remark}\label{Sec:PDP_Rem_Y_identification}
	By Proposition \ref{Sec:PDP_Prop_uniq_max_sol}, the  maximal solution to the constrained BSDE \eqref{Sec:PDP_BSDE}-\eqref{Sec:PDP_BSDE_constraint} is the pointwise limit of the solution to the penalized BSDE \eqref{Sec:PDP_BSDE_penalized}. Then, as a byproduct of  Lemma \ref{Sec:PDP_Lem_identification_vn} we have  the following identification  for $v$: $\P^{x,a}$-a.s.,	\begin{equation}\label{Sec:PDP_ident_vdelta1}
	v(X_s, I_s) = Y_s^{x,a}, \quad (x,a)\in E \times A,\,s \geqslant 0.
	\end{equation}
\end{remark}

\subsection{The non-dependence of the function $v$ on the variable $a$.}
We claim that the function $v$ in \ref{Sec:PDP_def_v} does not depend on its last argument:
\begin{equation}\label{Sec:PDP_vdelta_not_dep_a}
v(x,a)= v(x,a'), \quad  a,a' \in A,\quad \textup{for any}\,\, x \in E.
\end{equation}
We recall that,
by \eqref{Sec:PDP_Vstar_Y0} and \eqref{Sec:PDP_def_v}, $v$ coincides with the  value function $V^{\ast}$ of the dual control problem introduced in Section \ref{Sec:PDP_Section_dual_optimal_control}.
Therefore, \eqref{Sec:PDP_vdelta_not_dep_a} holds if we
prove that $V^{\ast}(x,a)$ does not depend on $a$. 
This is insured by the following result.
\begin{proposition}\label{Sec:PDP_P_J_no_dip_a}
	Assume that Hypotheses \textup{\textbf{(H$\textup{h$\lambda$Q}$)}}, \textup{\textbf{(H$\lambda_0$)}} 
	and  \textup{\textbf{(H$\textup{f}$)}} hold. Fix $x\in E$, $a, a' \in A$, and 
	$\nu 
	\in \mathcal V$.
	Then,  there exists a sequence  $(\nu^{\varepsilon})_{\varepsilon}: \Omega \times \R_+ \times A\rightarrow (0,\infty)$,  
	$\nu^{\varepsilon} \in \mathcal V$ for every $\varepsilon>0$,
	such that,
	\begin{equation}\label{Sec:PDP_claim_PDMP}
	\lim_{\varepsilon \rightarrow 0^+} J(x,a',\nu^{\varepsilon}) = J(x,a,\nu). 
	\end{equation}
\end{proposition}
\proof
See Appendix \ref{Sec:PDP_A_proof_Prop}.
\endproof
Identity \eqref{Sec:PDP_claim_PDMP}
implies that
\[
V^{\ast}(x,a') \geq J(x,a,\nu) 
\quad x \in E, \,\,a,a' \in A,
\]
and by the arbitrariness of $\nu$ we conclude that
\[
V^{\ast}(x,a') \geq V^{\ast}(x,a) 
\quad x \in E, \,\,a,a' \in A.
\]
In other words $V^\ast(x,a)=v(x,a)$ does not depend on $a$, and \eqref{Sec:PDP_vdelta_not_dep_a} holds.

\subsection{Viscosity properties of the function $v$.}
Taking into account \eqref{Sec:PDP_vdelta_not_dep_a}, by misuse of notation, we define  the function $v$ on $E$ by
\begin{equation}\label{Sec:PDP_redefine_vdelta}
v(x):= v(x,a),\quad \forall x \in E, \quad 	\textup{for any}\,\, a \in A.
\end{equation}
We shall  prove that the function $v$ in \eqref{Sec:PDP_redefine_vdelta} provides a viscosity solution to \eqref{Sec:PDP_HJB}.
We separate the proof of viscosity subsolution and supersolution properties, which are different. In particular the supersolution property is more delicate and should take into account the maximality property of $Y^{x,a}$.
\begin{remark}
	Identity \eqref{Sec:PDP_ident_vdelta1} in Remark \ref{Sec:PDP_Rem_Y_identification} gives
	\begin{equation}\label{Sec:PDP_ident_vdelta}
	v(X_s) = Y_s^{x,a}, \quad \forall x\in E,\,s \geqslant 0, \quad 	\textup{for any}\,\, a \in A.
	\end{equation}
\end{remark}
\paragraph{Proof of the viscosity subsolution property to \eqref{Sec:PDP_HJB}.}
\begin{proposition}\label{Sec:PDP_Prop_visc_subsol_property_vdelta}
	Let assumptions \textup{\textbf{(H$\textup{h$\lambda$Q}$)}}, \textup{\textbf{(H$\lambda_0$)}} 
	 and  \textup{\textbf{(H$\textup{f}$)}} hold. Then, the function $v$ in \eqref{Sec:PDP_redefine_vdelta} is a viscosity subsolution to  \eqref{Sec:PDP_HJB}.
\end{proposition}
\proof
Let $\bar{x} \in E$,  and let $\varphi \in C^1(E)$ be a test function such that
\begin{equation}\label{Sec:PDP_max_property}
	0= (v^{\ast}-\varphi)(\bar{x})=\max_{x \in E} (v^{\ast}-\varphi)(x).
\end{equation}
By the  definition of $v^{\ast}(\bar{x})$, there exists a sequence $(x_m)_m$ in $E$ such that
$$
x_m \rightarrow \bar{x} \,\,\textup{and}\,\,v(x_m) \rightarrow v^{\ast}(\bar{x})
$$
when $m$ goes to infinity.
By the continuity of $\varphi$, and taking into account \eqref{Sec:PDP_max_property}, it follows that
$$
\gamma_m := \varphi(x_m)  - v(x_m) \rightarrow 0,
$$
when $m$ goes to infinity.
Let $\eta$ be a fixed positive constant and $\tau_m :=\inf \{t \geqslant 0: |\phi(t,x_m)-x_m|\geqslant \eta\}$. Let moreover $(h_m)_m$ be a strictly positive sequence such that
$$
h_m \rightarrow 0 \,\, \textup{and}\,\,\frac{\gamma_m}{h_m} \rightarrow 0,
$$
when $m$ goes to infinity.

We notice that there exists $M \in \N$ such that, for every $m > M$, $h_m \wedge \tau_m=h_m$.
Let us introduce $\bar \tau:=\inf \{t \geqslant 0: |\phi(t,\bar x)- \bar x|\geqslant \eta\}$.
Clearly $\bar \tau>0$.
We show that it does not exists a subsequence $\tau_{n_k}$ of $\tau_n$ such that $\tau_{n_{k}}\rightarrow \tau_0\in [0,\,\bar \tau)$.
Indeed, let $\tau_{n_{k}}\rightarrow \tau_0\in [0,\,\bar \tau)$. In particular $|\phi(\tau_{n_{k}},\bar x)- \bar x|\geqslant \eta$. Then, by the continuity of $\phi$ it follows that $|\phi(\tau_0,\bar x)- \bar x|\geqslant \eta$, and this is in contradiction with the definition of $\bar \tau$.

Let us now fix $a \in A$, and let $Y^{x_m,a}$ be the unique maximal solution to \eqref{Sec:PDP_BSDE}-\eqref{Sec:PDP_BSDE_constraint}  under $\P^{x_m,a}$.
We apply the It\^o formula to $e^{-\delta t} \,Y_t^{x_m,a}$ between $0$ and $\theta_m := \tau_m \wedge h_m \wedge T_1$, where $T_1$ denotes the first jump time of $(X,I)$.
Using the identification \eqref{Sec:PDP_ident_vdelta}, from the constraint  \eqref{Sec:PDP_BSDE_constraint} and the fact that $K$ is a nondecreasing process it follows that $\P^{x_m,a}$-a.s.,
\begin{eqnarray*}
	v(x_m) &\leqslant& e^{-\delta \theta_m}\,v(X_{\theta_m}
	)
	+ \int_0^{\theta_m} e^{-\delta r}\,f(X_{r}
	, I_r
	)\,dr\\
	&& -\int_0^{\theta_m}e^{-\delta r}\,\int_{E } (v(y)-v(X_r
	))\,\tilde{q}(dr\,dy),
\end{eqnarray*}
where $\tilde{q}(dr\,dy)= p(dr\,dy)- \lambda(X_r
,I_r
)\,Q(X_r
,I_r
,dy)\,\,dr$.
In particular
\begin{eqnarray*}
	v(x_m) \leqslant \E^{x_m,a}\left[e^{-\delta \theta_m}\,v(X_{\theta_m}
	) + \int_0^{\theta_m} e^{-\delta r}\,f(X_{r}
	, I_r 
	)\,dr\right].
\end{eqnarray*}
Equation \eqref{Sec:PDP_max_property} implies that $v \leqslant v^{\ast} \leqslant \varphi$, and therefore
\begin{eqnarray*}
	\varphi(x_m) - \gamma_m \leqslant \E^{x_m,a}\left[e^{-\delta \theta_m}\,\varphi(X_{\theta_m}) + \int_0^{\theta_m} e^{-\delta r}\,f(X_{r}
	, I_r
	)\,dr\right].
\end{eqnarray*}
At this point, applying It\^o's formula to $e^{-\delta r}\,\varphi(X_{r} 
)$ between $0$ and $\theta_m$, we get
\begin{equation}\label{Sec:PDP_sub_sol_int_ineq}
	-\frac{\gamma_m}{h_m}  + \E^{x_m,a}\left[\int_0^{\theta_m} \frac{1}{h_m}\,e^{-\delta r}\,[\delta\,\varphi(X_{r}
	) -\mathcal L^{I_r
	} \varphi(X_{r}
	)-f(X_{r}
	, I_r
	)]\,dr\right]\leqslant 0,
\end{equation}
where $\mathcal L^{I_r
} \varphi (X_{r}
) = \int_{E}(\varphi (y)-\varphi (X_{r}
))\,\lambda(X_{r}
,I_r
)\,Q(X_{r}
,I_r
,dy)$.
Now we notice that,  $\P^{x_m,a}$-a.s., $(X_r,I_r) = (\phi(r,x_m),a)$ for $r \in [0,\,\theta_m]$.	Taking into account the continuity of the map $(y,b)\mapsto \delta\,\varphi(y)-\mathcal L^{b}\varphi(y)-f(y,b)$, we see that for any $\varepsilon>0$,
\begin{equation}\label{Sec:PDP_sub_sol_int_ineq2}
	-\frac{\gamma_m}{h_m}  + (\varepsilon + \delta\,\varphi(x_m) -\mathcal L^{a} \varphi(x_m)-f(x_m,a)) \,\E^{x_m,a}\left[\frac{\theta_m\,e^{-\delta\,\theta_m}}{h_m}\right]\leqslant 0,
\end{equation}
Let $f_{T_1}(s)$ 
denote the distribution  of $T_1$ under $\P^{x_m,a}$, see \eqref{Sec:PDP_Abis}.
Taking $m > M$, we have
\begin{align}
	\E^{x_m,a}\left[\frac{g(\theta_m)}{h_m}\right]&= \frac{1}{h_m}\int_0^{h_m}\!\!\!\!\!
	s\,e^{-\delta\,s}\,f_{T_1}(s)\,ds + \frac{h_m\,e^{-\delta\,h_m}}{h_m}\,\P^{x_m,a}[T_1> h_m]\nonumber\\
	&= \frac{1}{h_m}\int_0^{h_m}\!\!\!\!\!
	s\,e^{-\delta\,s}\,(\lambda(\phi(r,x_m),a)+ \lambda_0(A))\,e^{-\int_0^s(\lambda(\phi(r,x_m),a)+ \lambda_0(A))\,dr}\,ds\nonumber\\
	& + e^{-\delta\,h_m}\,e^{-\int_0^{h_m}(\lambda(\phi(r,x_m),a)+ \lambda_0(A))\,dr}.\label{Sec:PDP_2_terms}
\end{align}
By the boundedness of $\lambda$ and $\lambda_0$, it is easy to see that the two terms in the right-hand side of \eqref{Sec:PDP_2_terms} converge respectively to zero and one when $m$ goes to infinity. 
Thus, passing into the limit in \eqref{Sec:PDP_sub_sol_int_ineq2} as $m$ goes to infinity,
we obtain
$$
\delta\,\varphi(\bar{x}) -\mathcal L^{a} \varphi(\bar{x})-f(\bar{x}, a) \leqslant 0.
$$
From the arbitrariness of  $a \in A$ we  conclude that $v$ is a viscosity  subsolution to \eqref{Sec:PDP_HJB} in the sense of Definition \ref{Sec:PDP_Def_viscosity_sol_HJB}.
\endproof

\paragraph{Proof of the viscosity supersolution property to \eqref{Sec:PDP_HJB}.}
\begin{proposition}\label{Sec:PDP_Prop_visc_supersol_property_vdelta}
	Let assumptions \textup{\textbf{(H$\textup{h$\lambda$Q}$)}},
	\textup{\textbf{(H$\lambda_0$)}}, 
	and  \textup{\textbf{(H$\textup{f}$)}} hold. Then, the function $v$ in \eqref{Sec:PDP_redefine_vdelta} is a viscosity supersolution to  \eqref{Sec:PDP_HJB}.
\end{proposition}
\proof
Let $\bar{x} \in E$,  and let $\varphi \in C^1(E)$ be a test function such that
\begin{equation}\label{Sec:PDP_min_property}
	0= (v_{\ast}-\varphi)(\bar{x})=\min_{x \in E} (v_{\ast}-\varphi)(x).
\end{equation}
Notice that can assume w.l.o.g. that $\bar{x}$ is strict minimum of $v_{\ast}-\varphi$. As a matter of fact, one can subtract to $\varphi$ a positive cut-off function which behaves as $|x-\bar{x}|^2$ when $|x-\bar{x}|^2$ is  small, and that regularly converges to $1$ as  $|x-\bar{x}|^2$ increases to $1$.

Then, for every $\eta >0$, 
we can define
\begin{equation}\label{Sec:PDP_min_beta}
	0 < \beta(\eta) := \inf_{x \notin B(\bar{x},\eta)} (v_{\ast}-\varphi)(x).
\end{equation}
We will show the result by contradiction. Assume thus that
$$
H^{\varphi}(\bar{x}, \varphi, \nabla \varphi) < 0.
$$
Then by the continuity of $H$, there exists $\eta >0$, $\beta(\eta) >0$ and $\varepsilon \in (0,\,\beta(\eta) \delta]$ such that
$$
H^{\varphi}(y, \varphi, \nabla \varphi) \leqslant -\varepsilon,
$$
for all $y \in B(\bar{x}, \eta)=\{y \in E: |\bar{x}-y|<\eta\}$.
By definition of $v_{\ast}(\bar{x})$, there exists a sequence $(x_m)_m$ taking values in $B(\bar{x}, \eta)$ such that
$$
x_m \rightarrow \bar{x} \,\,\textup{and}\,\,v(x_m) \rightarrow v_{\ast}(\bar{x})
$$
when $m$ goes to infinity.
By the continuity of $\varphi$ and by \eqref{Sec:PDP_min_property} it follows that
$$
\gamma_m := v(x_m)- \varphi(x_m) \rightarrow 0,
$$
when $m$ goes to infinity.
Let us fix $T >0$ and let us define $\theta := \tau \wedge T$, where $\tau =\inf \{t \geqslant 0: X_{t} \notin B(\bar{x}, \eta)\}$.


At this point, let us fix  $a \in A$, and consider   the  solution  $Y^{n,x_m,a,\delta}$ to the penalized \eqref{Sec:PDP_BSDE_penalized}, under the probability $\P^{x_m,a}$.
Notice that
$$
\P^{x_m,a}\{\tau=0\}= \P^{x_m,a}\{X_0\notin B(\bar x, \eta)\}= 0.
$$
We apply the It\^o formula to $e^{-\delta t} \,Y_t^{n,x_m,a}$ between $0$ and 
$\theta$.
Then, proceeding as   in the proof of Lemma \ref{Sec:PDP_Lemma_rep_Y_n} 
we   get the following inequality:
\begin{eqnarray}\label{Sec:PDP_ineq_Yn_nu}
	Y_{0}^{n,x_m,a} \geqslant \inf_{\nu \in \mathcal{ V}^n} \E^{x_m,a}_{\nu}\left[e^{-\delta \theta 
	}\,Y_{\theta 
}^{n,x_m,a} + \int_0^{\theta 
} e^{-\delta r}\,f(X_{r}, I_r)\,dr\right].
\end{eqnarray}
Since $Y^{n,x_m,a}$ converges decreasingly to the maximal solution $Y^{x_m,a}$ to the constrained BSDE \eqref{Sec:PDP_BSDE}-\eqref{Sec:PDP_BSDE_constraint}, and recalling the identification \eqref{Sec:PDP_ident_vdelta}, \eqref{Sec:PDP_ineq_Yn_nu} leads to the corresponding  inequality for $v(x_m)$:
\begin{eqnarray*}
	v(x_m) \geqslant \inf_{\nu \in \mathcal{ V}} \E^{x_m,a}_{\nu}\left[e^{-\delta \theta 
	}\,v(X_{\theta 
}) + \int_0^{\theta 
} e^{-\delta r}\,f(X_{r}, I_r)\,dr\right].
\end{eqnarray*}
In particular, there exists a strictly  positive, predictable and bounded function $\nu^m$ such that
\begin{eqnarray}\label{Sec:PDP_ineq_num}
	v(x_m) \geqslant  \E^{x_m,a}_{\nu_m}\left[e^{-\delta \theta 
	}\,v(X_{\theta 
}) + \int_0^{\theta 
} e^{-\delta r}\,f(X_{r}, I_r)\,dr\right] - \frac{\varepsilon}{2 \delta}.
\end{eqnarray}
Now, from equation \eqref{Sec:PDP_min_property}
and \eqref{Sec:PDP_min_beta} we get
\begin{eqnarray*}
	\varphi(x_m) + \gamma_m \geqslant  \E^{x_m,a}_{\nu_m}\left[e^{-\delta \theta 
	}\,\varphi(X_{\theta 
}) + \beta\,e^{-\delta \theta 
}\,
\one_{\{\tau 
	\leqslant T\}}  + \int_0^{\theta 
} e^{-\delta r}\,f(X_{r}, I_r)\,dr\right] - \frac{\varepsilon}{2\,\delta}.
\end{eqnarray*}
At this point, applying It\^o's formula to $e^{-\delta r}\,\varphi(X_{r})$ between $0$ and $\theta 
$, we get
\begin{align}
	{\gamma_m}  &+  \E^{x_m,a}_{\nu_m}\left[\int_0^{\theta 
	} e^{-\delta r}\,[\delta\,\varphi(X_{r}) -\mathcal L^{I_r} \varphi(X_{r})-f(X_{r}, I_r)]\,dr -\beta\,e^{-\delta \theta 
}\,
\one_{\{\tau 
	\leqslant T\}}\right]\nonumber \\
&+ \frac{\varepsilon}{2}\geqslant 0,\label{Sec:PDP_super_sol_int_ineq}
\end{align}
where $\mathcal L^{I_r} \varphi (X_{r}) = \int_{E}(\varphi (y)-\varphi (X_{r}))\,\lambda(X_{r},I_r)\,Q(X_{r},I_r,dy)$.
Noticing that, for $0 \leqslant r \leqslant \theta 
$, 
\begin{eqnarray*}
	\delta\,\varphi(X_{r}) -\mathcal L^{I_r} \varphi(X_{r})-f(X_{r}, I_r) &\leqslant& \delta\,\varphi(X_{r}) -\inf_{b \in A}\{\mathcal L^{b} \varphi(X_{r})-f(X_{r},b)\}\\
	&=& H^{\varphi}(X_{r}, \varphi, \nabla \varphi)\\
	&\leqslant& -\varepsilon,
\end{eqnarray*}
from \eqref{Sec:PDP_super_sol_int_ineq} we obtain
\begin{eqnarray*}
	0 &\leqslant& {\gamma_m} + \frac{\varepsilon}{2\,\delta} +  \E^{x_m,a}_{\nu_m}\left[-\varepsilon \int_0^{\theta 
	} e^{-\delta r}\,dr -\beta\,e^{-\delta \theta 
}\,
\one_{\{\tau 
	\leqslant T\}}\right] \\
&=& {\gamma_m} - \frac{\varepsilon}{2\,\delta} +  \E^{x_m,a}_{\nu_m}\left[\left(\frac{\varepsilon}{\delta}- \beta\right)e^{-\delta \theta 
} \one_{\{\tau 
\leqslant T\}} + \frac{\varepsilon}{\delta} e^{-\delta \theta 
}\,\one_{\{\tau 
> T\}}\right]\\
&\leqslant& {\gamma_m} - \frac{\varepsilon}{2\,\delta} +  \frac{\varepsilon}{\delta}\,\E^{x_m,a}_{\nu_m}\left[e^{-\delta \theta 
}\,\one_{\{\tau 
> T\}}\right]\\
&=& {\gamma_m} - \frac{\varepsilon}{2\,\delta} +  \frac{\varepsilon}{\delta}\,\E^{x_m,a}_{\nu_m}\left[e^{-\delta T}\,\one_{\{\tau 
	> T\}}\right]\\
&\leqslant & {\gamma_m} - \frac{\varepsilon}{2\,\delta} +  e^{-\delta T}.
\end{eqnarray*}
Letting $T$ and $m$ go to infinity
we achieve the contradiction:
$0 \leqslant - \frac{\varepsilon}{2\,\delta}$.
\endproof
%

\appendix
\renewcommand\thesection{Appendix}
\section{}
\renewcommand\thesection{\Alph{subsection}}
\renewcommand\thesubsection{\Alph{subsection}}
\subsection{Proof of Proposition
	\ref{Sec:PDP_P_J_no_dip_a}}\label{Sec:PDP_A_proof_Prop}

We start by giving a technical result. In the sequel, $\Pi^{n_1,n_2}$ and $\Gamma^{n_1,n_2}$ will denote respectively the random sequences
$(T_{n_1},E_{n_1},A_{n_1},T_{n_1+1},E_{n_1+1},A_{n_1+1},...,T_{n_2},E_{n_2},A_{n_2})$ and \\ $(T_{n_1},A_{n_1},T_{n_1+1},A_{n_1+1},...,T_{n_2},A_{n_2})$, $n_1, n_2 \in \N \setminus \{0\}$, $n_1 \leq n_2$, where $(T_k,E_k,A_k)_{k \geq 1}$ are the random variables introduced in Section \ref{Sec:PDP_Section_control_rand}.
\begin{lemma}\label{Sec:PDP_L_aux_no_dip_a}
	Assume that Hypotheses \textup{\textbf{(H$\textup{h$\lambda$Q}$)}}, \textup{\textbf{(H$\lambda_0$)}} 
	and  \textup{\textbf{(H$\textup{f}$)}} hold. Let $\nu^n : \Omega \times \R_+ \times (\R_+ \times A)^n \times A \rightarrow (0,\,\infty)$, $n >1$ (resp. $\nu^0 : \Omega \times \R_+  \times A \rightarrow (0,\,\infty)$), be some $\mathcal P \otimes \mathcal{B}((\R_+ \times A)^n) \otimes \mathcal A$-measurable maps,  uniformly bounded with respect to $n$ (resp.  a bounded  $\mathcal P  \otimes \mathcal A$-measurable map). Let moreover $g: \Omega \times A\rightarrow (0,\,\infty)$ be a bounded $\mathcal A$-measurable map.
Set
	\begin{align}
	\nu_t(b) &= 
	\nu^0_t(b)\,\one_{\{t \leqslant T_1\}} + \sum_{n =1}^{\infty} \nu^n_t( 
	\Gamma^{1,n},b)\,\one_{\{T_{n} < t \leqslant T_{n+1}\}},\label{Sec:PDP_nu_pred_form}\\
	\nu'_t(b) &= 
	g(b)\,\one_{\{t \leqslant T_1\}} + \nu^0_t(b)\,\one_{\{ T_1 < t \leqslant T_2\}}  + \sum_{n =2}^{\infty} \nu_t^{n-1}( 
	\Gamma^{2,n},b)\,\one_{\{T_{n} < t \leqslant T_{n+1}\}}.\label{Sec:PDP_nu_vareps_pred_form}
	\end{align}
	
	Fix $x \in E$, $a,a' \in A$.
	Then, for every $n >1$,  and  for every $\mathcal{B}((\R_+ \times E\times A)^n)$-measurable function $F:(\R_+ \times E\times A)^n \rightarrow \R$,
	\begin{align}\label{Sec:PDP_conjecture_PDMP}
	{\mathbb E}^{x,a'}_{\nu'} \left[ 
	F( 
	\Pi^{1,n})|\mathcal F_{T_1}\right]= \frac{\mathbb E^{x,a}_{\nu}\left[ \one_{\{T_1 >\tau\}}\,
		F(\tau,\chi,\xi,\Pi^{1,n-1}
		)
		\right]}{\P^{x,a}_{\nu}(T_1 >\tau)}\bigg|_{\tau = T_1,\chi=X_1,\, \xi = A_1}.
	\end{align}
	
\end{lemma}
\begin{remark}
		$\P^{x,a}_{\nu}$ (resp. $\P^{x,a'}_{\nu'}$)  is  the unique  probability measure on $(\Omega, \mathcal F_{\infty})$  under which the random measure  $\tilde{p}^\nu$ (resp. $\tilde{p}^{\nu'}$) in \eqref{Sec:PDP_dual_comp} is the compensator of the measure $p$ in \eqref{Sec:PDP_p_dual} on $(0,\,\infty)\times E \times A$, see Proposition \ref{P_Prob_infty}.
\end{remark}
\proof

Taking into account 
\eqref{Sec:PDP_A_kbis}, \eqref{Sec:PDP_B_kbis}, and 
\eqref{Sec:PDP_nu_vareps_pred_form},  we have:  for all $r \geqslant T_1$,
\begin{align}\label{Sec:PDP_T2A2E2_density}
&\P^{x,a'}_{\nu'}\left[T_{2} > r,E_{2} \in F, A_{2} \in C|\mathcal{F}_{T_1} \right]\nonumber\\
&=\int_{r}^{\infty}\int_F\exp\left(-\int_{T_1}^s\lambda(\phi(t-T_1,E_1,A_1),A_1)\,dt-\int_{T_1}^{s}\int_A\nu^0_t(b)\,\lambda_0(db)\,dt\right)\cdot\nonumber\\
&\cdot \lambda(\phi(s-T_1,E_1,A_1),A_1)\,Q(\phi(s-T_1,E_1,A_1),A_1,dy)\,ds\nonumber\\
&+\int_{r}^{\infty}\int_C\exp\left(-\int_{T_1}^s\lambda(\phi(t-T_1,E_1,A_1),A_1)\,dt-\int_{T_1}^{s}\int_A\nu^0_t(b)\,\lambda_0(db)\,dt\right)\,\nu^0_s(b)\,\lambda_0(db)\,ds,
\end{align}
and, for all $r \geqslant T_n$,  $n >2$,
\begin{align}\label{Sec:PDP_TnAnEn_density}
&\P^{x,a}_{\nu'}\left[T_{n+1} > r,E_{n+1} \in F, A_{n+1} \in C|\mathcal{F}_{T_n} \right]\nonumber\\
&=\int_{r}^{\infty}\int_F\exp\bigg(-\int_{T_n}^s\lambda(\phi(t-T_n,E_n,A_n),A_n)\,dt-\int_{T_n}^{s}\int_A\nu^{n-1}_t( 
\Gamma^{2,n},b)\,\lambda_0(db)\,dt\bigg)\cdot\nonumber\\
&\cdot \lambda(\phi(s-T_n,E_n,A_n),A_n)\,Q(\phi(s-T_n,E_n,A_n),A_n,dy)\,ds\nonumber\\
&+\int_{r}^{\infty}\int_C\exp\bigg(-\int_{T_n}^s\lambda(\phi(t-T_n,E_n,A_n),A_n)\,dt-\int_{T_n}^{s}\int_A\nu^{n-1}_t( 
\Gamma^{2,n},b)\,\lambda_0(db)\,dt\bigg) \cdot\nonumber\\
& \cdot
\nu^{n-1}_s( 
\Gamma^{2,n},b)\,\lambda_0(db)\,ds.
\end{align}

We will  prove identity \eqref{Sec:PDP_conjecture_PDMP} by induction.
Let us start by showing that \eqref{Sec:PDP_conjecture_PDMP} holds in  the case $n=2$, namely that,  for every $\mathcal{B}((\R_+ \times E\times A)^2) $-measurable function $F:(\R_+ \times E\times A)^2 \rightarrow \R$,
\begin{align}\label{Sec:PDP_conjecturePDP_n=2}
{\mathbb E}^{x,a'}_{\nu'} \left[  F( 
\Pi^{1,2})|\mathcal F_{T_1}\right]
= \frac{{\mathbb E}^{x,a}_{\nu}\left[ \one_{\{T_1 >\tau\}}\,F(\tau,\chi,\xi, 
	\Pi^{1,1})\right]}{\P^{x,a}_{\nu}(T_1 >\tau)}\bigg|_{\tau = T_1,\chi=X_1,\, \xi = A_1}.
\end{align}
From \eqref{Sec:PDP_T2A2E2_density} we get
\begin{align*}
&{\mathbb E}^{x,a'}_{\nu'} \left[ F( 
\Pi^{1,2})| \mathcal{F}_{T_1} \right]={\mathbb E}^{x,a'}_{\nu'} \left[ F(T_1,E_1,A_1, T_2,E_2,A_2)| \mathcal{F}_{T_1} \right]\\
&= \int_{T_1}^{\infty}\int_{E}F(T_1,E_1,A_1,s,y,A_1)\,\exp\left(-\int_{T_1}^s\lambda(\phi(t-T_1,E_1,A_1),A_1)\,dt-\int_{T_1}^{s}\int_A\nu^0_t(b)\,\lambda_0(db)\,dt\right)\cdot\nonumber\\
&\cdot \lambda(\phi(s-T_1,E_1,A_1),A_1)\,Q(\phi(s-T_1,E_1,A_1),A_1,dy)\,ds\nonumber\\
&+\int_{T_1}^{\infty}\int_A F(T_1,E_1,A_1,s,\phi(s-T_1,E_1,A_1),b)\cdot\\
&\cdot\exp\left(-\int_{T_1}^s\lambda(\phi(t-T_1,E_1,A_1),A_1)\,dt-\int_{T_1}^{s}\int_A\nu^0_t(b)\,\lambda_0(db)\,dt\right)\,\nu^0_s(b)\,\lambda_0(db)\,ds.
\end{align*}
On the other hand,
$$
\P^{x,a}_\nu(T_1 >\tau)= \exp\left(-\int_{0}^\tau\lambda(\phi(t-\tau,\chi,\xi),\xi)\,dt-\int_{0}^{\tau}\int_A\nu^0_t(b)\,\lambda_0(db)\,dt\right),
$$
and
\begin{align*}
&{\mathbb E}^{x,a}_\nu \left[  \one_{\{T_1 >\tau\}}\,F(\tau,\chi,\xi,\Pi^{1,1}) \right]={\mathbb E}^{x,a}_\nu \left[  \one_{\{T_1 >\tau\}}\,F(\tau,\chi,\xi,T_1,E_1,A_1) \right]\\
&=
\int_\tau^\infty\int_{E} \one_{\{s >\tau\}}\,F(\tau,\chi,\xi,s,y,\xi)\,\exp\left(-\int_{0}^s\lambda(\phi(t-\tau,\chi,\xi),\xi)\,dt-\int_{0}^{s}\int_A\nu^0_t(b)\,\lambda_0(db)\,dt\right)\cdot\\
& \cdot\lambda(\phi(s-\tau,\chi,\xi),\xi)\,Q(\phi(s-\tau,\chi,\xi),\xi,dy)\,ds\\
&+\int_\tau^\infty\int_A \one_{\{s >\tau\}}\,F(\tau,\chi,\xi,s,\phi(s-\tau,\chi,\xi),b)\cdot\\
&\cdot \exp\left(-\int_{0}^s\lambda(\phi(t-\tau,\chi,\xi),\xi)\,dt-\int_{0}^{s}\int_A\nu^0_t(b)\,\lambda_0(db)\,dt\right)\,\nu^0_s(b)\,\lambda_0(db)\,ds.
\end{align*}
Therefore,
\begin{align*}
&\frac{{\mathbb E}^{x,a}_\nu \left[  \one_{\{T_1 >\tau\}}\,F(\tau,\chi,\xi,\Pi^{1,1}) \right]}{\P^{x,a}_\nu(T_1 >\tau)}\\
&= \exp\left(\int_{0}^\tau\lambda(\phi(t-\tau,\chi,\xi),\xi)\,dt+\int_{0}^{\tau}\int_A\nu^0_t(b)\,\lambda_0(db)\,dt\right)\cdot\\
&\cdot\int_\tau^\infty\int_{E} \one_{\{s >\tau\}}\,F(\tau,\chi,\xi,s,y,\xi)\, \exp\left(-\int_{0}^s\lambda(\phi(t-\tau,\chi,\xi),\xi)\,dt-\int_{0}^{s}\int_A\nu^0_t(b)\,\lambda_0(db)\,dt\right)\cdot\\
& \cdot\lambda(\phi(s-\tau,\chi,\xi),\xi)\,Q(\phi(s-\tau,\chi,\xi),\xi,dy)\,ds\\
&+\exp\left(\int_{0}^\tau\lambda(\phi(t-\tau,\chi,\xi),\xi)\,dt+\int_{0}^{\tau}\int_A\nu^0_t(b)\,\lambda_0(db)\,dt\right)\cdot\\
&\cdot\int_\tau^\infty\int_A \one_{\{s >\tau\}}\,F(\tau,\chi,\xi,s,\phi(s-\tau,\chi,\xi),b)\cdot\\
&\cdot \exp\left(-\int_{0}^s\lambda(\phi(t-\tau,\chi,\xi),\xi)\,dt-\int_{0}^{s}\int_A\nu^0_t(b)\,\lambda_0(db)\,dt\right)\,\nu^0_s(b)\,\lambda_0(db)\,ds\\
&=\int_\tau^\infty\int_{E} \one_{\{s >\tau\}}\,F(\tau,\chi,\xi,s,y,\xi)\,\exp\left(-\int_{\tau}^s\lambda(\phi(t-\tau,\chi,\xi),\xi)\,dt-\int_{\tau}^{s}\int_A\nu^0_t(b)\,\lambda_0(db)\,dt\right)\cdot\\
& \cdot\lambda(\phi(s-\tau,\chi,\xi),\xi)\,Q(\phi(s-\tau,\chi,\xi),\xi,dy)\,ds\\
&+\int_\tau^\infty\int_A \one_{\{s >\tau\}}\,F(\tau,\chi,\xi,s,\phi(s-\tau,\chi,\xi),b)\cdot\\
& \cdot \exp\left(-\int_{\tau}^s\lambda(\phi(t-\tau,\chi,\xi),\xi)\,dt-\int_{\tau}^{s}\int_A\nu^0_t(b)\,\lambda_0(db)\,dt\right)\,\nu^0_s(b)\,\lambda_0(db)\,ds,
\end{align*}
and \eqref{Sec:PDP_conjecturePDP_n=2} follows.

Assume now 
that \eqref{Sec:PDP_conjecture_PDMP} holds for  $n-1$,
namely that,  for every $\mathcal{B}((\R_+ \times E\times A)^{n-1}) $-measurable function $F:(\R_+ \times E\times A)^{n-1} \rightarrow \R$,
\begin{align}\label{Sec:PDP_conjecture_PDMP_induction}
{\mathbb E}^{x,a'}_{\nu'}
\left[ 
F( 
\Pi^{1,n-1}
)|\mathcal F_{T_1}\right]
= \frac{\mathbb E^{x,a}_{\nu}\left[ \one_{\{T_1 >\tau\}}\,
	F(\tau,\chi,\xi, \Pi^{1,n-2}
	)\right]}{\P^{x,a}_{\nu}(T_1 >\tau)}\bigg|_{\tau = T_1,\chi=X_1,\, \xi = A_1}.
\end{align}
We have to prove that \eqref{Sec:PDP_conjecture_PDMP_induction} implies that, for every $\mathcal{B}((\R_+ \times E\times A)^{n}) $-measurable function $F:(\R_+ \times E\times A)^{n} \rightarrow \R$, 
\begin{align}\label{Sec:PDP_conjecture_PDMP_induction2}
{\mathbb E}^{x,a'}_{\nu'}
\left[ 
F( 
\Pi^{1,n}
)|\mathcal F_{T_1}\right]
= \frac{\mathbb E^{x,a}_{\nu}\left[ \one_{\{T_1 >\tau\}}\,
	F(\tau,\chi,\xi, \Pi^{1,n-1}
	)\right]}{\P^{x,a}_{\nu}(T_1 >\tau)}\bigg|_{\tau = T_1,\chi=X_1,\, \xi = A_1}.
\end{align}
Using  \eqref{Sec:PDP_TnAnEn_density}, 
we get
\begin{align}\label{Sec:PDP_phi_per_passo_induttivo}
&{\mathbb E}^{x,a'}_{\nu'}  \left[ 
F( 
\Pi^{1,n})| \mathcal{F}_{T_1} \right]\nonumber\\
& = {\mathbb E}^{x,a'}_{\nu'} \left[	{\mathbb E}^{x,a'}_{\nu'} \left[ 
F( 
\Pi^{1,n})| \mathcal{F}_{T_{n-1}} \right]\big| \mathcal{F}_{T_1}\right]\nonumber\\
&= {\mathbb E}^{x,a'}_{\nu'} \bigg[
\int_{T_{n-1}}^{\infty}\int_{E} 
F( 
\Pi^{1,n-1},s,y,A_{n-1})\cdot\nonumber\\
& \cdot\exp\bigg(-\int_{T_{n-1}}^s\lambda(\phi(t-T_{n-1},E_{n-1},A_{n-1}),A_{n-1})\,dt
-\int_{T_{n-1}}^{s}\int_A\nu^{n-2}_t 
( 
\Gamma^{1,n-1},b)\,\lambda_0(db)\,dt\bigg)\cdot\nonumber\\
&\cdot \lambda(\phi(s-T_{n-1},E_{n-1},A_{n-1}),A_{n-1})\,Q(\phi(s-T_{n-1},E_{n-1},A_{n-1}),A_{n-1},dy)\,ds 
\nonumber\\
& +
\int_{T_{n-1}}^{\infty}\int_{A} 
F( 
\Pi^{1,n-1},s,\phi(s-T_{n-1},E_{n-1},A_{n-1}),b)\cdot\nonumber\\
&\cdot\exp\bigg(-\int_{T_{n-1}}^s\lambda(\phi(t-T_{n-1},E_{n-1},A_{n-1}),A_{n-1})\,dt -\int_{T_{n-1}}^{s}\int_A\nu^{n-2}_t 
( 
\Gamma^{1,n-1},b)\,\lambda_0(db)\,dt\bigg)\cdot\nonumber\\
&\cdot \nu^{n-2}_s 
( 
\Gamma^{1,n-1},b)\,\lambda_0(db)\,ds
\bigg| \mathcal{F}_{T_1}\bigg].
\end{align}
 At this point we observe that the  term in the conditional expectation in the right-hand side of \eqref{Sec:PDP_phi_per_passo_induttivo} only depends on the random sequence $\Pi^{1,n-1}$. 
 For any sequence of  random variables   $(S_i, W_i, V_i)_{i\in [1,n-1]}$ with values in $([0,\,\infty) \times E \times A)^{n-1}$, $S_{i-1} \leq S_{i}$ for every $i\in [1,n-1]$,  we set
 \begin{align*} 
 &\psi(S_1, W_1, V_1,..., S_{n-1}, W_{n-1}, V_{n-1}):=\nonumber\\
 &\int_{S_{n-1}}^{\infty}\int_{E} F(S_1,W_1,...,V_{n-1},S_{n-1},W_{n-1},s,y,A_{n-1})
 \cdot\nonumber\\
 & \cdot\exp\bigg(-\int_{S_{n-1}}^s\lambda(\phi(t-S_{n-1},W_{n-1},V_{n-1}),V_{n-1})\,dt \nonumber\\
 &\qquad \quad -\int_{S_{n-1}}^{s}\int_A\nu^{n-2}_t (S_1,V_1,...,S_{n-1},V_{n-1},b)
 \,\lambda_0(db)\,dt\bigg)\cdot\nonumber\\
 &\cdot \lambda(\phi(s-S_{n-1},W_{n-1},V_{n-1}),V_{n-1})\,Q(\phi(s-S_{n-1},W_{n-1},V_{n-1}),V_{n-1},dy)\,ds 
 \nonumber\\
 & +
 \int_{S_{n-1}}^{\infty}\int_{A} 
 F( S_1,W_1,V_1, ...,S_{n-1},W_{n-1},V_{n-1},
 ,s,\phi(s-S_{n-1},W_{n-1},V_{n-1}),b)\cdot\nonumber\\
 &\cdot\exp\bigg(-\int_{S_{n-1}}^s\lambda(\phi(t-S_{n-1},W_{n-1},V_{n-1}),V_{n-1})\,dt  \nonumber\\
 &\qquad \quad -\int_{S_{n-1}}^{s}\int_A\nu^{n-2}_t (S_1,V_1,...,S_{n-1},V_{n-1},b)
 \,\lambda_0(db)\,dt\bigg)\cdot\nonumber\\
 &\cdot \nu^{n-2}_s (S_1,V_1,...,S_{n-1},V_{n-1},b)
 \,\lambda_0(db)\,ds.
 \end{align*}
 Identity \eqref{Sec:PDP_phi_per_passo_induttivo}  can then  be rewritten as
 \begin{align}
 {\mathbb E}^{x,a'}_{\nu'}  \left[ 
 F( 
 \Pi^{1,n})| \mathcal{F}_{T_1} \right]
 ={\mathbb E}^{x,a'}_{\nu'}\left[	
 \psi( 
 \Pi^{1,n-1})\Big| \mathcal{F}_{T_1}\right].
 \end{align}
Then, by applying  the inductive step \eqref{Sec:PDP_conjecture_PDMP_induction}, we get
\begin{align}\label{Sec:PDP_change_prob}
&{\mathbb E}^{x,a'}_{\nu'}  \left[ 
F( 
\Pi^{1,n})| \mathcal{F}_{T_1} \right]\nonumber\\
& ={\mathbb E}^{x,a'}_{\nu'}\left[	
\psi( 
\Pi^{1,n-1})\Big| \mathcal{F}_{T_1}\right]\nonumber\\
& =(\P^{x,a}_{\nu}[T_1 >\tau])^{-1}\,\mathbb E^{x,a}_{\nu}\bigg[ \one_{\{T_1 >\tau\}}\,
\psi(\tau,\chi,\xi, 
\Pi^{1,n-2})
\bigg]\bigg|_{\tau = T_1,\chi=X_1,\, \xi = A_1}.
\end{align}
Since
\begin{align*}
&\psi(\tau,\chi,\xi, 
\Pi^{1,n-2}) \\
&=\int_{T_{n-2}}^\infty\int_{E}
F(\tau,\chi,\xi, 
\Pi^{1,n-2},s,y,A_{n-2})
\cdot\nonumber\\
& \cdot\exp\bigg(-\int_{T_{n-2}}^s\lambda(\phi(t-T_{n-2},E_{n-2},A_{n-2}),A_{n-2})\,dt -\int_{T_{n-2}}^{s}\int_A\nu^{n-2}_t
( 
\Gamma^{1,n-2},b)
\,\lambda_0(db)\,dt\bigg)\cdot\nonumber\\
& \cdot \lambda(\phi(s-T_{n-2},E_{n-2},A_{n-2}),A_{n-2})\,Q(\phi(s-T_{n-2},E_{n-2},A_{n-2}),A_{n-2},dy)\,ds\\
&+\int_{T_{n-2}}^\infty\int_{A}
F(\tau,\chi,\xi, 
\Pi^{1,n-2},s,\phi(s-T_{n-2},E_{n-2},A_{n-2}),b)
\cdot\nonumber\\
&\cdot \exp\bigg(-\int_{T_{n-2}}^s\lambda(\phi(t-T_{n-2},E_{n-2},A_{n-2}),A_1)\,dt-\int_{T_{n-2}}^{s}\int_A\nu^{n-2}_t
( 
\Gamma^{1,n-2},b)\,\lambda_0(db)\,dt\bigg)\cdot\nonumber\\
&\cdot\nu^{n-2}_s
( 
\Gamma^{1,n-2},b)
\,\lambda_0(db)\,ds\\
&= \E^{x,a}_{\nu}[F(\tau, \chi, \xi, \Pi^{1,n-1})|\mathcal{F}_{T_{n-2}}],
\end{align*}
identity \eqref{Sec:PDP_change_prob} reads
\begin{align}
&{\mathbb E}^{x,a'}_{\nu'}  \left[ 
F( 
\Pi^{1,n})| \mathcal{F}_{T_1} \right]\nonumber\\
& =(\P^{x,a}_{\nu}[T_1 >\tau])^{-1}\,\mathbb E^{x,a}_{\nu}\bigg[ \one_{\{T_1 >\tau\}}\,
\E^{x,a}_{\nu}[F(\tau, \chi, \xi, \Pi^{1,n-1})|\mathcal{F}_{T_{n-2}}]
\bigg]\bigg|_{\tau = T_1,\chi=X_1,\, \xi = A_1}\nonumber\\
& =\frac{\mathbb E^{x,a}_{\nu}\left[\one_{\{T_1 >\tau\}}\,
	F(\tau,\chi,\xi, 
	\Pi^{1,n-1})
	\right]}{\P^{x,a}_{\nu}(T_1 >\tau)}\bigg|_{\tau = T_1,\,\chi=E_1,\,\xi=A_1}.\label{Sec:PDP_conjecturePDP_n=3}
\end{align}
This concludes the proof of the Lemma.
\qed


\paragraph{Proof of Proposition
	\ref{Sec:PDP_P_J_no_dip_a}}
 We start by noticing that,
 \begin{equation*}
 J(x,a,\nu)={\mathbb E}^{x,a}_{\nu}\left[ 
 F(T_1,E_1,A_1,T_2,E_2,A_2,...)
 \right],
 \end{equation*}
 where
 \begin{align}\label{Sec:PDP_def_psi}
 & 
 F(T_1,E_1,A_1,T_2,E_2,A_2,...)\nonumber\\
 &= \int_0^\infty e^{-\delta t} f(X_t,I_t)\,dt\nonumber\\
 &= \int_0^{T_1} e^{-\delta t} f(\phi(t,X_0,I_0),I_0)\,dt  
 + \sum_{n=2}^{\infty}\int_{T_{n-1}}^{T_n} e^{-\delta t}f(\phi(t-T_{n-1},E_{n-1},A_{n-1}), A_{n-1})\,dt.
 \end{align}
 We aim at  constructing a sequence of controls  $(\nu^{\varepsilon})_{\varepsilon >0}$, $\nu_{\varepsilon} \in\mathcal V$ for every $\varepsilon$, such that
 \begin{align}\label{Sec:PDP_claim_PDMP2}
 J(x,a',\nu^{\varepsilon})&={\mathbb E}^{x,a'}_{\nu^{\varepsilon}}\left[ 
 F(T_1,E_1,A_1,T_2, E_2, A_2,... 
 ) \right]\nonumber\\
 &\quad  \underset{ \varepsilon \rightarrow 0^+}{\longrightarrow } \mathbb E^{x,a}_{\nu}\left[
 F(T_1,E_1,A_1, T_2, E_2, A_2, ...
 )
 \right]=	 J(x,a,\nu).
 \end{align}
 Since $\nu \in \mathcal V$, then there exists a $\P^{x,a}$-null set $N$ such that  $\nu$ admits the
 representation
 \begin{align}\label{Sec:PDP_nu_pred_form2}
 \nu_t(b) 
 &=\nu^0_t(b)\,\one_{\{t \leqslant T_1\}} + \sum_{n =1}^{\infty} \nu^n_t(T_1,A_1,T_2,A_2,...,T_n,A_n,b)\,\one_{\{T_{n} < t \leqslant T_{n+1}\}}
 \end{align}
 for all $(\omega, t) \in \Omega \times \R_+$,  $\omega \notin N$,
 for some  $\nu^n : \Omega \times \R_+ \times (\R_+ \times A)^n \times A \rightarrow (0,\,\infty)$, $n >1$ (resp. $\nu^0 : \Omega \times \R_+  \times A \rightarrow (0,\,\infty)$)  $\mathcal P \otimes \mathcal{B}((\R_+ \times A)^n) \otimes \mathcal A$-measurable maps uniformly bounded with respect to $n$ (resp. bounded $\mathcal P  \otimes \mathcal A$-measurable map), see e.g. Definition 26.3 in \cite{Da}.
 %
 %
 
 Let $\bar B(a,\varepsilon)$ be the closed ball centered in $a$ with radius $\varepsilon$.
 We notice that   $\varepsilon \mapsto \lambda_0(\bar B(a,\varepsilon))$ defines a nonegative,   right-continuous, nondecreasing funtion, satisfying  
 $$
 \lambda_0(\bar B(a,0)) = \lambda_0(\{a\}) \geq 0,\quad \lambda_0(\bar B(a,\varepsilon))>0\quad \forall \varepsilon >0.
 $$
 If $\lambda_{0}(\{a\})>0$, we set $h(\varepsilon)= \varepsilon$ for every $\varepsilon >0$. Otherwise, if $\lambda_{0}(\{a\})=0$, we define $h$ as the 
 right inverse  function of  $\varepsilon \mapsto \lambda_0(\bar B(a,\varepsilon))$, namely 
 $$
 h(p)= \inf\{\varepsilon>0: \lambda_0(\bar B(a,\varepsilon))\geq p\},\,\,p\geq 0. 
 $$
 From  Lemma 1.37  in \cite{chineseBook} the following property holds:
 \begin{equation}\label{Sec:PDP_inversePROP}
 \forall p\geq 0, 
 \quad \lambda_0(\bar B(a,h(p)))\geq p.
 \end{equation}
 
 At this point, we  introduce the following family of processes, parametrized by $\varepsilon$:
 \begin{align}\label{Sec:PDP_nuvarepsilon}
 \nu_t^{\varepsilon}(b) 
 &= \frac{1}{\varepsilon}\,\frac{1}{\lambda_0(\bar B(a,h(\varepsilon)))}\one_{\{b \in \bar B(a,h(\varepsilon))\}}\one_{\{t \leqslant T_1\}} + \nu^0_t(b)\,\one_{\{ T_1 < t \leqslant T_2\}} \nonumber \\
 & + \sum_{n =2}^{\infty} \nu_t^{n-1} 
 (T_2, A_2, ..., T_{n},A_{n},b)
 \,\one_{\{T_{n} < t \leqslant T_{n+1}\}}.
 \end{align}
 With this choice, for all $r >0$,
 \begin{align}\label{Sec:PDP_T_1A_1E_1_distr_epsilon}
 &\P^{x,a'}_{\nu^{\varepsilon}}(T_1 >r, E_1 \in F,A_1 \in C)\nonumber\\
 & = \int_r^\infty\int_F \exp\left(-\int_0^s\lambda(\phi(t,x,a'),a')\,dt-\frac{s}{\varepsilon}\right)\,\lambda(\phi(s,x,a'),a')\,Q(\phi(s,x,a'),a',dy)\,ds\nonumber\\
 &+\int_r^\infty\int_C \exp\left(-\int_0^s\lambda(\phi(t,x,a'),a')\,dt-\frac{s}{\varepsilon}\right)\,\frac{1}{\varepsilon}\,\frac{1}{\lambda_0(\bar B(a,h(\varepsilon)))}\one_{\{b \in \bar B(a,h(\varepsilon))\}}\,\lambda_0(db)\,ds.
 \end{align}
 To prove \eqref{Sec:PDP_claim_PDMP2}, it is enough to show that, for every $k >1$,
 \begin{equation}\label{Sec:PDP_claim_PDMPn}
 {\mathbb E}^{x,a'}_{\nu^{\varepsilon}}[ 
 \bar F(\Pi^{1,k} 
 )]  \underset{ \varepsilon \rightarrow 0}{\longrightarrow } \mathbb E^{x,a}_{\nu}[
 \bar F( 
 \Pi^{1,k})],
 \end{equation}
 where
 \begin{align}\label{barF}
 \bar F(S_1, W_1, V_1,..., S_k, W_k, V_k)&= \int_0^{S_1} e^{-\delta t} f(\phi(t,X_0,I_0),I_0)\,dt \nonumber\\
 &
 + \sum_{n=2}^{k}\int_{S_{n-1}}^{S_n} e^{-\delta t}f(\phi(t-S_{n-1},W_{n-1},V_{n-1}), V_{n-1})\,dt,
 \end{align}
 for any sequence of random variables $(S_n, W_n, V_n)_{n \in [1,k]}$ with values in $([0,\,\infty) \times E \times A)^n$, $S_{n-1} \leq S_{n}$ for every $n \in [1,k]$.
 As a matter of fact, the remaining term 
 $$
 R(\varepsilon,k):= \E^{x,a'}_{\nu^{\varepsilon}}\left[\int_{T_{k}}^{\infty} e^{-\delta t}f(X_t,I_t)\,dt\right]
 $$
 converges to zero, uniformly in $\varepsilon$, as $k$ goes to infinity.
 To see it, we notice that
 \begin{align}\label{Sec:PDP_estimate_rest}
 |R(\varepsilon,k)| \leq \frac{M_f}{\delta}\,\E^{x,a'}_{\nu^{\varepsilon}}\left[ e^{-\delta T_k}\right]= \frac{M_f}{\delta}\,\E^{x,a'}\left[ L^{\nu^{\varepsilon}}_{T_k}\,e^{-\delta T_k}\right],
 \end{align}
 where  $L^{\nu}$ denotes the Dol\'eans-Dade exponential local martingale defined in \eqref{Sec:PDP_Lnu}. 
 On the other hand, taking into account \eqref{Sec:PDP_nuvarepsilon} and \eqref{Sec:PDP_inversePROP}, we get 
 \begin{align*}
 \E^{x,a'}\left[ L^{\nu^{\varepsilon}}_{T_k}\,e^{-\delta T_k}\right]\leq\E^{x,a'}\left[ \frac{e^{T_1\,\lambda_0(A)}\,e^{-T_1\,\frac{1}{\varepsilon}}}{\varepsilon^2 
 }\,L^{\bar \nu}_{T_k}\,e^{-\delta T_k}\right]\leq \frac{4}{e^2}\, \E^{x,a'}\left[ \frac{e^{T_1\,\lambda_0(A)}}{T_1^2}\,L^{\bar \nu}_{T_k}\,e^{-\delta T_k}\right]
 \end{align*}	
 where   
 \begin{align*}
 \bar \nu(t,b) 
 := \one_{\{t \leqslant T_1\}} + \nu^0_t(b)\,\one_{\{ T_1 < t \leqslant T_2\}} + \sum_{n =2}^{\infty} \nu_t^{n-1} 
 (T_2, A_2, ..., T_{n},A_{n},b)
 \,\one_{\{T_{n} < t \leqslant T_{n+1}\}}.
 \end{align*}
 Since $\bar \nu \in \mathcal V$, by Proposition \ref{P_Prob_infty} there exists a unique  probability $\P^{x,a'}_{\bar \nu}$ on $(\Omega, \mathcal F_{\infty})$ such that its restriction on $(\Omega, \mathcal F_{T_k})$ is  $L_{T_k}^{\bar \nu}\,\P^{x,a'}$. Then \eqref{Sec:PDP_estimate_rest} reads 
 \begin{align}\label{Sec:PDP_estimate_rest2}
 |R(\varepsilon,k)|
 \leq\frac{4\,M_f}{\delta\,e^2}\, \E^{x,a'}_{\bar \nu}\left[ \frac{e^{T_1\,\lambda_0(A)}}{T_1^2}\,e^{-\delta T_k}\right], 
 \end{align}
 and the conclusion follows by the Lebesgue dominated convergence theorem.
 
 Let us now prove \eqref{Sec:PDP_claim_PDMPn}.
 By Lemma \ref{Sec:PDP_L_aux_no_dip_a} with $\nu'= \nu^{\varepsilon}$, taking into account
 \eqref{Sec:PDP_T_1A_1E_1_distr_epsilon},  we achieve
 \begin{align}\label{Sec:PDP_intermediate}
 &{\mathbb E}^{x,a'}_{\nu^{\varepsilon}}
 [\bar F( 
 \Pi^{1,k})] \nonumber\\
 &={\mathbb E}^{x,a'}_{\nu^{\varepsilon}}\left[ {\mathbb E}^{x,a'}_{\nu^{\varepsilon}}\left[ 
 \bar F( 
 \Pi^{1,k})
 \right|\mathcal{F}_{T_1}]\right]\nonumber\\
 &={\mathbb E}^{x,a'}_{\nu^{\varepsilon}}\left[\frac{\mathbb E^{x,a}_{\nu}\left[\one_{\{T_1 >\tau\}}\,
 	\bar F(s,y,b, 
 	\Pi^{1,k-1})
 	\right]}{\P^{x,a}_{\nu}(T_1 >\tau)}\bigg|_{s = T_1,\,y=X_1,\,b=A_1}\right]\nonumber\\
 &=\int_0^\infty\int_{E}\frac{\mathbb E^{x,a}_{\nu}\left[\one_{\{T_1 >s\}}\,
 	\bar F(s,y,a', 
 	\Pi^{1,k-1})
 	\right]}{\P^{x,a}_{\nu}(T_1 >s)}\cdot\nonumber\\ &\cdot\exp\left(-\int_0^s\lambda(\phi(t,x,a'),a')\,dt-\frac{s}{\varepsilon}\right)\,\lambda(\phi(s,x,a'),a')\,Q(\phi(s,x,a'),a',dy)\,ds\nonumber\\
 &+\int_0^\infty\int_A \frac{\mathbb E^{x,a}_{\nu}\left[\one_{\{T_1 >s\}}\,
 	\bar F(s,\phi(s,x,a'),b, 
 	\Pi^{1,k-1})
 	\right]}{\P^{x,a}_{\nu}(T_1 >s)}\cdot\nonumber\\ 	 &\cdot\exp\left(-\int_0^s\lambda(\phi(t,x,a'),a')\,dt-\frac{s}{\varepsilon}\right)\,\frac{1}{\varepsilon}\,\frac{1}{\lambda_0(\bar B(a,h(\varepsilon)))}\one_{\{b \in \bar B(a,h(\varepsilon))\}}\,\lambda_0(db)\,ds.
 \end{align}
 At this point, we set
 \begin{equation}
 \label{Sec:PDP_def_varphi}
 \varphi(s,y,b):=\frac{\mathbb E^{x,a}_{\nu}\left[\one_{\{T_1 >s\}}\,
 	\bar F(s,y,b, 
 	\Pi^{1,k-1})
 	\right]}{\P^{x,a}_{\nu}(T_1 >s)}, \quad s \in [0,\infty),\,y \in E, \,b \in A.
 \end{equation}
 Notice that, for every $(y,b) \in E \times A$,
  \begin{align*}
  \bar F(s,y,b,\Pi^{1,k-1})&= \int_0^{s} e^{-\delta t} f(\phi(t,X_0,I_0),I_0)\,dt +
  \int_{s}^{T_1} e^{-\delta t}f(\phi(t-s,y,b),b)\,dt\\
    &
  + \sum_{n=2}^{k-1}\int_{T_{n-1}}^{T_n} e^{-\delta t}f(\phi(t-T_{n-1},E_{n-1},A_{n-1}), A_{n-1})\,dt,
  \end{align*}
  so that 
 \begin{align}\label{Sec:PDP_estimate_varphi}
 |\varphi(s,y,b)| \leq \frac{M_f}{\delta}. 
 \end{align}
 Identity \eqref{Sec:PDP_intermediate} becomes
 \begin{align*}
 &{\mathbb E}^{x,a'}_{\nu^{\varepsilon}}[\bar F( 
 \Pi^{1,k})] \\
 &=\int_0^\infty\int_{E}\varphi(s,y,a')\, \exp\left(-\int_0^s\lambda(\phi(t,x,a'),a')\,dt-\frac{s}{\varepsilon}\right)\cdot\\
 &\cdot\lambda(\phi(s,x,a'),a')\,Q(\phi(s,x,a'),a',dy)\,ds\nonumber\\
 &+\int_0^\infty\int_A \varphi(s,\phi(s,x,a'),b)\, \exp\left(-\int_0^s\lambda(\phi(t,x,a'),a')\,dt-\frac{s}{\varepsilon}\right)\cdot\\
 &\cdot\frac{1}{\varepsilon}\,\frac{1}{\lambda_0(\bar B(a,h(\varepsilon)))}\one_{\{b \in \bar B(a,h(\varepsilon))\}}\,\lambda_0(db)\,ds\\
 &=:I_1(\varepsilon)+ I_2(\varepsilon).
 \end{align*}
 Using  the change of variable $s=\varepsilon\,z$, we have
 \begin{align*}
 I_1(\varepsilon)&=\int_0^\infty\int_{E}f_{\varepsilon}(z,y)\,\lambda(\phi(\varepsilon\,z,x,a'),a')\,Q(\phi(\varepsilon\,z,x,a'),a',dy)\,\,dz,\\
 I_2(\varepsilon)&=\int_0^\infty\int_A g_{\varepsilon}(z,b)\,\lambda_0(db)\,dz,
 \end{align*}
 where
 \begin{align*}
 &f_{\varepsilon}(z,y):=\varepsilon\,\varphi(\varepsilon\,z,y,a')\, \exp\left(-\int_0^{\varepsilon\,z}\lambda(\phi(t,x,a'),a')\,dt-z\right),\\
 &g_{\varepsilon}(z,b):=\varphi(\varepsilon\,z,\phi(\varepsilon\,z,x,a'),b)\, \exp\left(-\int_0^{\varepsilon\,z}\lambda(\phi(t,x,a'),a')\,dt-z\right)\,\frac{1}{\lambda_0(\bar B(a,h(\varepsilon)))}\one_{\{b \in \bar B(a,h(\varepsilon))\}}.
 \end{align*}
 Exploiting the continuity properties of  $\lambda$, $Q$, $\phi$ and $f$,   
 we get
 \begin{align}\label{Sec:PDP_I2_conv}
 I_2(\varepsilon) &\limit_{\varepsilon \rightarrow 0}
 \varphi(0,x,a),
 \end{align}
 where we have used  that $\phi(0,x,b)=x$ for every $b \in A$.
 On the other hand, from estimate \eqref{Sec:PDP_estimate_varphi}, it follows that
 \begin{align*}
 |f_{\varepsilon}(z,y)| 
 &\leq \frac{M_f}{\delta}\,\,e^{-z}\,\varepsilon. 
 \end{align*}
 Therefore
 \begin{equation}\label{Sec:PDP_I1_conv}
 |I_1(\varepsilon)|\leq \frac{M_f}{\delta}\, \varepsilon \, ||\lambda||_{\infty}\int_0^\infty e^{-z}\,dz=\frac{M_f}{\delta}\, \varepsilon \, ||\lambda||_{\infty} \limit_{\varepsilon \rightarrow 0} 0.
 \end{equation}
 Collecting \eqref{Sec:PDP_I1_conv} and \eqref{Sec:PDP_I2_conv}, we conclude that
 \begin{align}
 {\mathbb E}^{x,a'}_{\nu^{\varepsilon}}[\bar F( 
 \Pi^{1,k})]
 \underset{\varepsilon \rightarrow 0}{\longrightarrow}\,
 \varphi(0,x,a).
 \end{align}
 Recalling the definitions of  $\varphi$
 and $\bar F$
 given respectively in \eqref{Sec:PDP_def_varphi} and \eqref{barF},
 we see that
 \begin{align*}
 &\varphi(0,x,a)\\
 &=(\P^{x,a}_{\nu}(T_1 >0))^{-1}\,\mathbb E^{x,a}_{\nu}\left[\one_{\{T_1 >0\}}\,
 \bar F(0,x,a,\Pi^{1,k-1} 
 )
 \right]\\
 &=\mathbb E^{x,a}_{\nu}\left[
 \bar F(0,x,a,\Pi^{1,k-1} 
 )
 \right]\\
 &=\mathbb E^{x,a}_{\nu} \bigg[ \int_{0}^{T_1} e^{-\delta t} f(\phi(t,x,a),a)\,dt
 + \sum_{n=2}^{k}\int_{T_{n-2}}^{T_{n-1}} e^{-\delta t}f(\phi(t-T_{n-1},E_{n-1},A_{n-1}), A_{n-1})\,dt\bigg]\\
 &= {\mathbb E^{x,a}_{\nu}}\left[ 
 \bar F( 
 \Pi^{1,k})
 \right],
 \end{align*}
 and this concludes the proof.
 \qed

\paragraph{Acknowledgements.}
The  author  partially benefited
 from the support of the  Italian MIUR-PRIN 2010-11 "Evolution differential problems: deterministic and stochastic approaches and their interactions".

\end{document}